\documentclass[12pt]{amsart}
\usepackage{latexsym,amsxtra,eucal}

\calclayout

\let\le\undefined
\DeclareMathSymbol{\le}{\mathrel}{AMSa}{"36}         
\let\ge\undefined
\DeclareMathSymbol{\ge}{\mathrel}{AMSa}{"3E}         
\DeclareMathSymbol{\boxt}{\mathbin}{AMSa}{"02}       

\newcommand{\Z}{{\mathbb Z}}
\newcommand{\Q}{{\mathbb Q}}
\renewcommand{\L}{{\mathbb L}}
\newcommand{\R}{{\mathbb R}}

\newcommand{\maps}{\longmapsto}
\newcommand{\rarrow}{\longrightarrow}
\newcommand{\lrarrow}{\.\relbar\joinrel\relbar\joinrel\rightarrow\.}
\newcommand{\ot}{\otimes}
\newcommand{\bt}{\mathbin{\text{\smaller$\boxt$}}}

\newcommand{\+}{\protect\nobreakdash-}
\renewcommand{\:}{\colon}
\renewcommand{\.}{\mskip .5\thinmuskip\relax}
\renewcommand{\;}{,\medspace}
\renewcommand{\d}{\partial}
\newcommand{\bu}{\bullet}
\newcommand{\dsb}{\dotsb}
\newcommand{\dsc}{\dotsc}

\DeclareMathOperator{\Hom}{Hom}
\DeclareMathOperator{\Tor}{Tor}
\DeclareMathOperator{\Ext}{Ext}
\DeclareMathOperator{\Tot}{Tot}
\DeclareMathOperator{\Br}{Bar}
\DeclareMathOperator{\Cb}{Cob}
\DeclareMathOperator{\Hoch}{Hoch}
\DeclareMathOperator{\Id}{Id}
\DeclareMathOperator{\id}{id}
\DeclareMathOperator{\cone}{cone}

\newcommand{\modl}{\text{-}\mathrm{mod}}
\newcommand{\modld}{\text{-}\mathrm{mod}^{\mathrm{dg}}}
\newcommand{\modlc}{\text{-}\mathrm{mod}^{\mathrm{cdg}}}
\newcommand{\modlq}{\text{-}\mathrm{mod}^{\mathrm{qdg}}}
\newcommand{\modrd}{\mathrm{mod}^{\mathrm{dg}}\text{-}}
\newcommand{\modrc}{\mathrm{mod}^{\mathrm{cdg}}\text{-}}
\newcommand{\modrq}{\mathrm{mod}^{\mathrm{qdg}}\text{-}}
\newcommand{\modrcfp}{\mathrm{mod}^{\mathrm{cdg}}_
                     {\mathrm{fgp}}\text{-}}
\newcommand{\modrcfg}{\mathrm{mod}^{\mathrm{cdg}}_
                     {\mathrm{fg}}\text{-}}
\newcommand{\modrqfp}{\mathrm{mod}^{\mathrm{qdg}}_
                     {\mathrm{fgp}}\text{-}}
\newcommand{\modrcfl}{\mathrm{mod}^{\mathrm{cdg}}_
                     {\mathrm{fl}}\text{-}}
\newcommand{\modrcffd}{\mathrm{mod}^{\mathrm{cdg}}_
                     {\mathrm{ffd}}\text{-}}
\newcommand{\modrdfl}{\mathrm{mod}^{\mathrm{dg}}_
                     {\mathrm{fl}}\text{-}}
\newcommand{\modrdffd}{\mathrm{mod}^{\mathrm{dg}}_
                     {\mathrm{ffd}}\text{-}}
\newcommand{\comodlc}{\text{-}\mathrm{comod}^{\mathrm{cdg}}}

\newcommand{\prj}{{\mathrm{prj}}}
\newcommand{\inj}{{\mathrm{inj}}}
\newcommand{\fl}{{\mathrm{fl}}}
\newcommand{\co}{{\mathrm{co}}}
\newcommand{\ctr}{{\mathrm{ctr}}}
\newcommand{\abs}{{\mathrm{abs}}}
\newcommand{\cmp}{{\mathrm{cmp}}}
\newcommand{\fpd}{{\mathrm{fpd}}}
\newcommand{\fid}{{\mathrm{fid}}}
\newcommand{\fp}{{\mathrm{fgp}}}
\newcommand{\fg}{{\mathrm{fg}}}
\renewcommand{\b}{{\mathrm{b}}}
\newcommand{\Sing}{{\mathit{Sing}}}
\newcommand{\kfl}{{k\text{-}\mathrm{fl}}}
\renewcommand{\ss}{{\mathrm{ss}}}
\newcommand{\gr}{\mathrm{gr}}

\newcommand{\DD}{{\mathrm D}}
\newcommand{\CC}{{\mathcal C}}

\newcommand{\Pre}{\mathrm{Pre}}
\newcommand{\Comex}{\mathrm{Com}_{\mathrm{ex}}}
\newcommand{\Comab}{\mathrm{Com}_{\mathrm{ab}}}
\newcommand{\Ac}{\mathrm{Ac}}

\newcommand{\op}{{\mathrm{op}}}

\newcommand{\Section}[1]{\bigskip
\section{#1}
\medskip}

\theoremstyle{plain}
\newtheorem*{thm}{Theorem}
\newtheorem*{thmA}{Theorem A}
\newtheorem*{thmB}{Theorem B}
\newtheorem*{thmC}{Theorem C}
\newtheorem*{thmD}{Theorem D}
\newtheorem*{prop}{Proposition}
\newtheorem*{propA}{Proposition A}
\newtheorem*{propB}{Proposition B}
\newtheorem*{propC}{Proposition C}

\newtheorem*{lemA}{Lemma A}
\newtheorem*{lemB}{Lemma B}
\newtheorem*{lemC}{Lemma C}
\newtheorem*{lemD}{Lemma D}
\newtheorem*{cor}{Corollary}
\newtheorem*{corA}{Corollary A}
\newtheorem*{corB}{Corollary B}
\newtheorem*{corC}{Corollary C}
\newtheorem*{corE}{Corollary E}
\theoremstyle{definition}
\newtheorem*{rem}{Remark}

\begin{document}

\title{Hochschild (co)homology of the second kind I}
\author{Alexander Polishchuk \ and \ Leonid Positselski}

\address{Department of Mathematics, University of Oregon, 
Eugene, OR 97403, USA} 
\email{apolish@uoregon.edu}

\address{Sector of Algebra and Number Theory, Institute for
Information Transmission Problems, Bolshoy Karetny per.~19 str.~1,
Moscow 127994, Russia}
\email{posic@mccme.ru}

\maketitle

\tableofcontents

\section*{Introduction}
\medskip

 CDG\+algebras (where ``C'' stands for ``curved'') were introduced
in connection with nonhomogeneous Koszul duality in~\cite{Pcurv}.
 Several years earlier, (what we would now call)
$A_\infty$\+algebras with curvature were considered in~\cite{GJ}
as natural generalizations of the conventional $A_\infty$\+algebras.
 In fact, \cite{GJ}~appears to be the first paper where the Hochschild
(and even cyclic) homology of curved algebras was discussed.

 Recently, the interest to these algebras was rekindled by their
connection with the categories of matrix
factorizations~\cite{Seg,Dyck,PV,CT,Tu}.
 In these studies, beginnings of the theory of Hochschild
(co)homology for CDG\+algebras have emerged.
 The aim of the present paper is to work out the foundations of
the theory on the basis of the general formalism of \emph{derived
categories of the second kind} as developed in the second author's
paper~\cite{Pkoszul}.
 The terminology, and the notion of a \emph{differential derived
functor of the second kind}, which is relevant here, go
back to the classical paper~\cite{HMS}.

 The subtle but crucial difference between the differential derived
functors of the first and the second kind lies in the way one
constructs the totalizations of bicomplexes: one can take either
direct sums or direct products along the diagonals.
 The construction of the differential $\Tor$ and $\Ext$ of
the first kind, which looks generally more natural at the first
glance, leads to trivial functors in the case of a CDG\+algebra
with nonzero curvature over a field.
 So does the (familiar) definition of Hochschild (co)homology of
the first kind.

 On the other hand, with a CDG\+algebra $B$ one can associate
the DG\+category $C$ of right CDG\+modules over $B$, projective
and finitely generated as graded $B$\+modules.
 For the DG\+category $C$, the Hochschild (co)homology of the first
kind makes perfect sense.
 The main problem that we address in this paper is the problem of
comparison between the Hochschild (co)homology of the first kind of
the DG\+category $C$ and the Hochschild (co)homology of the second kind
of the original CDG\+algebra $B$ (defined using the differential
$\Tor$/$\Ext$ of the second kind).

 We proceed in two steps: first, compare the Hochschild (co)homology
of the second kind for $B$ and $C$, and then deal with the two kinds
of Hochschild (co)homology of~$C$.
 The first step is relatively easy: our construction of an isomorphism
works, at least, for any CDG\+algebra $B$ over a field~$k$
(see Section~\ref{dg-of-cdg-subsect}).
 However, a trivial counterexample shows that the two kinds of
Hochschild (co)homology of $C$ are \emph{not} isomorphic in general
(see Section~\ref{counterexample}).
 There are natural maps between the two kinds of Hochschild
(co)homology, though.
 A sufficient condition for these maps to be isomorphisms is formulated
in terms of the derived categories of the second kind of CDG\+bimodules
over~$B$.
 In the maximal generality that we have been able to attain, this is
a kind of ``resolution of the diagonal'' condition for
the CDG\+bimodule $B$ over~$B$
(see Theorems~\ref{comparison-dg-of-cdg}.C\+-D and
Corollaries~\ref{cdg-koszul}.B, \ref{noetherian-cdg-rings}.B,
and~\ref{matrix-factorizations}).

 Let us say a few more words about the first step.
 There is no obvious map between the Hochschild complexes of $B$
and $C$, so one cannot directly compare their cohomology.
 Instead, we construct a third complex (both in the homological
and the cohomological versions) endowed with natural maps
from/to these two complexes, and show that these maps are
quasi-isomorphisms.
 To obtain the intermediate complex, we embed both $B$ and $C$
into a certain larger differential category.

 The idea of these embeddings goes back to A.~Schwarz's
work~\cite{Sch}.
 The starting observation is that a CDG\+algebra is not
a CDG\+module over itself in any natural way (even though it is
naturally a CDG\+bimod\-ule over itself).
 It was suggested in~\cite{Sch}, however, that one can relax
the conditions on differential modules over CDG\+algebras
(called ``$Q$\+algebras'' in~\cite{Sch}) thereby making the modules
carry their own curvature endomorphisms.
 In recognition of A.~Schwarz's vision, we partly borrow his
terminology by calling such modules \emph{QDG\+modules}.

 Any CDG\+algebra is naturally both a left and a right QDG\+module
over itself.
 While CDG\+modules form a DG\+category, QDG\+modules form
a \emph{CDG\+category}.
 Both a CDG\+algebra $B$ (considered as a CDG\+category with
a single object) and the DG\+category $C$ of CDG\+modules over
it embed naturally into the CDG\+category $D$ of QDG\+modules
over $B$, so the Hochschild complex of $D$ provides
an intermediate object for comparison between the Hochschild
complexes of $B$ and~$C$.

 Now let us turn to the second step.
 The (conventional) derived category of DG\+modules over a DG\+algebra
is defined as the localization of the homotopy category of DG\+modules
by the class of quasi-isomorphisms, or equivalently, by
the thick subcategory of acyclic DG-modules.
 This does not make sense for CDG\+modules, since their differentials
have nonzero squares, so their cohomology cannot be defined.
 Indeed, the subcategory of acyclic DG\+modules is not even invariant
under CDG\+isomorphisms between
DG\+algebras~\cite[Examples~9.4]{Pkoszul}.

 The definition of the \emph{derived categories of the second kind},
various species of which are called the \emph{coderived},
the \emph{contraderived}, the \emph{absolute derived}, and
the \emph{complete derived categories}, for DG- and CDG-modules are
not based on any notion of cohomology of a differential module.
 Rather, the classes of \emph{coacyclic}, \emph{contraacyclic},
\emph{absolutely acyclic}, and \emph{completely acyclic} CDG\+modules
are built up starting from short exact sequences of CDG\+modules
(with closed morphisms between them).

 For reasons related to the behavior of tensor products with
respect to infinite direct sums and products of vector spaces,
the derived categories and functors of the second kind work better
for \emph{coalgebras} than for algebras, even though one is forced
to use them for algebras if one is interested in curved algebras
and modules.
 (For derived categories and functors of the first kind, it is
the other way.)
 That is why one has to impose additional conditions like finiteness
of homological dimension, Noetherianness, etc.,\ on the underlying
graded algebras of one's CDG\+algebras in order to make
the derived categories of the second kind well-behaved and
the relation between them and the derived functors of the second kind
working properly.
 We did our best to make such additional conditions as weak as possible
in this paper, but the price of generality is technical complexity.

 Unlike the $\Tor$ and $\Ext$, the Hochschild (co)homology is
essentially an invariant of a pair (a field or commutative ring,
an algebra over it).
 It is \emph{not} preserved when the ground field or ring is changed.
 In this paper, we always work over an arbitrary commutative ring~$k$,
or a commutative ring of finite homological dimension, as needed.
 The only exceptions are some examples depending on the Koszul duality
results from~\cite{Pkoszul}, which are established only over a field.
 Working over a commutative ring involves all kinds of $k$\+flatness
or $k$\+projectivity conditions that need to be imposed on
the algebras and modules, both in order to define the Hochschild
(co)homology and to compute various (co)homology theories in terms
of standard complexes.

 Recent studies of the categories of matrix factorizations and of
the associated CDG\+algebras showed the importance of developing
the relevant homological algebra using only $\Z/2$\+grading (as
opposed to the conventional $\Z$\+grading).
 In this paper we work with CDG\+algebras and CDG\+categories
graded by an arbitrary abelian group $\Gamma$ endowed with
some additional data that is needed to define $\Gamma$\+graded
complexes and perform operations with them.
 The behavior of our (co)homology theories with respect to
a replacement of the grading group $\Gamma$ is discussed
in detail (see Section~\ref{change-grading-group}).

 We exhibit several classes of examples of DG\+algebras and
DG\+categories for which the two kinds of $\Tor$, $\Ext$, and
Hochschild (co)homology coincide.
 These examples roughly correspond to the classes of DG\+algebras
for which the derived categories of the first and second kind
are known to coincide~\cite[Section~9.4]{Pkoszul}.
 In particular, one of these classes is that of the DG\+categories
that are cofibrant with respect to G.~Tabuada's model category
structure (see Section~\ref{cofibrant-subsect}).

 Examples of CDG\+algebras $B$ such that the two kinds of $\Tor$
and $\Ext$ for the corresponding DG\+category $C$ of CDG\+modules
over $B$, finitely generated and projective as graded $B$\+modules,
are known to coincide are fewer; and examples when we can show that
the two kinds of Hochschild (co)homology for this DG\+category $C$
coincide are fewer still.
 Among the former are all the CDG\+rings $B$ whose underlying
graded rings are Noetherian of finite homological dimension
(see Section~\ref{noetherian-cdg-rings}).
 In the latter class we have some CDG\+algebras over fields
admitting Koszul filtrations of finite homological dimension
(see Section~\ref{cdg-koszul}), curved commutative local algebras
describing germs of isolated hypersurface singularities (due to
the results of~\cite{Dyck}), and curved commutative smooth algebras
over perfect fields with the curvature function having no other
critical values but zero (due to the recent results of~\cite{LP};
see Section~\ref{matrix-factorizations}).

 Our discussion of the Hochschild (co)homology of the DG\+categories
of matrix factorizations is finished in Section~\ref{direct-sum},
where we show that the Hochschild (co)homology of the second kind
of the DG\+category of matrix factorizations over a smooth affine
variety over an algebraically closed field of characteristic zero
is isomorphic to the direct sum of the Hochschild (co)homology of
the first kind of the similar DG\+categories corresponding to all
the critical values of the potential.

 We are grateful to Anton Kapustin, Ed Segal, Daniel Pomerleano,
Kevin Lin, and Junwu Tu for helpful conversations.
 A.~P. is partially supported by the NSF grant DMS-1001364.
 L.~P. is partially supported by a grant from P.~Deligne 2004
Balzan prize and an RFBR grant.

\Section{CDG-Categories and QDG-Functors}

 This section is written in the language of CDG\+categories.
 Expositions in the generality of CDG\+rings, which might be
somewhat more accessible to an inexperienced reader, can be
found in~\cite{Pcurv,Pkoszul,Sch}.
 For a discussion of DG\+categories, we refer to~\cite{Kel},
\cite{Toen}, and~\cite[Section~1.2]{Pkoszul}.

\subsection{Grading group} \label{grading-group}
 Let $\Gamma$ be an abelian group endowed with a symmetric bilinear
form $\sigma\:\Gamma\times\Gamma\rarrow \Z/2$ and a fixed element
$\boldsymbol{1}\in\Gamma$ such that $\sigma(\boldsymbol{1},
\boldsymbol{1}) = 1\bmod 2$.
 We will use $\Gamma$ as the group of values for the gradings of
our complexes.
 The differentials will raise the degree by~$\boldsymbol{1}$, and
signs like $(-1)^{\sigma(a,b)}$ will appear in the sign rules.

 For example, in the simplest cases one may have $\Gamma=\Z$, \
$\boldsymbol{1}=1$, and $\sigma(a,b)=ab\bmod 2$ for $a$, $b\in\Gamma$,
or, alternatively, $\Gamma=\Z/2$, \ $\boldsymbol{1}=1\bmod 2$, and
$\sigma(a,b)=ab$.
 One can also take $\Gamma$ to be any additive subgroup of $\Q$,
containing $\Z$ and consisting of fractions with odd denominators,
$\boldsymbol{1}=1$, and $\sigma(a,b)=ab\bmod 2$.
 Of course, it is also possible that $\Gamma=\Z^d$ for any finite
or infinite~$d$, etc.
 When working over a commutative ring~$k$ containing
the field~$\mathbb F_2$, we will not need the form~$\sigma$,
and so $\Gamma=\Q$ or $\Gamma=0$ become admissible choices as well.

 From now on, we will assume a grading group data $(\Gamma,\sigma,
\boldsymbol{1})$ to be fixed.
 When appropriate, we will identify the integers with their images
under the natural map $\Z\rarrow \Gamma$ sending $1$
to~$\boldsymbol{1}$ without presuming this map to be injective,
and denote $\sigma(a,b)$ simply by $ab$ for $a$, $b\in\Gamma$.
 So we will write simply $1$ instead of~$\boldsymbol{1}$, etc.
 This map $\Z\rarrow\Gamma$ will be also used when constructing
the total complexes of polycomplexes some of whose gradings are
indexed by the integers and the other ones by elements of
the group~$\Gamma$.
 Conversely, to any $a\in\Gamma$ one assigns the class
$\sigma(\boldsymbol{1},a)\in\Z/2$, which we will denote simply
by~$a$ in the appropriate contexts.

\subsection{CDG-categories}  \label{cdg-categories-subsect}
 A \emph{CDG\+category} $C$ is a category whose sets of morphisms
$\Hom_C(X,Y)$ are $\Gamma$\+graded abelian groups (i.~e., $C$ is
a $\Gamma$\+graded category) endowed with homogeneous endomorphisms
$d\:\Hom_C(X,Y)\rarrow\Hom_C(X,Y)$ of degree~$1$ and fixed elements
$h_X\in\Hom_C(X,X)$ of degree~$2$ for all objects $X$, $Y\in C$.
 The endomorphisms $d$ are called the \emph{differentials} and
the elements $h_X$ are called the \emph{curvature elements}.
 The following equations have to be satisfied: $d(fg)=d(f)g+
(-1)^{|f|}fd(g)$ for any composable homogeneous morphisms $f$
and $g$ in $C$ of the degrees $|f|$ and $|g|\in\Gamma$, \
$d^2(f)=h_Yf-fh_X$ for any morphism $f\:X\rarrow Y$ in~$C$,
and $d(h_X)=0$ for any object $X\in C$.

 The simplest example of a CDG\+category is the category $\Pre(A)$ of
\emph{precomplexes} over an additive category~$A$.
 The objects of $\Pre(A)$ are $\Gamma$\+graded objects $X$ in $A$
endowed with an endomorphism $d_X\:X\rarrow X$ of degree~$1$.
 The $\Gamma$\+graded abelian group of morphisms $\Hom_{\Pre(A)}(X,Y)$
is the group of homogeneous morphisms $X\rarrow Y$ of
$\Gamma$\+graded objects.
 The differentials $d\:\Hom(X,Y)\rarrow\Hom(X,Y)$ are given
by the rule $d(f)=d_Y f-(-1)^{|f|}fd_X$, and the curvature
elements are $h_X=d_X^2$.
 In particular, when $A=Ab$ is the category of abelian groups,
we obtain the CDG\+category of precomplexes of abelian groups
$\Pre(Ab)$.

 A CDG\+category with a single object is another name for
a \emph{CDG\+ring}.
 A CDG\+ring $(B,d,h)$ is a $\Gamma$\+graded ring $B$ endowed with
an odd derivation~$d$ of degree~$1$ and a curvature element
$h\in B^2$ such that $d^2(b)=[h,b]$ for any $b\in B$ and $d(h)=0$.

 An \emph{isomorphism} between objects $X$ and $Y$ of a CDG\+category
$C$ is (an element of) a pair of morphisms $i\:X\rarrow Y$ and
$j\:Y\rarrow X$ of degree~$0$ such that $ji=\id_X$, \ $ij=\id_Y$,
and $d(i)=0=d(j)$; any one of the latter two equations implies
the other one.
 It also follows that $jh_Yi=h_X$.

 Let $X$ be an object of a CDG\+category $C$ and $\tau\in\Hom_C(X,X)$
be its homogeneous endomorphism of degree~$1$.
 An object $Y\in C$ is called the \emph{twist} of an object $X$ with
an endomorphism~$\tau$ (the notation: $Y=X(\tau)$) if homogeneous
morphisms $i\:X\rarrow Y$ and $j\:Y\rarrow X$ of degree~$0$ are
given such that $ji=\id_X$, \ $ij=\id_Y$, and $jd(i)=\tau$.
 In this case one has $jh_Yi=h_X+d\tau+\tau^2$.
 For any object $X\in C$ and an element $n\in\Gamma$, an object
$Y\in C$ is called the \emph{shift} of $X$ with the grading~$n$
(the notation: $Y=X[n]$) if homogeneous morphisms $i\:X\rarrow Y$
and $j\:Y\rarrow X$ of the degrees $n$ and $-n$, respectively,
are given such that $ji=\id_X$, \ $ij=\id_Y$, and $d(i)=0=d(j)$.
 In this case one has $jh_Yi=h_X$.

 An object $X\in C$ is called the \emph{direct sum} of a family of
objects $X_\alpha\in C$ if homogeneous morphisms $i_\alpha\:X_\alpha
\rarrow X$ of degree~$0$ are given such that the induced map
$\Hom_C(X,Y)\rarrow\prod_\alpha\Hom_C(X_\alpha,Y)$ is an isomorphism
of $\Gamma$\+graded abelian groups for any object $Y\in C$, and
$di_\alpha=0$.
 In this case one has $h_X i_\alpha=i_\alpha h_{X_\alpha}$, so
the endomorphism $h_X$ corresponds to the family of morphisms
$i_\alpha h_{X_\alpha}$ under the above isomorphism for $Y=X$.
 The (\emph{direct}) \emph{product} of a family of object is defined
in the dual way.
 An object $X$ is the direct sum of a finite family of objects
$X_\alpha\in C$ if and only if it is their direct product.
 Of course, the notions of a shift and a direct sum/product of objects
make sense in a (nondifferential) $\Gamma$\+graded category, too;
one just drops the conditions involving $d$ and~$h$.

 Twists, shifts, direct sums, and products of objects
of a CDG\+category are unique up to a unique isomorphism
whenever they exist.

 A \emph{DG\+category} is a CDG\+category in which all the curvature
elements are zero.

 The \emph{opposite CDG\+category} to a CDG\+category $C$ is
constructed as follows.
 The class of objects of $C^\op$ coincides with the class of objects
of~$C$.
 For any objects $X$, $Y\in C$ the $\Gamma$\+graded abelian group
$\Hom_{C^\op}(X^\op,Y^\op)$ is identified with $\Hom_C(Y,X)$, and
the differential $d^\op$ on this group coincides with~$d$.
 The composition of morphisms in $C^\op$ differs from that in $C$
by the sign rule, $f^\op g^\op = (-1)^{|f||g|}(gf)^\op$.
 Finally, the curvature elements in $C^\op$ are $h_{X^\op}=-h_X$.
 In particular, this defines the CDG\+ring
$B^\op=(B^\op,d^\op,h^\op)$ opposite to a CDG\+ring $B=(B,d,h)$.

 Now let $k$ be a commutative ring.
 A \emph{k\+linear} CDG\+category is a CDG\+category whose
$\Gamma$\+graded abelian groups of morphisms are endowed with
$\Gamma$\+graded $k$\+module structures so that the compositions
are $k$\+bilinear and the differentials are $k$\+linear.
 The \emph{tensor product} $C\ot_kD$ of two $k$\+linear CDG\+categories
$C$ and $D$ is constructed as follows.
 The objects of $C\ot_kD$ are pairs $(X',X'')$ of objects $X'\in C$
and $X''\in D$.
 The $\Gamma$\+graded $k$\+module of morphisms
$\Hom_{C\ot_kD}((X',X''),(Y',Y''))$ is the tensor product
$\Hom_C(X',Y')\ot_k\Hom_D(X'',Y'')$; the differential~$d$ on
this module is defined by the formula $d(f'\ot f'')=d(f')\ot f''
+ (-1)^{|f'|}f'\ot d(f'')$.
 The curvature elements are $h_{(X',X'')}=h_{X'}\ot\id_{X''}+
\id_{X'}\ot h_{X''}$.

\subsection{QDG-functors}
 Let $C$ and $D$ be CDG\+categories.
 A \emph{covariant CDG\+functor} $F\:C\rarrow D$ is a homogeneous
additive functor between the $\Gamma$\+graded categories $C$ and $D$,
endowed with fixed elements $a_X\in \Hom_D(F(X),F(X))$ of degree~$1$
for all objects $X\in C$ such that $F(df) = dF(f) + a_Y F(f) -
(-1)^{|f|} F(f) a_X$ for any morphism $f\:X\rarrow Y$ in $C$ and
$F(h_X) = h_{F(X)} + da_X + a_X^2$ for any object~$X$.
 A contravariant CDG\+functor $C\rarrow D$ is defined as
a covariant CDG\+functor $C^\op\rarrow D$, or equivalently,
a covariant CDG\+functor $C\rarrow D^\op$.
 The \emph{opposite} CDG\+functor $F^\op\:C^\op\rarrow D^\op$ to
a covariant CDG\+functor $F\:C\rarrow D$ is defined
by the rule $(F,a)^\op=(F^\op,-a)$.

 (Covariant or contravariant) CDG\+functors $C\rarrow D$ are
objects of the \emph{DG\+category of CDG\+functors}.
 The $\Gamma$\+graded abelian group of morphisms between covariant
CDG\+functors $F$ and $G$ is the $\Gamma$\+graded group of
homogeneous morphisms, with the sign rule, between $F$ and $G$
considered as functors between $\Gamma$\+graded categories.
 More precisely, a morphism $f\:F\rarrow G$ of degree $n\in\Gamma$
is a collection of morphisms $f_X\:F(X)\rarrow G(X)$ of degree~$n$
in $D$ for all objects $X\in C$ such that $f_Y F(g)=(-1)^{n|g|}
G(g) f_X$ for any morphism $g\:X\rarrow Y$ in~$C$. 
 The differential~$d$ on the $\Gamma$\+graded group $\Hom(F,G)$ of
morphisms between CDG-functors $F=(F,a)$ and $G=(G,b)$ is defined
by the rule $(df)_X = d(f_X) + b_X f_X - (-1)^{|f|} f_X a_X$.

 A (covariant or contravariant) \emph{QDG\+functor} $F$ between
CDG\+categories $C$ and $D$ is the same set of data as
a CDG\+functor satisfying the same equations, except for
the equation connecting $F(h_X)$ with $h_{F(X)}$, which is omitted.
 QDG\+functors $C\rarrow D$ are objects of the \emph{CDG\+category
of QDG\+functors}.
 The $\Gamma$\+graded abelian group of morphisms between
QDG\+functors and the differential on it are defined exactly in
the same way as in the CDG\+functor case.
 The curvature element of a QDG\+functor $F\:C\rarrow D$ is
the endomorphism $h_F\:F\rarrow F$ of degree~$2$ defined by
the formula $(h_F)_X = h_{F(X)} + da_X + a_X^2 - F(h_X)$
for all $X\in C$.

 The composition of QDG\+functors $(F,a)\:C\rarrow D$ and
$(G,b)\:D\rarrow E$ is the QDG\+functor $(G\circ F\;c)$, where
$c_X=G(a_X)+b_{F(X)}$ for any object $X\in C$.
 A CDG\+functor or QDG\+functor $F=(F,a)\:C\rarrow D$ is said to be
\emph{strict} if $a_X=0$ for all objects $X\in C$.
 The identity CDG\+functor $\Id_C$ of a CDG\+category $C$ is
the strict CDG\+functor $(\Id_C,0)$.
 The composition of strict QDG\+functors is a strict QDG\+functor,
and the composition of (strict) CDG\+functors is a (strict)
CDG\+functor.

 Two CDG\+functors $F\:C\rarrow D$ and $G\:D\rarrow C$ between
CDG\+categories $C$ and $D$ are called mutually inverse
\emph{equivalences} of CDG\+categories if they are equivalences
of the $\Gamma$\+graded categories such that the adjunction
isomorphisms $i\:GF\rarrow\Id_C$ and $j\:FG\rarrow\Id_D$ are
closed morphisms of CDG\+functors, i.~e., $d(i)=0=d(j)$
(any one of the two equations implies the other one).
 A CDG\+functor $F\:C\rarrow D$ is an equivalence if and only if it
is fully faithful as a functor between $\Gamma$\+graded categories
and any object $Y\in D$ is a twist of an object $F(X)$ for some
$X\in C$.

 An equivalence $(F,G)$ between CDG\+categories $C$ and $D$ is called
a \emph{strict equivalence} if the CDG\+functors $F$ and $G$ are
strict.
 A strict CDG\+functor $F\:C\rarrow D$ is a strict equivalence if and
only if it is fully faithful as a functor between $\Gamma$\+graded
categories and any object $Y\in D$ is isomorphic to an object $F(X)$
for some $X\in C$.

 A strict CDG\+functor between DG\+categories is called
a \emph{DG\+functor}.
 An \emph{equivalence} of DG\+categories is their strict equivalence
as CDG\+categories.

 If all objects of the category $D$ admit twists with all of their
endomorphisms of degree~$1$, then the embedding of the DG\+category
of strict CDG\+functors $C\rarrow D$ into the DG\+category of all
CDG\+functors is an equivalence of DG\+categories, and the embedding
of the CDG\+category of strict QDG\+functors $C\rarrow D$ into
the CDG\+category of all QDG\+functors is a strict equivalence of
CDG\+categories.

 A QDG\+functor between $k$\+linear CDG\+categories is \emph{k\+linear}
if its action on the $\Gamma$\+graded $k$\+modules of morphisms in
the CDG\+categories is $k$\+linear.
 Given three $k$\+linear CDG\+categories $C$, $D$, $E$, the functor
of composition of $k$\+linear QDG\+functors $C\rarrow D$ and
$D\rarrow E$ is a strict $k$\+linear CDG\+functor on the tensor
product of the $k$\+linear CDG\+categories of QDG\+functors.
 The composition (on either side) with a fixed CDG\+functor is
a strict CDG\+functor between the CDG\+categories of QDG\+functors,
and the composition with a fixed QDG\+functor is a strict
QDG\+functor between such CDG\+categories.

 Given two $k$\+linear QDG\+functors $F'=(F',a')\:C'\rarrow D'$ and
$F''=(F'',a'')\:C''\allowbreak\rarrow D''$, their tensor product
$(F'\ot F''\;a)\:C'\ot_k C''\rarrow D'\ot_k D''$ is defined by
the rule $(F'\ot F'')(X',X'')=(F'(X'),F''(X''))$ on the objects,
$(F'\ot F'')(f'\ot f'')=F'(f')\ot F''(f'')$ on the morphisms, and
$a_{X'\ot X''}=a_{X'}\ot\id_{X''}+\id_{X'}\ot a_{X''}$.
 The tensor product of strict QDG\+functors is a strict QDG\+functor,
and the tensor product of (strict) CDG\+functors is a (strict)
CDG\+functor.

\subsection{QDG\+modules}  \label{qdg-modules-subsect}
 A \emph{left QDG\+module} over a small CDG\+category $C$ is
a strict covariant QDG\+functor $C\rarrow\Pre(Ab)$.
 Analogously, a right QDG\+module over $C$ is a strict contravariant
QDG\+functor $C^\op\rarrow\Pre(Ab)$.
 (Left or right) CDG\+modules over a CDG\+category $C$ are 
similarly defined in terms of strict CDG\+functors with values in
the CDG\+category $\Pre(Ab)$.
 The CDG\+categories of left and right QDG\+modules over $C$ are
denoted by $C\modlq$ and $\modrq C$; the DG\+categories of left
and right CDG\+modules over $C$ are denoted by $C\modlc$ and
$\modrc C$.
 Since the CDG\+category $\Pre(Ab)$ admits arbitrary twists,
one obtains (strictly) equivalent (C)DG\+categories by considering
not necessarily strict QDG\+ or CDG\+functors.

 Given a CDG\+ring or CDG\+category $C$, we will denote by $C^\#$
the underlying $\Gamma$\+graded ring or category.
 For a QDG\+module $M$ over $C$, we similarly denote by $M^\#$
the underlying $\Gamma$\+graded $C^\#$\+module (i.~e., homogeneous
additive functor from $C^\#$ to the $\Gamma$\+graded category of
$\Gamma$\+graded abelian groups) of~$M$.

 If $k$ is a commutative ring and $C$ is a $k$\+linear
CDG\+category, then any QDG\+functor $C\rarrow\Pre(Ab)$ can be
lifted to a $k$\+linear QDG\+functor $C\rarrow\Pre(k\modl)$ in
a unique way, where $k\modl$ denotes the abelian category of
$k$\+modules.
 So the CDG\+category $C\modlq$ can be also described as
the CDG\+category of (strict) $k$\+linear QDG\+functors
$C\rarrow\Pre(k\modl)$.
 Notice that another notation for the CDG\+category $\Pre(k\modl)$
is $k\modlq$, where $k$ is considered as a CDG\+ring concentrated
in degree~$0$ with the trivial differential and curvature, while
$k\modlc$ is a notation for the DG\+category of complexes of
$k$\+modules.

 Let $C$ be a small $k$\+linear CDG\+category, $N$ be a right
QDG\+module over $C$, and $M$ be a left QDG\+module.
 The tensor product $N^\#\ot_{C^\#} M^\#$ is a $\Gamma$\+graded
$k$\+module defined as the quotient module of the direct sum of
$N(X)\ot_k M(X)$ over all objects $X\in C$ by the sum of the images
of the maps $N(Y)\ot_k M(X)\rarrow N(X)\ot_k M(X)\oplus
N(Y)\ot_k M(Y)$ over all homogeneous morphisms $X\rarrow Y$ in~$C$.
 There is a natural differential on $N^\#\ot_{C^\#}M^\#$ defined
by the usual formula $d(n\ot m)=d(n)\ot m + (-1)^{|n|} n\ot d(m)$.
 The precomplex of $k$\+modules so obtained is denoted by
$N\ot_CM$.
 The tensor product over~$C$ is a strict CDG\+functor
$$
 \ot_C\:\modrq C\times C\modlq\lrarrow k\modlq,
$$
and its restriction to the DG\+subcategories of CDG\+modules
is a DG\+functor
$$
 \ot_C\:\modrc C\times C\modlc\lrarrow k\modlc.
$$

 A QDG\+functor between CDG\+categories $F\:C\rarrow D$ induces
a strict QDG\+func\-tor of inverse image (restriction of scalars)
$F^*\:D\modlq\rarrow C\modlq$.
 Here we use the natural strict equivalence between
the CDG\+categories of arbitrary and strict QDG\+functors
$C\rarrow\Pre(Ab)$.
 When $F$ is a CDG\+functor, the functor $F^*$ is a strict
CDG\+functor, and it restricts to a DG\+functor
$D\modlc\rarrow C\modlc$.
 For any right QDG\+module $N$ and left QDG\+module $M$ over
a $k$\+linear CDG\+category $D$ and a $k$\+linear CDG\+functor
$F\:C\rarrow D$ there is a natural map of precomplexes
of $k$\+modules $F^*(N)\ot_C F^*(M)\rarrow N\ot_D M$,
commuting with the differentials.

 For any CDG\+category $B$ there is a natural strict CDG\+functor
$B\rarrow\modrq B$ assigning to an object $X\in B$ the right
QDG\+module $R_X\:Y\longmapsto\Hom_B(Y,X)$ over~$B$.
 Here the differential on $R_X(Y)$ coincides with the differential
on $\Hom_B(Y,X)$.

 A CDG\+module over a DG\+category is called a \emph{DG\+module}.
 The DG\+categories of left and right DG\+modules over a small
DG\+category $C$ are denoted by $C\modld$ and $\modrd C$.
 In particular, $k\modld$ is yet another notation for
the DG\+category of complexes of $k$\+modules for a commutative
ring~$k$.
 If $C$ is a $k$\+linear DG\+category, then the objects of
$C\modld$ can be viewed as DG\+functors $C\rarrow k\modld$,
and the objects of $\modrd C$ can be viewed as DG\+functors
$C^\op\rarrow k\modld$.

 Given left QDG\+modules $M'$ and $M''$ over $k$\+linear
CDG\+categories $B'$ and $B''$, their tensor product $M'\ot_k M''$
is the QDG\+module over $B'\ot_k B''$ defined as the composition
of the tensor product of strict QDG\+functors
$M'\ot M''\:B'\ot_k B''\rarrow \Pre(k\modl)\ot_k\Pre(k\modl)$
with the strict CDG\+functor of tensor product of precomplexes
$\ot_k\:\Pre(k\modl)\ot_k\Pre(k\modl)\rarrow\Pre(k\modl)$.
 The latter functor assigns to two precomplexes of $k$\+modules
their tensor product as $\Gamma$\+graded $k$\+modules, endowed
with the differential defined by the usual formula.
 The tensor product of CDG\+modules is a CDG\+module.

\subsection{Pseudo-equivalences}  \label{pseudo-equi-subsect}
 Let us call a homogeneous additive functor $F^\#\:C^\#\rarrow D^\#$
between $\Gamma$\+graded additive categories $C^\#$ and $D^\#$
a \emph{pseudo-equivalence} if $F$ is fully faithful and any object
$Y\in D^\#$ can be obtained from objects $F(X)$, \ $X\in C^\#$, using
the operations of finite direct sum, shift, and passage to a direct
summand. {\hfuzz=3pt\par}
 A CDG\+functor between CDG\+categories $F\:C\rarrow D$ is called
a \emph{pseudo-equivalence} if it is fully faithful as a functor
between the $\Gamma$\+graded categories and any object $Y\in D$
can be obtained from objects $F(X)$, \ $X\in C$, using the operations
of finite direct sum, shift, twist, and passage to a direct summand.

 The category of (left or right) $\Gamma$\+graded modules over
a small $\Gamma$\+graded category $C^\#$ is abelian.
 Let us call a right $\Gamma$\+graded module $N$ over $C^\#$
(\emph{finitely generated}) \emph{free} if it is a (finite) direct
sum of representable modules $R_X$, where $X\in C^\#$.
 A $\Gamma$\+graded module $P$ over $C^\#$ is a projective object
in the abelian category of $\Gamma$\+graded modules if and only if
it is a direct summand of a free $\Gamma$\+graded module.
 A $\Gamma$\+graded module $P$ is a compact projective object
(i.~e., a projective object representing a covariant functor
preserving infinite direct sums on the category of modules) if and
only if it is a direct summand of a finitely generated free
$\Gamma$\+graded module.
 In this case, a $\Gamma$\+graded modules $P$ is said to be
\emph{finitely generated projective}.

 Given a CDG\+category $B$, denote by $\modrcfp B$ and $\modrqfp B$
the DG\+category of right CDG\+modules and the CDG\+category of right
QDG\+modules over $B$, respectively, which are finitely generated
projective as $\Gamma$\+graded modules.
 The representable QDG\+modules $R_X$ are obviously objects of
$\modrqfp B$, so there is a strict CDG\+functor $R\:B\rarrow
\modrqfp B$.
 There is also the strict CDG\+functor of tautological embedding
$I\:\modrcfp B\rarrow\modrqfp B$.

\begin{lemA}
 The CDG\+functors $R$ and $I$ are pseudo-equivalences. 
\end{lemA}

\begin{proof}
 First of all notice that any two objects of a CDG\+category $C$
that are isomorphic in the $\Gamma$\+graded category $C^\#$
are each other's twists.
 In particular, so are any two QDG\+modules over a CDG\+category
$B$ that are isomorphic as $\Gamma$\+graded $B^\#$\+modules.
 Hence in order to prove that $R$ is a pseudo-equivalence, it suffices
to show that any (finitely generated) projective $\Gamma$\+graded
right $B^\#$\+module $P$ admits a QDG\+module structure.
 Indeed, if there is a $\Gamma$\+graded right $B^\#$\+module $Q$
such that the $\Gamma$\+graded module $P\oplus Q$ admits
a differential~$d$ making it a QDG\+module, and $\iota\:P \rarrow
P\oplus Q$ and $\pi\: P\oplus Q\rarrow P$ are the embedding of and
the projection onto the direct summand $P$ in $P\oplus Q$, then
the differential $\pi d\iota$ on $P$ makes it a QDG\+module.

 To prove that $I$ is a pseudo-equivalence, it suffices to show that
the $\Gamma$\+graded $B^\#$\+module $P\oplus P[-1]$ admits
a CDG\+module structure for any (finitely generated) projective
right $B^\#$\+module $P$.
 Define the right CDG\+module $Q$ over $B$ with the group $Q(X)$
consisting of formal expressions of the form $p'+d(p'')$, \ 
$p'$, $p''\in P(X)$, with $P(X)$ embedded into $Q(X)$ as the set
of all expressions $p+d(0)$.
 The differential $d$ on $Q(X)$, being restricted to $P(X)$, maps
$p+d(0)$ to $0+d(p)$ and $B$ acts on $P\subset Q$ as it acts on~$P$.
 The action of $B$ is extended from $P$ to $Q$ in the unique way
making the Leibniz rule satisfied, and the differential~$d$ is
extended from $P$ to $Q$ in the unique way making the equation on
$d^2$ hold (see~\cite[proof of Theorem~3.6]{Pkoszul} for explicit
formulas).
 There is a natural exact sequence of $\Gamma$\+graded
$B^\#$\+modules $0\rarrow P^\#\rarrow Q^\#\rarrow P^\#[-1]\rarrow 0$,
which splits, since $P^\#$ is projective.
\end{proof}

\begin{lemB}
 If $F\:C\rarrow D$ is a pseudo-equivalence of small CDG\+categories,
then the induced strict CDG\+functors $F^*\:D\modlq\rarrow C\modlq$
and\/ $\modrq D\rarrow\modrq C$ are strict equivalences
of CDG\+categories.
 For any QDG\+modules $N\in\modrq D$ and $M\in D\modlq$, the natural
map $N\ot_D M\rarrow F^*(N)\ot_C F^*(M)$ is an isomorphism of
precomplexes. 
 Besides, the induced DG\+functors $F^*\:D\modlc\rarrow C\modlc$ and\/
$\modrc D\rarrow\modrc C$ are equivalences of DG\+categories.
\end{lemB}

\begin{proof}
 First of all, it is obvious that if $F\:C^\#\rarrow D^\#$ is
a pseudo-equivalence of $\Gamma$\+graded categories, then
the induced functor of restriction of scalars between the categories
of (left or right) $\Gamma$\+graded modules over $D$ and $C$ is
an equivalence of $\Gamma$\+graded categories.
 These equivalences transform the functor of tensor product of
$\Gamma$\+graded modules over $D$ into the functor of tensor product
of $\Gamma$\+graded modules over~$C$.
 Thus, it remains to check that any QDG\+module over $C$ can be
extended to a QDG\+module over~$D$.
 And this is also straightforward.
\end{proof}

 More generally, one can see that the assertions of Lemma~B hold
for any CDG\+functor $F\:C\rarrow D$ that is a pseudo-equivalence
\emph{as a $\Gamma$\+graded functor} $C^\#\rarrow D^\#$.
 The assertions of both Lemmas A and~B remain valid if one replaces
finitely generated projective modules with finitely generated
free ones.

\begin{lemC} \hfuzz=4pt
 If a CDG\+functor $F\:C\rarrow D$ is a pseudo-equivalence of
CDG\+categories, then so is the CDG\+functor $F^\op\:
C^\op\rarrow D^\op$.
 If $k$\+linear CDG\+functors $F'\:C'\allowbreak\rarrow D'$ and
$F''\:C''\rarrow D''$ are pseudo-equivalences of CDG\+categories,
then so is the CDG\+functor $F'\ot F''\:C'\ot_k C''\rarrow
D'\ot_k D''$. \qed
\end{lemC}

\Section{Ext and Tor of the Second Kind}

 This section contains an exposition of the classical theory of
the two kinds of differential derived functors, largely
following~\cite{HMS}, except that we deal with CDG\+categories
rather than DG\+(co)algebras.
 The classical theory allows to establish an isomorphism
between the Hochschild (co)homology of the second kind of
a CDG\+category $B$ and the DG\+category $C$ of right
CDG\+modules over $B$ that are finitely generated and projective
as graded $B$\+modules.
 We also construct a natural map between the two kinds of Hochschild
(co)homology of any DG\+category~$C$ linear over a field~$k$.
 
\subsection{Ext and Tor of the first kind}  \label{ext-tor-first-kind}
 Given a DG\+category $D$, denote by $Z^0(D)$ the category whose
objects are the objects of $D$ and whose morphisms are the closed
(i.~e., annihilated by the differential) morphisms of degree~$0$
in~$D$.
 Let $H^0(D)$ denote the category whose objects are the objects
of $D$ and whose morphisms are the elements of the cohomology groups
of degree~$0$ of the complexes of morphisms in~$D$.
 The categories $Z^0(D)$ and $H^0(D)$ have preadditive category
structures (i.~e., the abelian group structures on the sets of
morphisms).
 In addition, these categories are endowed with the shift functors
$X\maps X[n]$ for all $n\in\Gamma$, provided that shifts of all
objects exist in~$D$ (see~\ref{cdg-categories-subsect}).
 Finally, let $H(D)$ denote the $\Gamma$\+graded category whose
objects are the objects of $D$ and whose morphisms are
the $\Gamma$\+graded groups of cohomology of the complexes of
morphisms in~$D$.

 Let $k$ be a commutative ring and $C$ be a small $k$\+linear
DG\+category.
 Let us endow the additive categories $Z^0(C\modld)$ and $Z^0(\modrd C)$
with the following exact category structures.
 A short sequence $M'\rarrow M\rarrow M''$ of DG\+modules and closed
morphisms between them  is exact if and only if \emph{both} the short
sequence of $\Gamma$\+graded $C^\#$\+modules $M'{}^\#\rarrow M^\#\rarrow
M''{}^\#$ and the short sequence of $\Gamma$\+graded $H(C)$\+modules
of cohomology $H(M')\rarrow H(M)\rarrow H(M'')$ are exact in the abelian
categories of $\Gamma$\+graded modules and their homogeneous
morphisms of degree~$0$.
 In other words, for any object $X\in C$ the sequence
$M'(X)\rarrow M(X)\rarrow M''(X)$ must be a short exact sequence of
complexes of $k$\+modules whose $\Gamma$\+graded cohomology modules 
also form a short exact sequence (i.~e., the boundary maps vanish).

 Denote the additive category $Z^0(k\modld)$ of $\Gamma$\+graded
complexes of $k$\+modules with its exact category structure
defined above by $\Comex(k\modl)$.
 Let $d$~denote the differentials on objects of $\Comex(k\modl)$.
 We will be interested in the derived categories
$\DD^-(\Comex(k\modl))$ and $\DD^+(\Comex(k\modl))$ of complexes,
bounded from above or below, over the exact category $\Comex(k\modl)$.
 The differential acting between the terms of a complex over
$\Comex(k\modl)$ will be denoted by~$\d$.

 The objects of $\DD^-(\Comex(k\modl))$ can be viewed as bicomplexes
with one grading by the integers bounded from above and the other
grading by elements of the group~$\Gamma$.
 (The differential~$d$ preserves the grading by the integers, while
changing the $\Gamma$\+valued grading; and the differential~$\d$
raises the grading by the integers by~$1$, while preserving
the $\Gamma$\+valued grading.)
 To any such bicomplex, one can assign its $\Gamma$\+graded total
complex, constructed by taking infinite direct sums along
the diagonals.
 This defines a triangulated functor from $\DD^-(\Comex(k\modl))$
to the unbounded derived category of $\Gamma$\+graded complexes of
$k$\+modules,
$$
 \Tot^\oplus\:\DD^-(\Comex(k\modl))\lrarrow \DD(k\modl).
$$

 Analogously, the objects of $\DD^+(\Comex(k\modl))$ can be viewed as
bicomplexes with one grading by the integers bounded from below and
the other grading by elements of the group~$\Gamma$.
 To any such bicomplex, one can assign its $\Gamma$\+graded total
complex, constructed by taking infinite products along the diagonals.
 This defines a triangulated functor
$$
 \Tot^\sqcap\:\DD^+(\Comex(k\modl))\lrarrow \DD(k\modl).
$$
 Any complex over $\Comex(k\modl)$ bounded from above (resp.,\ below)
that becomes exact (with respect to the differential~$\d$) after
passing to the cohomology of the $\Gamma$\+graded complexes of
$k$\+modules (with respect to the differential~$d$) is annihilated
by the functor $\Tot^\oplus$ (resp.,\ $\Tot^\sqcap$).

\begin{rem}
 The latter assertion does not hold for the total complexes of
unbounded complexes over $\Comex(k\modl)$, constructed by taking
infinite direct sums or products along the diagonals.
 That is the reason why we define the functors $\Tot^\oplus$ and
$\Tot^\sqcap$ for bounded complexes only.
 The assertion holds, however, for the functor of ``Laurent
totalization'' of unbounded complexes, which coincides with
$\Tot^\oplus$ for complexes bounded from above and with $\Tot^\sqcap$
for complexes bounded from below.
 See~\cite{HMS} and the introduction to~\cite{Pkoszul}
(cf.\ Remark~\ref{second-kind-general}).
\end{rem}

 Now consider the functor of two arguments
(see~\ref{qdg-modules-subsect})
\begin{equation} \label{dg-tensor-product}
 \ot_C\:Z^0(\modrd C)\times Z^0(C\modld)\lrarrow\Comex(k\modl).
\end{equation}
 We would like to construct its left derived functor
$$
 \ot_C^\L\:Z^0(\modrd C)\times Z^0(C\modld)
 \lrarrow\DD^-(\Comex(k\modl)).
$$
 For this purpose, notice that both exact categories $Z^0(\modrd C)$
and $Z^0(C\modld)$ have enough projective objects.
 Specifically, for any object $X\in C$ the representable DG\+module
$R_X\in Z^0(\modrd C)$ is projective, and so is the the cone of
the identity endomorphism of~$R_X$ (taken in the DG\+category
$\modrd C$).
 Any object of $Z^0(\modrd C)$ is the image of an admissible
epimorphism acting from an (infinite) direct sum of shifts of objects 
of the above two types.

 Given a right DG\+module $N$ and a left DG\+module $M$ over $C$,
choose a left projective resolution $Q_\bu$ of $N$ and a left
projective resolution $P_\bu$ of $M$ in the exact categories
$Z^0(\modrd C)$ and $Z^0(C\modld)$.
 When substituted as one of the arguments of the functor~$\ot_C$,
any projective object of one of the exact categories of DG\+modules
makes this functor an exact functor from the other exact
category of DG\+modules to the exact category $\Comex(k\modl)$.
 This allows to define $N\ot_C^\L M\in\DD^-(\Comex(k\modl))$
as the object represented either by the complex $Q_\bu\ot_C M$,
or by the complex $N\ot_C P_\bu$, or by the total complex of
the bicomplex $Q_\bu\ot_C P_\bu$.

 Analogously, consider the functor of two arguments
\begin{equation} \label{dg-hom}
 \Hom^C\:Z^0(C\modld)^\op\times Z^0(C\modld)\rarrow\Comex(k\modl),
\end{equation}
assigning to any two left DG\+modules over $C$ the complex of
morphisms between them as DG\+functors $C\rarrow k\modld$.
 We would like to construct its right derived functor
$$
 \R\!\Hom^C\:Z^0(C\modld)^\op\times Z^0(C\modld)
 \rarrow\DD^+(\Comex(k\modl)).
$$
 Notice that the exact category $Z^0(C\modld)$ has enough injective
objects.
 For any projective object $Q\in Z^0(\modrd C)$ and an injective
$k$\+module $I$, the object $\Hom_k(Q,I)\in Z^0(C\modld)$ is
injective, and any injective object in the exact category
$Z^0(C\modld)$ is a direct summand of an object of this type.
 To prove these assertions, it suffices to check that for any
DG\+module $M\in Z^0(\modrd C)$, any object $X\in C$, and any element
of $M(X)$ or $H(M)(X)$ there is a DG\+module $Q$ as above and
a closed morphism of DG\+modules $M\rarrow\Hom_k(Q,I)$ that is
injective on the chosen element.

 Given left DG\+modules $L$ and $M$ over $C$, choose a left
projective resolution $P_\bu$ of $L$ and a right injective
resolution $J^\bu$ of $M$ in the exact category $Z^0(C\modld)$.
 Substituting a projective object as the first argument or
an injective object as the second argument of the functor
$\Hom^C$, one obtains an exact functor from the exact category
of DG\+modules in the other argument to the exact category
$\Comex(k\modl)$.
 This allows to define $\R\!\Hom^C(L,M)\in\DD^+(\Comex(k\modl))$
as the object represented either by the complex $\Hom^C(P_\bu,M)$,
or by the complex $\Hom^C(L,J^\bu)$, or by the total complex of
the bicomplex $\Hom^C(P_\bu,J^\bu)$.

 Composing the derived functor $\ot_C^\L$ with the functor
$\Tot^\oplus$, we obtain the derived functor
$$
 \Tor^C\:Z^0(\modrd C)\times Z^0(C\modld)\lrarrow\DD(k\modl).
$$
 Similarly, composing the derived functor $\R\!\Hom^C$ with
the functor $\Tot^\sqcap$, we obtain the derived functor
$$
 \Ext_C\:Z^0(C\modld)^\op\times Z^0(C\modld)\lrarrow\DD(k\modl).
$$

 One can compute the derived functors $\Tor^C$ and $\Ext_C$ using
resolutions of a more general type than above.
 Specifically, let $N$ and $M$ be a left and a right DG\+module
over~$C$.
 Let $\dsb\rarrow F_2\rarrow F_1\rarrow F_0\rarrow M$ be a complex
of left DG\+modules over $C$ and (closed morphisms between them)
such that the complex of $\Gamma$\+graded $H(C)$\+modules
$\dsb\rarrow H(F_2)\rarrow H(F_1)\rarrow H(F_0)\rarrow H(M)\rarrow0$
is exact. 
 Assume that the DG\+modules $F_i$ are \emph{h\+flat} (homotopy
flat), i.~e., for any $i\ge0$ and any right DG\+module $R$ over $C$
such that $H(R)=0$ one has $H(R\ot_C F_i)=0$.
 Let $Q_\bu$ be a left projective resolution of the DG\+module $N$
in the exact category of right DG\+modules over~$C$.
 Then the natural maps $\Tot^\oplus(Q_\bu\ot_C F_\bu)\rarrow\Tot^\oplus
(Q_\bu\ot_C M)$ and $\Tot^\oplus(Q_\bu\ot_C F_\bu)\rarrow\Tot^\oplus
(N\ot_C F_\bu)$ are quasi-isomorphisms, so the $\Gamma$\+graded
complex of $k$\+modules $\Tot^\oplus(N\ot_C F_\bu)$ represents
the object $\Tor^C(N,M)$ in $\DD(k\modl)$.

 Analogously, let $L$ and $M$ be left DG\+modules over~$C$.
 Let $\dsb\rarrow P_2\rarrow P_1\rarrow P_0\rarrow L$ be a complex
of left DG\+modules over $C$ which becomes exact after passing to
the $\Gamma$\+graded cohomology modules.
 Assume that the DG\+modules $P_i$ are \emph{h\+projective} (homotopy
projective), i.~e., for any $i\ge0$ and any left DG\+module $R$
over $C$ such that $H(R)=0$ one has $H(\Hom^C(P_i,R))=0$.
 Then the complex of $k$\+modules $\Tot^\sqcap(\Hom^C(P_\bu,M))$
represents the object $\Ext_C(L,M)$ in $\DD(k\modl)$.
 Similarly, let $M\rarrow J^0\rarrow J^1\rarrow J^2\rarrow\dsb$ be
a complex of left DG\+modules over $C$ which becomes exact after
passing to the cohomology modules.
 Assume that the DG\+modules $J^i$ are \emph{h\+injective}, i.~e.,
for any $i\ge0$ and any left DG\+module $R$ over $C$ such that
$H(R)=0$ one has $H(\Hom^C(R,J^i))=0$.
 Then the complex of $k$\+modules $\Tot^\sqcap(\Hom^C(L,J^\bu))$
represents the object $\Ext_C(L,M)$.

 In particular, it follows that the functors $\Tor^C$ and $\Ext_C$
transform quasi-isomorphisms of DG\+modules (i.~e., morphisms of
DG\+modules inducing isomorphisms of the $\Gamma$\+graded
cohomology modules) in any of their arguments into isomorphisms
in $\DD(k\modl)$.

 Furthermore, consider the case when the complex of morphisms
between any two objects of $C$ is an h\+flat complex of
$k$\+modules.
 Then for any left DG\+module $M$ over $C$ such that the complex
of $k$\+modules $M(X)$ is h\+flat for any object $X\in C$,
the bar-construction
\begin{multline*}\textstyle
 \dsb\lrarrow\bigoplus_{Y,Z\in C}C(X,Y)\ot_k C(Y,Z)
 \ot_k M(Z) \\ \textstyle
 \lrarrow \bigoplus_{Y\in C} C(X,Y)\ot_k M(Y)\lrarrow M(X),
\end{multline*}
where we use the simplifying notation $C(X,Y)=\Hom_C(Y,X)$ for
any objects $X$, $Y\in C$, defines a left resolution of
the DG\+module $M$ which consists of h\+flat DG\+modules over $C$
and remains exact after passing to the cohomology modules.
 Thus, for any right DG\+module $N$ over $C$ the total complex of
the bar-complex
$$\textstyle
 \dsb\lrarrow\bigoplus_{Y,Z\in C} N(Y)\ot_k C(Y,Z)\ot_k M(Z)
 \lrarrow\bigoplus_{X\in C} N(Y)\ot_k M(Y),
$$
constructed by taking infinite direct sums along the diagonals,
represents the object $\Tor^C(N,M)$ in $\DD(k\modl)$.
 The h\+flatness condition on the DG\+module $M$ can be replaced
with the similar condition on the DG\+module~$N$.

 Analogously, assume that the complex of morphisms between any two
objects of $C$ is an h\+projective complex of $k$\+modules.
 Let $L$ and $M$ be left DG\+modules over $C$ such that either
the complex of $k$\+modules $L(X)$ is h\+projective for any
object $X\in C$ or the complex of $k$\+modules $M(X)$ is h\+injective
for any object $X\in C$.
 Then the total complex of the cobar-complex
$$\textstyle
 \prod_{X\in C}\Hom_k(L(X),M(X))\mskip-.1\thinmuskip
 \lrarrow\mskip -.1\thinmuskip
 \prod_{X,Y\in C}\Hom_k(C(X,Y)\ot_k L(Y)\;M(X))
 \mskip-.1\thinmuskip\lrarrow\mskip-.1\thinmuskip\dsb,
$$
constructed by taking infinite products along the diagonals,
represents the object $\Ext_C(L,M)$ in $\DD(k\modl)$.

 Given a $k$\+linear DG\+functor $F\:C\rarrow D$, a right DG\+module
$N$ over $D$, and a left DG\+module $M$ over $D$, there is a natural
morphism
\begin{equation} \label{tor-first-kind-F-star}
 \Tor^C(F^*N,F^*M)\lrarrow\Tor^D(N,M)
\end{equation}
in $\DD(k\modl)$.
 Analogously, given a $k$\+linear DG\+functor $F\:C\rarrow D$ and
left DG\+modules $L$ and $M$ over $D$, there is a natural morphism 
\begin{equation} \label{ext-first-kind-F-star}
 \Ext_D(L,M)\lrarrow\Ext_C(F^*L,F^*M)
\end{equation}
in $\DD(k\modl)$.
 If the functor $H(F)\:H(C)\rarrow H(D)$ is a pseudo-equivalence of
$\Gamma$\+graded categories, then the natural morphisms between
the objects $\Tor$ and $\Ext$ over $C$ and $D$ are isomorphisms
for any DG\+modules $L$, $M$, and~$N$.
 This follows from the fact that the similar morphisms between
the objects $\Tor$ and $\Ext$ over $H(C)$ and $H(D)$ are isomorphisms.

\subsection{Ext and Tor of the second kind: general case}
\label{second-kind-general}
 Let $B$ be a small $k$\+linear CDG\+category.
 Then the categories $Z^0(B\modlc)$ and $Z^0(\modrc B)$ of
(left and right) CDG\+mod\-ules over $B$ and closed morphisms of
degree~$0$ between them are abelian.
 In particular, consider the abelian category $Z^0(k\modlc)$ of
$\Gamma$\+graded complexes of $k$\+modules and denote it by
$\Comab(k\modl)$.
 We will be interested in the derived categories
$\DD^-(\Comab(k\modl))$ and $\DD^+(\Comab(k\modl))$ of complexes,
bounded from above or below, over the abelian category
$\Comab(k\modl)$.

 The objects of $\DD^-(\Comab(k\modl))$ can be viewed as bicomplexes
with one grading by the integers bounded from above and the other
grading by elements of the group~$\Gamma$.
 To any such bicomplex, one can assign its $\Gamma$\+graded total
complex, constructed by taking infinite products along the diagonals.
 This defines a triangulated functor
$$
 \Tot^\sqcap\:\DD^-(\Comab(k\modl))\lrarrow\DD(k\modl).
$$
 Analogously, the objects of $\DD^+(\Comab(k\modl))$ can be viewed
as bicomplexes with one grading by the integers bounded from below
and the other grading by elements of the group~$\Gamma$.
 To any such bicomplex, one can assign its $\Gamma$\+graded total
complex, constructed by taking infinite direct sums along
the diagonals.
 This defines a triangulated functor
$$
 \Tot^\oplus\:\DD^+(\Comab(k\modl))\lrarrow\DD(k\modl).
$$

\begin{rem} \emergencystretch=3em\hbadness=5000
 The functors of total complexes of unbounded complexes over
$\Comab(k\modl)$, constructed by taking infinite direct sums or
infinite products along the diagonals, are not well-defined on
the derived category $\DD(\Comab(k\modl))$.
 The procedure of ``Laurent totalization'' of unbounded complexes,
which coincides with $\Tot^\sqcap$ for complexes bounded from
above and with $\Tot^\oplus$ for complexes bounded from below,
defines a functor on $\DD(\Comab(k\modl))$, though.
 Notice that this Laurent totalization is different from
the one discussed in Remark~\ref{ext-tor-first-kind}
(the chosen direction along the diagonals is opposite in
the two cases).
\end{rem}

 Now consider the functor of two arguments
(see~\ref{qdg-modules-subsect})
\begin{equation} \label{cdg-tensor-product}
 \ot_B\:Z^0(\modrc B)\times Z^0(B\modlc)\lrarrow\Comab(k\modl).
\end{equation}
 We would like to construct its left derived functor
$$
 \ot_B^\L\:Z^0(\modrc B)\times Z^0(B\modlc)\lrarrow
 \DD^-(\Comab(k\modl)).
$$

 Notice that the abelian categories $Z^0(\modrc B)$ and $Z^0(B\modlc)$
have enough projective objects.
 More precisely, for any projective left $\Gamma$\+graded module $P$
over $B^\#$ the corresponding freely generated CDG\+module $Q$,
as constructed in the proof of Lemma~\ref{pseudo-equi-subsect}.A,
is a projective object of $Z^0(B\modlc)$.
 Any projective object in $Z^0(B\modlc)$ is a direct summand of
an object of this type.
 For any projective object $Q$ in $Z^0(B\modlc)$, the underlying
left $\Gamma$\+graded $B^\#$\+module $Q^\#$ is projective.

 Let us call a left $\Gamma$\+graded $B^\#$\+module $P^\#$
\emph{flat} if the functor of tensor product with $P^\#$ over $B^\#$
is exact on the abelian category of right $\Gamma$\+graded
$B^\#$\+modules.
 Given a left CDG\+module $P$ over $B$, if the left $B^\#$\+module
$P^\#$ is flat, then the functor of tensor product with $P$ is
exact as a functor $Z^0(\modrc B)\rarrow\Comab(k\modl)$.
 Any projective $\Gamma$\+graded $B^\#$\+module is flat.

 Given a right CDG\+module $N$ and a left CDG\+module $M$ over $B$,
choose a left resolution $Q_\bu$ of $N$ in $Z^0(\modrc B)$ and
a left resolution $P_\bu$ of $M$ in $Z^0(B\modlc)$ such that
the $\Gamma$\+graded $B^\#$\+modules $Q_i^\#$ and $P_i^\#$ are flat.
 In view of the above remarks, we can define 
$N\ot_B^\L M\in\DD^-(\Comab(k\modl))$ as the object represented
either by the complex $Q_\bu\ot_B M$, or by the complex $N\ot_B P_\bu$,
or by the total complex of the bicomplex $Q_\bu\ot_B P_\bu$.

 Analogously, consider the functor of two arguments
\begin{equation}  \label{cdg-hom}
 \Hom^B\:Z^0(B\modlc)^\op\times Z^0(B\modlc)\lrarrow
 \Comab(k\modl),
\end{equation}
assigning to any two left CDG\+modules over $B$ the complex of
morphisms between them as strict CDG\+functors $B\rarrow k\modlc$.
 We would like to construct its right derived functor
$$
 \R\!\Hom^B\:Z^0(B\modlc)^\op\times Z^0(B\modlc)\lrarrow
 \DD^+(\Comab(k\modl)).
$$

 Notice that the abelian category $Z^0(B\modlc)$ has enough injective
objects.
 For any injective object $J$ in $Z^0(B\modlc)$, the underlying
left $\Gamma$\+graded $B^\#$\+module $J^\#$ is injective.
 One can construct these injective CDG\+modules as the duals to
projective (or flat) right CDG\+modules (see the discussion of
injective DG\+modules in~\ref{ext-tor-first-kind}) or obtain them
as the CDG\+modules cofreely cogenerated by injective
$\Gamma$\+graded $B^\#$\+modules (see the construction of
injective resolutions in~\cite[proof of Theorem~3.6]{Pkoszul}).

 Given left CDG\+modules $L$ and $M$ over $B$, choose a left
resolution $P_\bu$ of $L$ and a right resolution $J^\bu$ of $M$
in $Z^0(B\modlc)$ such that the $\Gamma$\+graded $B^\#$\+modules
$P_i^\#$ are projective and the $\Gamma$\+graded $B^\#$\+modules
$J^i{}^\#$ are injective.
 Define $\R\!\Hom^B(L,M)\in\DD^-(\Comab(k\modl))$ as the object
represented either by the complex $\Hom^B(P_\bu,M)$, or by
the complex $\Hom^B(L,J^\bu)$, or by the total complex of
the bicomplex $\Hom^B(P_\bu,J^\bu)$.

 Composing the derived functor $\ot^\L_B$ with the functor
$\Tot^\sqcap$, we obtain the derived functor
$$
 \Tor^{B,I\!I}\:Z^0(\modrc B)\times Z^0(B\modlc)\lrarrow
 \DD(k\modl).
$$
 Similarly, composing the derived functor $\R\!\Hom^B$ with
the functor $\Tot^\oplus$, we obtain the derived functor
$$
 \Ext^{I\!I}_B\:Z^0(B\modlc)^\op\times Z^0(B\modlc)\lrarrow
 \DD(k\modl).
$$
 The derived functors $\Tor^{B,I\!I}$ and $\Ext^{I\!I}_B$ are
called the \emph{Tor and Ext of the second kind} of CDG\+modules
over~$B$.

 Notice that the derived functors $\ot_B^\L$ and $\R\!\Hom^B$
assign distinghuished triangles to short exact sequences of
CDG\+modules in any argument, hence so do the derived functors
$\Tor^{B,I\!I}$ and $\Ext^{I\!I}_B$.

 Given a $k$\+linear CDG\+functor $F\:B\rarrow C$, a right
CDG\+module $N$ over $C$, and a left CDG\+module $M$ over $C$,
there is a natural morphism
\begin{equation} \label{tor-second-kind-F-star}
 \Tor^{B,I\!I}(F^*N,F^*M)\lrarrow\Tor^{C,I\!I}(N,M)
\end{equation}
in $\DD(k\modl)$.
 Analogously, given a $k$\+linear CDG\+functor $F\:B\rarrow C$ and
left CDG\+modules $L$ and $M$ over $C$, there is a natural morphism
\begin{equation} \label{ext-second-kind-F-star}
 \Ext_C^{I\!I}(L,M)\lrarrow\Ext_B^{I\!I}(F^*L,F^*M)
\end{equation}
in $\DD(k\modl)$.
 If the functor $F^\#\:B^\#\rarrow C^\#$ is a pseudo-equivalence
of $\Gamma$\+graded categories, then these natural morphisms
are isomorphisms for any CDG\+modules $L$, $M$, and~$N$.

 Now let $C$ be a small $k$\+linear DG\+category.
 Then the identity functors from the exact categories
$Z^0(\modrd C)$ and $Z^0(C\modld)$ to the abelian categories
$Z^0(\modrc C)$ and $Z^0(C\modlc)$ are exact, so any
resolution in $Z^0(\modrd C)$ or $Z^0(C\modld)$ is also
a resolution in $Z^0(\modrc C)$ or $Z^0(C\modlc)$.
 Besides, any DG\+module that is projective or injective in
the exact category $Z^0(\modrd C)$ or $Z^0(C\modld)$ is also
projective or injective as a $\Gamma$\+graded $C^\#$\+module.
 It follows that there are natural morphisms
\begin{align} \label{tor-first-second}
 \Tor^C(N,M)&\lrarrow\Tor^{C,I\!I}(N,M) \\
\intertext{and} \label{ext-first-second}
 \Ext_C^{I\!I}(L,M)&\lrarrow\Ext_C(L,M)
\end{align}
in $\DD(k\modl)$ for any DG\+modules $L$, $M$, and $N$ over~$C$.

\subsection{Flat/projective case}  \label{second-kind-flat}
 Let $B$ be a small $k$\+linear CDG\+category, $N$ a right
CDG\+module over $B$, and $M$ a left CDG\+module over $B$.
 Consider the $\Gamma$\+graded complex of $k$\+modules
$\Br^\sqcap(N,B,M)$ constructed in the following way.
 As a $\Gamma$\+graded $k$\+module, $\Br^\sqcap(N,B,M)$ is obtained by
totalizing a bigraded $k$\+module with one grading by elements of
the group $\Gamma$ and the other grading by nonpositive integers,
the totalizing being performed by taking infinite products along
the diagonals.
 The component of degree $-i\in\Z$ of that bigraded module is
the $\Gamma$\+graded $k$\+module
$$\textstyle
 \bigoplus_{X_0,\dsc,X_i\in B} N(X_0)\ot_k B(X_0,X_1)\ot_k
 \dsb\ot_k B(X_{i-1},X_i)\ot_k M(X_i),
$$
where, as in~\ref{ext-tor-first-kind}, we use the simplifying
notation $B(X,Y)=\Hom_B(Y,X)$.

 The differential on $\Br^\sqcap(N,B,M)$ is the sum of the three
components $\d$, $d$, and $\delta$ given by the formulas
\begin{align*}
 &\d(n\ot b_1\ot\dsb\ot b_i\ot m) = nb_1\ot b_2\ot\dsb\ot b_i\ot m
 - n\ot b_1b_2\ot b_3\ot\dsb\ot b_i\ot m \\ &+ \dsb + 
 (-1)^{i-1} n\ot b_1\ot\dsb\ot b_{i-2}\ot b_{i-1}b_i\ot m +
 (-1)^i n\ot b_1\ot\dsb \ot b_{i-1}\ot b_im,
\end{align*}
where the products $b_jb_{j+1}$ denote the composition
of morphisms in $B$ and the products $nb_1$ and $b_im$ denote
the action of morphisms in $B$ on the CDG\+modules,
\begin{align*}
 (-1)^id(n\ot b_1\ot\dsb\ot b_i\ot m) &= d(n)\ot b_1\ot\dsb\ot b_i\ot m
 \\ &+ (-1)^{|n|} n\ot d(b_1)\ot b_2\ot\dsb\ot b_i\ot m + \dsb 
 \\ &+ (-1)^{|n|+|b_1|+\dsb+|b_i|} n\ot b_1\ot\dsb\ot b_i\ot d(m),
\end{align*}
and
\begin{multline*}
 \delta(n\ot b_1\ot \dsb\ot b_i\ot m) = n\ot h\ot b_1\ot\dsb\ot b_i
 \ot m \\ - n\ot b_1\ot h\ot b_2\ot\dsb\ot b_i\ot m + \dsb +
 (-1)^i n\ot b_1\ot\dsb\ot b_i\ot h\ot m.
\end{multline*}

\begin{propA}
 Assume that all the $\Gamma$\+graded $k$\+modules $B^\#(X,Y)$
are flat, and either all the $\Gamma$\+graded $k$\+modules $N^\#(X)$
are flat, or all the $\Gamma$\+graded $k$\+modules $M^\#(X)$ are flat.
 Then the complex\/ $\Br^\sqcap(N,B,M)$ represents the object\/
$\Tor^{B,I\!I} (N,M)$ in the derived category\/ $\DD(k\modl)$.
\end{propA}

\begin{proof}
 Choose a left resolution $Q_\bu$ of the right CDG\+module $N$
and a left resolution $P_\bu$ of the left CDG\+module $M$ such that
the $\Gamma$\+graded $B^\#$\+modules $P_j^\#$ and $Q_j^\#$ are flat.
 Consider the tricomplex $\Br^\sqcap(Q_\bu,B,P_\bu)$ and construct
its $\Gamma$\+graded total complex by taking infinite products along
the diagonals.
 Then this total complex maps naturally to both the complex
$\Br^\sqcap(N,B,M)$ and the total complex
$\Tot^\sqcap(Q_\bu\ot_B P_\bu)$ of the tricomplex $Q_\bu\ot_B P_\bu$,
constructed also by taking infinite products along the diagonals.
 These morphisms of $\Gamma$\+graded complexes are both
quasi-isomorphisms.
 Cf.\ the proof of Proposition~\ref{hochschild-subsect}.A below,
where some additional details can be found.
\end{proof}

 Let $F\:B\rarrow C$ be a $k$\+linear CDG\+functor, $N$ be a right
CDG\+module over $C$, and $M$ be a left CDG\+module over~$C$.
 Then there is a natural morphism of complexes of $k$\+modules
$F_*\:\Br^\sqcap(F^*N,B,F^*M)\rarrow\Br^\sqcap(N,C,M)$ given by
the rule
\begin{multline} \label{bar-cdg-functorial}
\textstyle
 F_*(n\ot b_1\ot\dsb\ot b_i\ot m) = \sum_{j_0,\dsc,j_i=0}^\infty 
 (-1)^{\rho(j_0,\dsc,j_i;\.|n|,|b_1|,\dsc,|b_i|)} \\ 
 n\ot a^{\ot j_0}\ot F(b_1)\ot a^{\ot j_1}\ot\dsb\ot F(b_i)
 \ot a^{\ot j_i}\ot m,
\end{multline}
where
\begin{multline} \label{rho-sign-formula}
 \rho(j_0,\dsc,j_i;\.t_0,t_1,\dsc,t_i) = 
 (j_0+\dsb+j_i-1)(j_0+\dsb+j_i)/2 \\ +j_0(i+1)+j_1i+\dsb+j_i
 +j_0t_0+j_1(t_0+t_1)+\dsb+j_i(t_0+t_1+\dsb+t_i).
\end{multline}
 The image of an arbitrary element in $\Br^\sqcap(F^*N,B,F^*M)$ is
constructed as the sum of the images of (the infinite number of)
its bihomogeneous components, the sum being convergent
bidegree-wise in $\Br^\sqcap(N,C,M)$.

 Suppose the CDG\+categories $B$ and $C$ satisfy
the assumptions of Proposition~A, and so does one of
the CDG\+modules $N$ and~$M$.
 Then the morphism of bar-complexes $F_*$ represents
the morphism~\eqref{tor-second-kind-F-star} of the objects $\Tor$
in $\DD(k\modl)$.

 Now let $L$ and $M$ be left CDG\+modules over~$B$.
 Consider the $\Gamma$\+graded complex of $k$\+modules
$\Cb^\oplus(L,B,M)$ constructed as follows.
 As a $\Gamma$\+graded $k$\+module, $\Cb^\oplus(L,B,M)$ is obtained by
totalizing a bigraded $k$\+module with one grading by elements of
the group $\Gamma$ and the other grading by nonnegative integers,
the totalizing being done by taking infinite direct sums along
the diagonals.
 The component of degree $i\in\Z$ of that bigraded module is
the $\Gamma$\+graded $k$\+module
$$ \textstyle
 \prod_{X_0,\dsc,X_i\in B} \Hom_k(B(X_0,X_1)\ot_k\dsb\ot_k
 B(X_{i-1},X_i)\ot_k L(X_i)\;M(X_0)).
$$

 The differential on $\Cb^\oplus(L,B,M)$ is the sum of the three
components $\d$, $d$, and $\delta$ given by the formulas
\begin{gather*}
\begin{split}
 &(\d f)(b_1,\dsc,b_{i+1},l) = (-1)^{|f||b_1|} b_1f(b_2,\dsc,b_{i-1},l)
 - f(b_1b_2,b_3,\dsc,b_{i+1},l) \\ &+ \dsb + (-1)^i f(b_1,\dsc,b_{i-1},
 b_ib_{i+1},l) + (-1)^{i+1} f(b_1,\dsc,b_i,b_{i+1}l),
\end{split} \displaybreak[1]\\
\begin{split}
 &(-1)^i(df)(b_1,\dsc,b_i,l) = d(f(b_1,\dsc,b_i,l)) -
 (-1)^{|f|}f(db_1,b_2,\dsc,b_i,l) \\
 &- (-1)^{|f|+|b_1|} f(b_1,db_2,b_3,\dsc,b_i,l) - \dsb - 
 (-1)^{|f|+|b_1|+\dsc+|b_i|} f(b_1,\dsc,b_i,dl),
\end{split}
\end{gather*}
and
\begin{align*}
 (\delta f)(b_1,\dsc,b_{i-1},l) &=
 -f(h,b_1,\dsc,b_{i-1},l) \\ &+ f(b_1,h,\dsc,b_{i-1},l) - \dsb
 + (-1)^i f(b_1,\dsc,b_{i-1},h,l).
\end{align*}

\begin{propB}
 Assume that all the $\Gamma$\+graded $k$\+modules $B^\#(X,Y)$
are projective, and either all the $\Gamma$\+graded $k$\+modules
$L^\#(X)$ are projective, or all the $\Gamma$\+graded $k$\+modules
$M^\#(X)$ are injective.
 Then the complex\/ $\Cb^\oplus(L,B,M)$ represents the object\/
$\Ext_B^{I\!I}(L,M)$ in the derived category\/ $\DD(k\modl)$.
\end{propB}

\begin{proof}
 Choose a left resolution $P_\bu$ of the left CDG\+module
$L$ and a right resolution $J^\bu$ of the left CDG\+module $M$
such that the $\Gamma$\+graded $B^\#$\+modules $P_j^\#$ are
projective and the $\Gamma$\+graded $B^\#$\+modules $Q^j{}^\#$
are injective.
 Consider the tricomplex $\Cb^\oplus(P_\bu,B,J^\bu)$ and construct
its $\Gamma$\+graded total complex by taking infinite direct sums
along the diagonals.
 Both the complex $\Cb^\oplus(L,B,M)$ and the total complex
$\Tot^\oplus(\Hom^B(P_\bu,J^\bu))$ of the tricomplex
$\Hom^B(P_\bu,J^\bu)$ map quasi-isomorphically into the above
total complex.
\end{proof}

 Let $F\:B\rarrow C$ be a $k$\+linear CDG\+functor, and $L$ and $M$
be left CDG\+modules over~$C$.
 Then there is a natural morphism of complexes of $k$\+modules
$F^*\:\Cb^\op(L,C,\allowbreak M)\rarrow\Cb^\op(F^*L,B,F^*M)$ given by
the rule
\begin{multline} \label{cobar-cdg-functorial}\textstyle
 (F^*f)(b_1\ot\dsb b_i\ot l)=\sum_{j_0,\dsc,j_i=0}^\infty
 (-1)^{\lambda(j_0,\dsc,j_i;|f|,|b_1|,\dsc,|b_i|)} \\
 f(a^{\ot j_0}\ot F(b_1)\ot a^{\ot j_1}\ot\dsb
 \ot F(b_i)\ot a^{j_i}\ot l),
\end{multline}
where
\begin{align} \label{sigma-sign-formula} 
 \lambda(j_0,\dsc,j_i;\.t_0,t_1,\dsc,t_i) &= 
 j_0i+j_1(i-1)+\dsb+j_{i-1} \\ \notag
 &+j_0t_0+j_1(t_0+t_1)+\dsb+j_i(t_0+t_1+\dsb+t_i).
\end{align}

 Suppose the CDG\+categories $B$ and $C$ satisfy the assumptions
of Proposition~B, and so does one of the CDG\+modules $L$ and~$M$.
 Then the morphism of cobar-complexes $F^*$ represents
the morphism~\eqref{ext-second-kind-F-star} of the objects $\Ext$
in $\DD(k\modl)$.

 Denote by $\Br^\oplus(N,B,M)$ the $\Gamma$\+graded complex of
$k$\+modules constructed in the same way as $\Br^\sqcap(N,B,M)$,
except that the totalization is being done by taking infinite
direct sums along the diagonals.
 Similarly, denote by $\Cb^\sqcap(L,B,M)$ the $\Gamma$\+graded
complex of $k$\+modules constructed in the same way as
$\Cb^\oplus(L,B,M)$ except that the totalization is being done by
taking infinite products along the diagonals.
 
 Assume that $C$ is a small DG\+category in which the complex of
morphisms between any two objects is an h\+flat complex of flat
$k$\+modules, and either a right DG\+module $N$ or a left DG\+module
$M$ over $C$ is such that all the complexes of $k$\+modules $N(X)$
or $M(X)$ are h\+flat complexes of flat $k$\+modules.
 Then the natural map $\Br^\oplus(N,C,M)\rarrow\Br^\sqcap(N,C,M)$
represents the morphism $\Tor^C(N,M)\rarrow \Tor^{C,I\!I}(N,M)$
in $\DD(k\modl)$.

 Analogously, assume that the complex of morphisms between any two
objects in a DG\+category $C$ is an h\+projective complex of
projective $k$\+modules, and either a left DG\+module $L$ over $C$
is such that all the complexes of $k$\+modules $L(X)$ are
h\+projective complexes of projective $k$\+modules, or a left
DG\+module $M$ over $C$ is such that all the complexes of $k$\+modules
$M(X)$ are h\+injective complexes of injective $k$\+modules.
 Then the natural map $\Cb^\oplus(L,C,M)\rarrow\Cb^\sqcap(L,C,M)$
represents the morphism $\Ext_C^{I\!I}(L,M)\rarrow\Ext_C(L,M)$
in $\DD(k\modl)$.

 Notice that the complexes $\Br^\oplus(N,B,M)$ and $\Cb^\sqcap
(L,B,M)$ are \emph{not} functorial with respect to nonstrict
CDG\+functors between CDG\+categories $B$ because of the infinite
summation in the formulas \eqref{bar-cdg-functorial}
and~\eqref{cobar-cdg-functorial}.

\begin{propC}
 Let $B$ be a small $k$\+linear CDG\+category.
 Assume that the maps $k\rarrow\Hom_B(X,X)$ corresponding to
the curvature elements $h_X\in\Hom_B(X,X)$ admit $k$\+linear
retractions $\Hom_B(X,X)\rarrow k$, i.~e., they are embeddings
of $k$\+module direct summands.
 In particular, this holds when $k$ is a field and all the elements
$h_X$ are nonzero.
 Then for any CDG\+modules $L$, $M$ and $N$ the complexes
$\Br^\oplus(N,B,M)$ and $\Cb^\sqcap(L,B,M)$ are acyclic.
\end{propC}

\begin{proof}
 This follows from the fact that the differentials~$\delta$ on
the bigraded bar- and cobar-complexes are acyclic.
\end{proof}

\subsection{Hochschild (co)homology}  \label{hochschild-subsect}
 Let $B$ be a small $k$\+linear CDG\+category.
 Consider the CDG\+category $B\ot_k B^\op$; since it is naturally
isomorphic to its opposite CDG\+category, there is no need to
distinguish between the left and the right CDG\+modules over it.
 Furthermore, there is a natural (left) CDG\+module over
the CDG\+category $B\ot_k B^\op$ assigning to an object
$(X,Y^\op)\in B\ot_k B^\op$ the precomplex of $k$\+modules
$B(X,Y)=\Hom_B(Y,X)$.
 By an abuse of notation, we will denote this CDG\+module (as
well as the corresponding right CDG\+module) simply by~$B$.
 Assume that the $\Gamma$\+graded $k$\+modules $B^\#(X,Y)$
are flat for all objects $X$, $Y\in B$.

 The \emph{Hochschild homology of the second kind} $HH^{I\!I}_*(B,M)$
of a $k$\+linear CDG\+cate\-gory $B$ with coefficients in a (left)
CDG\+module $M$ over $B\ot_k B^\op$ is defined as the homology of
the object $\Tor^{B\ot_k B^\op\;I\!I}(B,M)\in\DD(k\modl)$.
 In particular, the Hochschild homology of the second kind of
the CDG\+module $M=B$ over $B\ot_k B^\op$ is called simply
the Hochschild homology of the second kind of the $k$\+linear
CDG\+category $B$ and denoted by $HH^{I\!I}_*(B,B) = HH^{I\!I}_*(B)$.
 
 The \emph{Hochschild cohomology of the second kind} $HH^{I\!I\;*}
(B,M)$ of a $k$\+linear CDG\+category $B$ with coefficients in
a (left) CDG\+module $M$ over $B\ot_k B^\op$ is defined as
the cohomology of the object
$\Ext_{B\ot_k B^\op}^{I\!I}(B,M)\in\DD(k\modl)$.
 In particular, the Hochschild cohomology of the second kind of
the CDG\+module $M=B$ over $B\ot_k B^\op$ is called simply
the Hochschild cohomology of the second kind of the $k$\+linear
CDG\+category $B$ and denoted by $HH^{I\!I\;*}(B,B)=HH^{I\!I\;*}(B)$.

\begin{rem}
 We define the Hochschild (co)homology of the second kind for
CDG\+cate\-gories $B$ satisfying the above flatness assumption only,
even though our definition makes sense without this requirement.
 In fact, this assumption is never used in this paper (except in
the discussion of explicit complexes below in this section, which
requires a stronger projectivity assumption in the cohomology case
anyway).
 However, we believe that our definition is not the \emph{right}
one without the flatness assumption, since one is not supposed to use
underived nonexact functors when defining (co)homology theories.
 So to define the Hochschild (co)homology of the second kind in
the general case one would need to replace a CDG\+category $B$ with
a CDG\+category, equivalent to it in some sense and satisfying
the flatness requirement.
 We do not know how such a replacement could look like.
 The analogue of this procedure for Hochschild (co)homology of
the first kind is well-known (in this case it suffices to replace
a DG\+category $C$ with a quasi-equivalent DG\+category with
h\+flat complexes of morphisms; see below).
\end{rem}

 By the result of~\ref{second-kind-flat}, the Hochschild homology
$HH^{I\!I}_*(B,M)$ is computed by the explicit bar-complex
$\Br^\sqcap(B\;B\ot_k B^\op\;M)$.
 When the $\Gamma$\+graded $k$\+modules $B^\#(X,Y)$ are
projective for all objects $X$, $Y\in B$, the Hochschild cohomology
$HH^{I\!I\;*}(B,M)$ is computed by the explicit cobar-complex
$\Cb^\oplus(B\;B\ot_k B^\op\;M)$.
 However, these complexes are too big and apparently not very useful.

 There are smaller and much more important complexes computing
the Hochschild (co)homology, namely, the Hochschild complexes.
 The homological Hochschild complex of the second kind
$\Hoch_\bu^\sqcap(B,M)$ is constructed in the following way.
 As a $\Gamma$\+graded $k$\+module, $\Hoch_\bu^\sqcap(B,M)$ is
obtained by taking infinite products along the diagonals of
a bigraded $k$\+module with one grading by elements of the group
$\Gamma$ and the other grading by nonpositive integers.
 The component of degree $-i\in\Z$ of that bigraded $k$\+module is
the $\Gamma$\+graded $k$\+module
$$\textstyle
 \bigoplus_{X_0,\dsc,X_i\in B} M(X_i,X_0^\op)\ot_k B(X_0,X_1)
 \ot_k\dsb\ot_k B(X_{i-1},X_i).
$$
 The differential on $\Hoch_\bu^\sqcap(B,M)$ is the sum of the three
components $\d$, $d$, and $\delta$ given by the formulas
\begin{gather*}
\begin{split}
 \d(m\ot b_1\ot\dsb \ot b_i) &= mb_1\ot b_2\ot\dsb\ot b_i -
 m\ot b_1b_2\ot b_3\ot\dsb\ot b_i \\&+ \dsb + (-1)^{i-1}
 m\ot b_1\ot\dsb\ot b_{i-2}\ot b_{i-1}b_i \\&+ (-1)^{i+|b_i|
 (|m|+|b_1|+\dsb+|b_{i-1}|)} b_im\ot b_1\ot \dsb\ot b_{i-1},
\end{split} \displaybreak[1]\\
\begin{split}
 (-1)^i d(m\ot b_1\ot\dsb\ot b_i) &= d(m)\ot b_1\ot\dsb\ot b_i +
 (-1)^{|m|} m\ot d(b_1)\ot b_2\ot\dsb\ot b_i \\&+ \dsb +
 (-1)^{|m|+|b_1|+\dsb+|b_{i-1}|} m\ot b_1\ot\dsb\ot b_{i-1}\ot d(b_i),
\end{split}
\end{gather*}
and
\begin{multline*}
 \delta(m\ot b_1\ot\dsb\ot b_i) = m\ot h\ot b_1\ot\dsb\ot b_i \\
 - m\ot b_1\ot h\ot b_2\ot\dsb\ot b_i + \dsb + (-1)^i m\ot b_1\ot\dsb
 \ot b_i\ot h.
\end{multline*}

\begin{propA}
 The homology of the complex\/ $\Hoch_\bu^\sqcap(B,M)$ is naturally
isomorphic to the Hochschild homology of the second kind
$HH_*^{I\!I}(B,M)$ as a $\Gamma$\+graded $k$\+module.
\end{propA}

\begin{proof}
 Choose a left resolution $P_\bu$ of the CDG\+module $M$
such that the $\Gamma$\+graded $B^\#\ot_k B^\#{}^\op$\+modules
$P_j^\#$ are flat.
 Consider the bicomplex $\Hoch_\bu^\sqcap(B,P_\bu)$ and construct
its total complex by taking infinite products along the diagonals.
 This total complex maps naturally to both the complex
$\Hoch_\bu^\sqcap(B,M)$ and the total complex of the bicomplex
$B\ot_{B\ot_k B^\op}P_\bu$, constructed by taking infinite
products along the diagonals.
 These morphisms of $\Gamma$\+graded complexes are both
quasi-isomorphisms.

 Indeed, the morphism $\Hoch_\bu^\sqcap(B,P_\bu)\rarrow
\Hoch_\bu^\sqcap(B,M)$ is a quasi-isomorphism, because the functor
$\Hoch_\bu^\sqcap(B,{-})$ transforms exact sequences of
CDG\+modules over $B\ot_k B^\op$ into exact sequences of
complexes.
 The morphism $\Hoch_\bu^\sqcap(B,P_\bu)\rarrow B\ot_{B\ot_k B^\op}
P_\bu$ is a quasi-isomorphism, since the morphism
$\Hoch_\bu^\sqcap(B,P)\rarrow P$ is a quasi-isomorphism for any
CDG\+module $P$ over $B\ot_k B^\op$ such that the $\Gamma$\+graded
$B^\#\ot_k B^\#{}^\op$\+module $P^\#$ is flat.
 The latter assertion follows from the similar statement for
the bigraded Hochschild complex of the $\Gamma$\+graded
$B^\#\ot_k B^\#{}^\op$\+module $P^\#$ with the differential~$\d$.
\end{proof}

 Let $F\:B\rarrow C$ be a $k$\+linear CDG\+functor and $M$ be
a CDG\+module over $C\ot_k C^\op$.
 Let us denote the CDG\+module $(F\ot F^\op)^*M$ over
$B\ot_k B^\op$ simply by $F^*M$.
 There is a natural morphism of complexes of $k$\+modules
$F_*\:\Hoch_\bu^\sqcap(B,F^*M)\rarrow\Hoch_\bu^\sqcap(C,M)$
defined by the rule
\begin{multline} \label{ho-hoch-cdg-functorial} \textstyle
 F_*(m\ot b_1\ot\dsb\ot b_i) = \sum_{j_0,\dsc,j_i=0}^\infty 
 (-1)^{\rho(j_0,\dsc,j_i;\.|m|,|b_1|,\dsc,|b_i|)} \\ 
 m\ot a^{\ot j_0}\ot F(b_1)\ot a^{\ot j_1}\ot\dsb\ot F(b_i)
 \ot a^{\ot j_i},
\end{multline}
where the value of~$\rho$ in the exponent is given by
the formula~\eqref{rho-sign-formula}.
 The image of an arbitrary element in $\Hoch_\bu^\sqcap(B,F^*M)$ is
constructed as the sum of the images of (the infinite number of)
its bihomogeneous components, the sum being convergent
bidegree-wise in $\Hoch^\bu(C,M)$.

 The morphism $F_*$ of Hochschild complexes computes the map of
Hochschild homology
\begin{equation} \label{ho-hoch-second-kind-F-star}
 HH_*^{I\!I}(B,F^*M)\lrarrow HH_*^{I\!I}(C,M)
\end{equation}
obtained by passing to the homology in the morphism of $\Tor$
objects~\eqref{tor-second-kind-F-star} for the CDG\+functor
$F\ot F^\op$.
 Furthermore, there is a natural closed morphism $B\rarrow F^*C$ of
CDG\+modules over $B\ot_k B^\op$, inducing a map of Hochschild
homology
\begin{equation} \label{ho-hoch-second-kind-of-itself-F-star}
 HH_*^{I\!I}(B)\lrarrow HH_*^{I\!I}(C)
\end{equation}
and a morphism of Hochschild complexes $F_*\:\Hoch_\bu^\sqcap(B,B)
\rarrow\Hoch_\bu^\sqcap(C,C)$ computing this homology map.

 The cohomological Hochschild complex of the second kind
$\Hoch^{\oplus,\bu}(B,M)$ is constructed as follows.
 As a $\Gamma$\+graded $k$\+module, $\Hoch_\bu^\sqcap(B,M)$ is
obtained by taking infinite direct sums along the diagonals
of a bigraded $k$\+module with one grading by elements of
the group $\Gamma$ and the other grading by nonnegative integers.
 The component of degree $i\in\Z$ of that bigraded $k$\+module
is the $\Gamma$\+graded $k$\+module
$$\textstyle
 \prod_{X_0,\dsc,X_i\in B}\Hom_k(B(X_0,X_1)\ot_k\dsb\ot_k
 B(X_{i-1},X_i)\;M(X_0,X_i^\op)).
$$

 The differential on $\Hoch^{\oplus,\bu}(B,M)$ is the sum of
the three components $\d$, $d$, and $\delta$ given by the formulas
\begin{gather*}
\begin{split}
(\d f)(b_1,\dsc,b_{i+1}) &= (-1)^{|f||b_1|}b_1 f(b_2,\dsc,b_{i+1})
- f(b_1b_2,b_3,\dsc,b_{i+1}) \\
&+ \dsb + (-1)^i f(b_1,\dsc,b_{i-1},b_ib_{i+1}) +
(-1)^{i+1}f(b_1,\dsc,b_i)b_{i+1},
\end{split} \displaybreak[1]\\
\begin{split}
(-1)^i(df)(b_1,\dsc,b_i) &= d(f(b_1,\dsc,b_i)) -
(-1)^{|f|}f(db_1,b_2,\dsc,b_i) \\ &- \dsb - 
(-1){}^{|f|+|b_1|+\dsb+|b_{i-1}|}f(b_1,\dsc,b_{i-1},db_i),
\end{split}
\end{gather*}
and
\begin{align*}
(\delta f)(b_1,\dsc,b_{i-1}) &= - f(h,b_1,\dsc,b_{i-1}) \\
&+ f(b_1,h,b_2,\dsc,b_{i-1}) - \dsb + (-1)^if(b_1,\dsc,b_{i-1},h).
\end{align*}

\begin{propB}
 Assume that all the $\Gamma$\+graded $k$\+modules $B^\#(X,Y)$
are projective.
 Then the cohomology of the complex\/ $\Hoch^{\oplus,\bu}(B,M)$ is
naturally isomorphic to the Hochschild cohomology of the second
kind $HH^{I\!I\;*}(B,M)$ as a $\Gamma$\+graded $k$\+module.
\end{propB}

\begin{proof}
 Choose a right resolution $J^\bu$ of the CDG\+module $M$
such that the $\Gamma$\+graded $B^\#\ot_k B^\#{}^\op$\+modules
$J^j{}^\#$ are injective.
 Consider the bicomplex $\Hoch^{\oplus,\bu}(B,J^\bu)$ and
construct its total complex by taking infinite direct sums along
the diagonals.
 Both the complex $\Hoch^{\oplus,\bu}(B,M)$ and the total complex
of the bicomplex $\Hom^{B\ot_k B^\op}(B,J^\bu)$ map
quasi-isomorphically into the above total complex. 
\end{proof}

 For a $k$\+linear CDG\+functor $F\:B\rarrow C$ and a CDG\+module
$M$ over $C\ot_k C^\op$, there is a natural morphism of complexes
of $k$\+modules $F^*\:\Hoch^{\oplus,\bu}(C,M)\rarrow
\Hoch^{\oplus,\bu}(B,F^*M)$ defined by the rule
\begin{multline} \label{coho-hoch-cdg-functorial}\textstyle
 (F^*f)(b_1\ot\dsb b_i)=\sum_{j_0,\dsc,j_i=0}^\infty
 (-1)^{\lambda(j_0,\dsc,j_i;|f|,|b_1|,\dsc,|b_i|)} \\
 f(a^{\ot j_0}\ot F(b_1)\ot a^{\ot j_1}\ot\dsb
 \ot F(b_i)\ot a^{j_i}), 
\end{multline}
where the value of~$\lambda$ in the exponent is given by
the formula~\eqref{sigma-sign-formula}.

 Suppose the CDG\+categories $B$ and $C$ satisfy the assumptions
of Proposition~B\hbox{}.
 Then the morphism $F^*$ of Hochschild complexes computes
the map of Hochschild cohomology
\begin{equation}  \label{coho-hoch-second-kind-F-star}
 HH^{I\!I\;*}(C,M)\lrarrow HH^{I\!I\;*}(B,F^*M)
\end{equation}
obtained by passing to the cohomology in the morphism of $\Ext$
objects~\eqref{ext-second-kind-F-star} for the CDG\+functor
$F\ot F^\op$.
 Notice that, unlike the Hochschild homology, the Hochschild
cohomology of CDG\+categories $HH^{I\!I\;*}(B)$ is \emph{not}
functorial with respect to arbitrary CDG\+functors
$F\:B\rarrow C$.
 It \emph{is} contravariantly functorial, however, with respect
to CDG\+functors $F$ for which the functor $F^\#\:B^\#\rarrow C^\#$
is fully faithful, since the closed morphism of CDG\+modules
$B\rarrow F^*C$ is an isomorphism in this case.
 
 The (co)homology of the complexes $\Hoch_\bu^\sqcap(B,M)$ and
$\Hoch^{\oplus,\bu}(B,M)$ are what is called
the ``Borel--Moore Hochschild homology'' and
the ``compactly supported Hochschild cohomology'' in~\cite{CT}.

 Now denote by $\Hoch_\bu^\oplus(B,M)$ the $\Gamma$\+graded complex
of $k$\+modules constructed in the same way as $\Hoch_\bu^\sqcap
(B,M)$, except that the totalization is being done by taking
infinite direct sums along the diagonals.
 Similarly, denote by $\Hoch^{\sqcap,\bu}(B,M)$ the $\Gamma$\+graded
complex of $k$\+modules constructed in the same way as
$\Hoch^{\oplus,\bu}(B,M)$, except that the totalization is being
done by taking infinite products.
 The complexes $\Hoch_\bu^\oplus(B,M)$ and $\Hoch^{\sqcap,\bu}(B,M)$
play an important role when $B$ is a DG\+category, but apparently
not otherwise, as we will see below.

 Let $C$ be a small $k$\+linear DG\+category.
 Assume that the complexes of $k$\+modules $C(X,Y)$ are h\+flat
for all objects $X$, $Y\in C$.
 The (conventional) Hochschild homology (of the first kind)
$HH_*(C,M)$ of a $k$\+linear DG\+category $C$ with coefficients
in a DG\+module $M$ over $C\ot_k C^\op$ is the homology of
the object $\Tor^{C\ot_k C^\op}(C,M)\in\DD(k\modl)$.
 In particular, the Hochschild homology of the DG\+module $M=C$
over $C$ is called simply the Hochschild homology of $C$ and
denoted by $HH_*(C,C)=HH_*(C)$.

 The (conventional) Hochschild cohomology (of the first kind)
$HH^*(C,M)$ of a $k$\+linear DG\+category $C$ with coefficients
in a DG\+module $M$ over $C\ot_k C^\op$ is the cohomology of
the object $\Ext_{C\ot_k C^\op}(C,M)\in\DD(k\modl)$.
 In particular, the Hochschild cohomology of the DG\+module $M=C$
over $C$ is called simply the Hochschild cohomology of $C$ and
denoted by $HH^*(C,C)=HH^*(C)$.

 Let $F\:C\rarrow D$ be a $k$\+linear DG\+functor between
DG\+categories whose complexes of morphisms are h\+flat complexes
of $k$\+modules.
 Then for any DG\+module $M$ over $D\ot_k D^\op$ passing to
the homology in the morphism of $\Tor$
objects~\eqref{tor-first-kind-F-star}
for the DG\+functor $F\ot F^\op$ provides a natural map
of $\Gamma$\+graded $k$\+modules
\begin{equation}  \label{ho-hoch-first-kind-F-star}
 HH_*(C,F^*M)\lrarrow HH_*(D,M).
\end{equation}
 Composing this map with the map induced by the closed morphism
$C\rarrow F^*D$ of DG\+modules over $C\ot_k C^\op$, we obtain
a natural map
\begin{equation}  \label{ho-hoch-first-kind-of-itself-F-star}
 HH_*(C)\lrarrow HH_*(D).
\end{equation}
 Passing to the cohomology in the morphism of $\Ext$
objects~\eqref{ext-first-kind-F-star}
for the DG\+functor $F\ot F^\op$ provides a natural map
\begin{equation}  \label{coho-hoch-first-kind-F-star}
 HH^*(D,M)\lrarrow HH^*(C,F^*M).
\end{equation}
 Unlike the Hochschild homology, the Hochschild cohomology of
DG\+categories $HH^*(C)$ is \emph{not} functorial with respect to
arbitrary DG\+functors $F\:C\rarrow D$.
 It \emph{is} contravariantly functorial, however, with respect to
DG\+functors $F$ such that the functor $H(F)\:H(C)\rarrow H(D)$ is
fully faithful, since the closed morphism of DG\+modules
$C\rarrow F^*D$ is a quasi-isomorphism in this case.

 When the functor $H(F)$ is a pseudo-equivalence of $\Gamma$\+graded
categories, the maps
(\ref{ho-hoch-first-kind-F-star}\+-\ref{coho-hoch-first-kind-F-star})
are isomorphisms, as is the natural map $HH^*(D)\rarrow HH^*(C)$.
 Indeed, under our assumptions on the DG\+categories $C$ and $D$
the $\Gamma$\+graded category $H(C\ot_k C^\op)$ is isomorphic to
$H(C)\ot_k H(C)^\op$ and similarly for $D$, so the assertion follows
from Lemma~\ref{pseudo-equi-subsect}.C and the results
of~\ref{ext-tor-first-kind}.

 Just as in~\ref{ext-tor-first-kind}, one shows that the complex
$\Hoch_\bu^\oplus(C,M)$ computes the Hochschild homology
$HH_*(C,M)$.
 The morphism of complexes $F_*\:\Hoch_\bu^\oplus(C,F^*M)\rarrow
\Hoch_\bu^\oplus(D,M)$ induced by a DG\+functor $F$ computes
the map of Hochschild homology~\eqref{ho-hoch-first-kind-F-star}.
 In particular the morphism of complexes $F_*\:\Hoch_\bu^\oplus
(C,C)\rarrow\Hoch_\bu^\oplus(D,D)$ induced by~$F$ computes
the map~\eqref{ho-hoch-first-kind-of-itself-F-star}.
 When all the complexes of morphisms in $C$ are h\+projective
complexes of $k$\+modules, the complex $\Hoch^{\sqcap,\bu}(C,M)$
computes the Hochschild cohomology $HH^*(C,M)$.
 When both DG\+categories $C$ and $D$ satisfy the same condition,
the morphism of complexes $F^*\:\Hoch^{\sqcap,\bu}(D,M)\rarrow
\Hoch^{\sqcap,\bu}(C,F^*M)$ computes
the map~\eqref{coho-hoch-first-kind-F-star}.

 When the complexes $C(X,Y)$ are h\+flat complexes of flat
$k$\+modules for all objects $X$, $Y\in C$, both the Hochschild
(co)homology of the first and the second kind are defined for
any DG\+module $M$ over $C\ot_k C^\op$.
 In this case, there are natural morphisms of $\Gamma$\+graded
$k$\+modules
\begin{equation}  \label{hoch-first-second-kind}
 HH_*(C,M)\lrarrow HH^{I\!I}_*(C,M)
 \quad\text{and}\quad
 HH^{I\!I\;*}(C,M)\lrarrow HH^*(C,M)
\end{equation}
and, in particular,
\begin{equation}  \label{hoch-of-itself-first-second-kind}
 HH_*(C)\lrarrow HH^{I\!I}_*(C)
 \quad\text{and}\quad
 HH^{I\!I\;*}(C)\lrarrow HH^*(C).
\end{equation}
 All of these are obtained from the comparison
morphisms~(\ref{tor-first-second}--\ref{ext-first-second})
for the two kinds of functors $\Tor$ and $\Ext$.
 The morphism $HH_*(C,M)\rarrow HH^{I\!I}_*(C,M)$ is computed
by the morphism of complexes $\Hoch_\bu^\oplus(C,M)\rarrow
\Hoch^\sqcap_\bu(C,M)$.
 When the complexes $C(X,Y)$ are h\+projective complexes
of projective $k$\+modules for all objects $X$, $Y\in C$,
the morphism $HH^{I\!I\;*}(C,M)\rarrow HH^*(C,M)$ is computed
by the morphism of complexes $\Hoch^{\oplus,\bu}(C,M)\rarrow
\Hoch^{\sqcap,\bu}(C,M)$.

 On the other hand, assume that the maps $k\rarrow\Hom_B(X,X)$
corresponding to the curvature elements $h_X\in\Hom_B(X,X)$
are embeddings of $k$\+linear direct summads.
 Then for any CDG\+module $M$ over $B\ot_k B^\op$ the complexes
$\Hoch_\bu^\oplus(B,M)$ and $\Hoch^{\sqcap,\bu}(B,M)$ are
acyclic~\cite[Lemma~3.9 and Theorem~4.2(a)]{CT}.
 Notice that these complexes are \emph{not} functorial with respect
to nonstrict CDG\+functors between CDG\+categories $B$ because
of the infinite summation in the formulas
\eqref{ho-hoch-cdg-functorial} and~\eqref{coho-hoch-cdg-functorial}.
 The complexes $\Hoch_\bu^\oplus(B,M)$ and $\Br^\oplus
(B\;B\ot_k B^\op\;M)$ are \emph{not} quasi-isomorphic in
general, even when $k$~is a field, $B$ is a CDG\+algebra considered
as a CDG\+category with a single object, and $M=B$.
 Neither are the complexes $\Hoch^{\sqcap,\bu}(B,M)$ and
$\Cb^\sqcap(B\;B\ot_k B^\op\;M)$.

\subsection{Change of grading group}  \label{change-grading-group}
 Let us first introduce some terminology that will be used
throughout the rest of the paper.
 A $\Gamma$\+graded module $N^\#$ over a $\Gamma$\+graded category
$B^\#$ is said to have \emph{flat dimension~$d$} if $d$~is
the minimal length of a left flat resolution of~$N^\#$ in the abelian
category of $\Gamma$\+graded $B^\#$\+modules, or equivalently,
the functor of tensor product with $N^\#$ over $B^\#$ has
the homological dimension~$d$.
\emph{Projective} and \emph{injective} dimensions of
$\Gamma$\+graded $B^\#$\+modules are defined in the similar way.
 The \emph{left homological dimension} of a $\Gamma$\+graded category
$B^\#$ is the homological dimension of the abelian category of
$\Gamma$\+graded left $B^\#$\+modules, and the \emph{weak homological
dimension} of $B^\#$ is the homological dimension of the functor
of tensor product of $\Gamma$\+graded modules over~$B^\#$.

 Let $(\Gamma,\sigma,\boldsymbol{1})$ and $(\Gamma',\sigma',
\boldsymbol{1}')$ be two different grading group data
(see~\ref{grading-group}) and $\phi\:\Gamma\rarrow\Gamma'$ be
a morphism of abelian groups taking $\boldsymbol{1}$ to
$\boldsymbol{1}'$ such that $\sigma$ is the pull-back of~$\sigma'$
by~$\phi$.
 Then to any $\Gamma'$\+graded $k$\+module $V'$ one can assign
a $\Gamma$\+graded $k$\+module $\phi^*V'$ defined by the rule
$(\phi^*V)^n=V'{}^{\phi(n)}$ for $n\in\Gamma$.

 The functor $\phi^*$ has a left adjoint functor~$\phi_!$ and
a right adjoint functor~$\phi_*$.
 The former assigns to a $\Gamma$\+graded $k$\+module $V$
the $\Gamma'$\+graded $k$\+module $V'$ constructed by taking
the direct sums of the grading components of $V$ over all
the preimages in $\Gamma$ of a given element $n'\in\Gamma'$,
while the latter involves taking direct products over
the preimages of~$n'$ in~$\Gamma$.

 All three functors $\phi_!$, $\phi_*$, and $\phi^*$ are exact.
 Besides, they transform (pre)complexes of $k$\+modules to
(pre)complexes of $k$\+modules and commute with passing to
the cohomology of the complexes of $k$\+modules.
 So they induce triangulated functors between the derived
categories $\DD_\Gamma(k\modl)$ and $\DD_{\Gamma'}(k\modl)$ of
$\Gamma$\+graded and $\Gamma'$\+graded complexes of $k$\+modules.

 Given a $\Gamma$\+graded $k$\+linear CDG\+category $B$, one can
apply the functor $\phi_!$ to all its precomplexes of morphisms,
obtaining a $\Gamma'$\+graded $k$\+linear CDG\+category $\phi_!B$.
 To a (left or right) CDG\+module $M'$ over $\phi_!B$ one can assign
a CDG\+module $\phi^*M'$ over $B$, and to a CDG\+module $M$ over
$B$ one can assign CDG\+modules $\phi_!M$ and $\phi_*M$
over~$\phi_!B$.

 The functors $\phi_!$, $\phi_*$, and $\phi^*$ are compatible with
the functors of tensor product and $\Hom$ of CDG\+modules in
the following sense.
 For any left CDG\+modules $L$, \ $M$ and right CDG\+module $N$
over $B$ there are natural isomorphisms
\begin{equation}  \label{phi-push-compatibility}
 \phi_!N\ot_{\phi_!B}\phi_!M\simeq\phi_!(N\ot_BM)
 \quad\text{and}\quad
 \Hom^{\phi_!B}(\phi_!L,\phi_*M)\simeq\phi_*\Hom^B(L,M).
\end{equation}
 For any left CDG\+modules $L'$ and $M'$ over $\phi_!B$ there are
natural isomorphisms
\begin{equation}  \label{phi-projection-formula}
\begin{aligned} 
 \phi^*(\phi_!N\ot_{\phi_!B}M')&\simeq N\ot_B\phi^*M', \\
 \phi^*\Hom^{\phi_!B}(\phi_!L,M')&\simeq\Hom^B(L,\phi^*M'), \\
 \phi^*\Hom^{\phi_!B}(L',\phi_*M)&\simeq\Hom^B(\phi^*L',M).
\end{aligned}
\end{equation}

 It follows from the isomorphisms~\eqref{phi-projection-formula}
that the functors $\phi_!$ preserve all the flatness and
projectivity properties of CDG- and DG-modules considered above
in this paper, while the functors $\phi_*$ preserve
the injectivity properties.
 Furthermore, the functors $\phi_!$ commute with the functors
$\Tot^\oplus$, while the functors $\phi_*$ commute with
the functors $\Tot^\sqcap$. 
 Therefore, in view of the isomorphisms~\eqref{phi-push-compatibility},
for any $\Gamma$\+graded DG\+category $C$ and DG\+modules $L$, $M$,
and $N$ over it there are natural isomorphisms
\begin{equation}
 \Tor^{\phi_!C}(\phi_!N,\phi_!M)\simeq\phi_!\Tor^C(N,M)
 \.\ \ \text{and}\ \ 
 \Ext_{\phi_!C}(\phi_!L,\phi_*M)\simeq\phi_*\Ext_C(L,M).
\end{equation}
in $\DD_{\Gamma'}(k\modl)$.

 Furthermore, the functor $\phi_!$ preserves tensor products
of $k$\+linear (C)DG\+categories.
 Thus, assuming that the complexes of morphisms in the DG\+category
$C$ are h\+flat complexes of $k$\+modules, for any DG\+module $M$
over $C\ot_k C^\op$ there are natural isomorphisms of Hochschild
(co)homology
\begin{equation}
 HH_*(\phi_!C,\phi_!M)\simeq\phi_!HH_*(C,M)
 \quad\text{and}\quad
 HH^*(\phi_!C,\phi_*M)\simeq\phi_*HH^*(C,M).
\end{equation}
 In particular, there is an isomorphism
\begin{equation}
 HH_*(\phi_!C)\simeq \phi_!HH_*(C).
\end{equation}
and a natural morphism
\begin{equation}
 HH^*(\phi_!C) \lrarrow HH^*(\phi_!C,\phi_*C)
 \.\simeq\. \phi_* HH^*(C)\.
\end{equation}
 The latter morphism is an isomorphism when the kernel of
the map $\phi\:\Gamma\rarrow\Gamma'$ is finite (so the functors
$\phi_!$ and $\phi_*$ are isomorphic).

 The analogous results for (co)homology theories of the second 
kind hold under more restrictive conditions, since the functor
$\phi_!$ does not commute with $\Tot^\sqcap$ in general, nor
does the functor $\phi_*$ commute with $\Tot^\oplus$.
 However, there are morphisms of functors $\phi_!\Tot^\sqcap
\rarrow\Tot^\sqcap\phi_!$ and $\Tot^\oplus\phi_*\rarrow
\phi_*\Tot^\oplus$.

 Hence for any $\Gamma$\+graded CDG\+category $B$ and CDG\+modules
$L$, $M$, and $N$ over it there are natural morphisms
\begin{gather}
 \label{tor-second-phi-push}
 \phi_!\Tor^B(N,M)\lrarrow\Tor^{\phi_!B}(\phi_!N,\phi_!M), \\
 \label{ext-second-phi-push}
 \Ext_{\phi_!B}(\phi_!L,\phi_*M)\lrarrow\phi_*\Ext_B(L,M).
\end{gather}
in $\DD_{\Gamma'}(k\modl)$.
 The morphisms
(\ref{tor-second-phi-push}\+-\ref{ext-second-phi-push})
are always isomorphisms when the kernel of the map
$\phi\:\Gamma\rarrow\Gamma'$ is finite.
 They are also isomorphisms when the derived functors in question
can be computed using \emph{finite} resolutions
(cf.~\ref{functors-second-kind}).
 So the morphism~\eqref{tor-second-phi-push} is an isomorphism
whenever one of the $\Gamma$\+graded $B^\#$\+modules $N^\#$ and
$M^\#$ has finite flat dimension.
 The morphism~\eqref{ext-second-phi-push} is an isomorphism
whenever either the $\Gamma$\+graded $B^\#$\+module $L^\#$ has
finite projective dimension, or the $\Gamma$\+graded
$B^\#$\+module has finite injective dimension.

 Thus, assuming that the $\Gamma$\+graded $k$\+modules of morphisms
in the category $B^\#$ are flat, for any CDG\+module $M$ over
$B\ot_k B^\op$ there are natural morphisms of Hochschild (co)homology
\begin{gather}
 \label{ho-hoch-second-phi-push}
 \phi_!HH_*^{I\!I}(B,M)\lrarrow HH_*^{I\!I}(\phi_!B,\phi_!M) \\
 \label{coho-hoch-second-phi-push}
 HH^{I\!I\;*}(\phi_!B,\phi_*M) \lrarrow\phi_*HH^{I\!I\;*}(B,M),
\end{gather}
which are always isomorphims when the kernel of the map
$\phi\:\Gamma\rarrow\Gamma'$ is finite.
 The map~\eqref{ho-hoch-second-phi-push} is an isomorphism
whenever one of the $\Gamma$\+graded $B^\#\ot_k B^\#{}^\op$\+modules
$B^\#$ and $M^\#$ has finite flat dimension.
 The map~\eqref{coho-hoch-second-phi-push} is an isomorphism
whenever either the $\Gamma$\+graded $B^\#\ot_k B^\#{}^\op$\+module
$B^\#$ has finite projective dimension, or the $\Gamma$\+graded
$B^\#\ot_k B^\#{}^\op$\+module $M^\#$ has finite injective dimension.

 In particular, there is a natural map
\begin{equation}
 \phi_!HH_*^{I\!I}(B)\lrarrow HH_*(\phi_!B),
\end{equation}
which is an isomorphism when either the kernel of the map
$\phi\:\Gamma\rarrow\Gamma'$ is finite, or the $\Gamma$\+graded
$B^\#\ot_k B^\#{}^\op$\+module $B^\#$ has finite flat dimension.
 There are also natural maps
\begin{equation}
 HH^{I\!I\;*}(\phi_!B)\lrarrow HH^{I\!I\;*}(\phi_!B,\phi_*B)
 \lrarrow \phi_*HH^{I\!I\;*}(B),
\end{equation}
which are both isomorphisms when the kernel of the map
$\phi\:\Gamma\rarrow\Gamma'$ is finite.

 Given a $\Gamma'$\+graded $k$\+linear CDG\+category $B'$, one can
apply the functor $\phi^*$ to all of its precomplexes of morphisms,
obtaining a $\Gamma$\+graded $k$\+linear CDG\+category $\phi^*B'$.
 To a (left or right) CDG\+module $M'$ over $B'$ one can assign
a CDG\+module $\phi^*M'$ over~$\phi^*B'$.
 Assume that the map $\phi\:\Gamma\rarrow\Gamma'$ is surjective.
 Then the functors $\phi^*\:Z^0(B'\modlc)\rarrow Z^0(\phi^*B'\modlc)$
and $Z^0(\modrc B')\rarrow Z^0(\modrc \phi^*B')$ are equivalences
of abelian categories.
 For any left CDG\+modules $L'$, \ $M'$ and right CDG\+module $N'$
over $B'$ there are natural isomorphisms
\begin{equation}
\begin{aligned}
 \phi^*N'\ot_{\phi^*B'}\phi^*M'&\simeq \phi^*(N'\ot_{B'}M'), \\
 \Hom^{\phi^*B'}(\phi^*L',\phi^*M')&\simeq\phi^*\Hom^{B'}(L',M').
\end{aligned}
\end{equation}
 Furthermore, the functor $\phi^*$ commutes with the functors
$\Tot^\oplus$ and $\Tot^\sqcap$ when applied to polycomplexes with
one grading by elements of the group $\Gamma'$ and the remaining 
gradings by the integers.

 Therefore, there are natural isomorphisms
\begin{align}
 \Tor^{\phi^*B'\;I\!I}(\phi^*N',\phi^*M')&\simeq
 \phi^*\Tor^{B'\;I\!I}(N',M'), \\
 \Ext_{\phi^*B'}^{I\!I}(\phi^*N',\phi^*M')&\simeq
 \phi^*\Ext_{B'}^{I\!I}(N',M')
\end{align}
in $\DD_\Gamma(k\modl)$, and similar isomorphisms for
the $\Tor$ and $\Ext$ of the first kind over a $k$\+linear
DG\+category~$C$.

 There is a natural strict CDG\+functor $\phi^*B'\ot_k\phi^*B'{}^\op
\rarrow\phi^*(B'\ot_k B'{}^\op)$.
 So, assuming that the $\Gamma'$\+graded $k$\+modules of morphisms
in the category $B'{}^\#$ are flat, for any CDG\+module $M'$ over
$B'\ot_k B'{}^\op$ there are natural maps
\begin{gather}  \label{hoch-second-phi-pull}
 HH_*^{I\!I}(\phi^*B',\phi^*M')\lrarrow \phi^*HH_*^{I\!I}(B',M'), \\
 \phi^*HH^{I\!I\;*}(B',M')\lrarrow HH^{I\!I\;*}(\phi^*B',\phi^*M'),
\end{gather}
and, in particular,
\begin{equation}  \label{hoch-second-of-itself-phi-pull}
 HH_*^{I\!I}(\phi^*B')\lrarrow \phi^*HH_*^{I\!I}(B')
 \quad\text{and}\quad
 \phi^*HH^{I\!I\;*}(B')\lrarrow HH^{I\!I\;*}(\phi^*B').
\end{equation}

 One can see that the maps~(\ref{hoch-second-phi-pull}\+-%
\ref{hoch-second-of-itself-phi-pull}) are isomorphisms whenever
the kernel $\Gamma''$ of the map $\phi\:\Gamma\rarrow\Gamma'$
is finite and its order $|\Gamma''|$ is invertible in~$k$.
 Indeed, the CDG\+category $\phi^*B'\ot_k\phi^*B'{}^\op$ is linear
over the group ring $k[\Gamma'']$ of the abelian group~$\Gamma''$, 
and the CDG\+category $\phi^*(B'\ot_K B'{}^\op)$ is strictly
equivalent (in fact, isomorphic) to $(\phi^*B'\ot_k\phi^*B'{}^\op)
\ot_{k[\Gamma'']}k$.
 The same assertions apply to Hochschild (co)homology of the first
kind of a $\Gamma'$\+graded DG\+category $C$ whose complexes of
morphisms are h\+flat complexes of $k$\+modules.

\subsection{DG-category of CDG-modules}  \label{dg-of-cdg-subsect}
 Let $B$ be a small $k$\+linear CDG\+category such that
the $\Gamma$\+graded $k$\+modules $B^\#(X,Y)$ are flat
for all objects $X$, $Y\in B$.
 Denote by $C=\modrcfp B$ the DG\+category of right CDG\+modules over
$B$, projective and finitely generated as $\Gamma$\+graded
$B^\#$\+modules, and by $D=\modrqfp B$ the CDG\+category of right
QDG\+modules over~$B$ satisfying the same condition.
 The results below also apply to finitely generated free modules
in place of finitely generated projective ones.

 There are strict $k$\+linear CDG\+functors $R\:B\rarrow D$ and
$I\:C\rarrow D$, and moreover, these CDG\+functors are
pseudo-equivalences of CDG\+categories (see~\ref{pseudo-equi-subsect}).
 Strictly speaking, the categories $C$ and $D$ as we have defined them
are only \emph{essentially} small rather than small, i.~e., they are
strictly equivalent to small CDG\+categories.
 So from now on we will tacitly assume that $C$ and $D$ have been
replaced with their small full subcategories containing at least
one object in every isomorphism class and such that the functors
$R$ and $I$ are still defined.

 The pseudo-equivalences $R$ and $I$ induce equivalences between
the DG\+categories of (left or right) CDG\+modules over
the CDG\+categories $B$, $C$, and $D$.
 Let $N$ be a right CDG\+module and $L$, $M$ be left CDG\+modules
over $B$; denote by $N_C$, $N_D$, $L_C$, $L_D$, etc.\ the corresponding
CDG\+modules over $C$ and~$D$ (which are defined uniquely up to
a unique isomorphism).
 By the results of~\ref{second-kind-general}
(see~(\ref{tor-second-kind-F-star}\+-\ref{ext-second-kind-F-star})),
the CDG\+functors $R$ and $I$ induce isomorphisms
\begin{align}
 \Tor^{B,I\!I}(N,M)\rarrow\Tor^{D,I\!I}(N_D,M_D)
 \quad&\text{and}\quad 
 \Tor^{C,I\!I}(N_C,M_C)\rarrow\Tor^{D,I\!I}(N_D,M_D) \notag \\
 \Ext^{I\!I}_D(L_D,M_D)\rarrow\Ext^{I\!I}_B(L,M)
 \ \ &\text{and}\ \
 \Ext^{I\!I}_D(L_D,M_D)\rarrow\Ext^{I\!I}_C(L_C,M_C).
\end{align}

 There are also the induced pseudo-equivalences $R\ot R^\op\:
B\ot_k B^\op\rarrow D\ot_k D^\op$ and $I\ot I^\op\:C\ot_k C^\op
\rarrow D\ot_k D^\op$.
 These pseudo-equivalences induce equivalences between
the DG\+categories of CDG\+modules over the CDG\+categories
$B\ot_k B^\op$, \ $C\ot_k C^\op$, and $D\ot_k D^\op$.
 In particular, the CDG\+module $B$ over $B\ot_k B^\op$ corresponds
to the CDG\+module $C$ over $C\ot_k C^\op$ and to the CDG\+module
$D$ over $D\ot_k D^\op$ under these equivalences of DG\+categories.
 Indeed, the closed morphisms $B\rarrow R^*D$ and $C\rarrow I^*D$
of CDG\+modules over $B\ot_k B^\op$ and $C\ot_k C^\op$ induced by
the functors $R$ and $I$ are isomorphisms, since the functors
$R^\#$ and $I^\#$ are fully faithful.

 Let $M$ be a CDG\+module over $B\ot_k B^\op$; denote by $M_C$
and $M_D$ the corresponding CDG\+modules over $C\ot_k C^\op$ and
$D\ot_k D^\op$.
 Then the CDG\+functors $R\ot R^\op$ and $I\ot I^\op$ induce
isomorphisms (see~\eqref{ho-hoch-second-kind-F-star},
\eqref{coho-hoch-second-kind-F-star})
\begin{align}
 HH^{I\!I}_*(B,M)\rarrow HH^{I\!I}_*(D,M_D)
 \quad&\text{and}\quad
 HH^{I\!I}_*(C,M_C)\rarrow HH^{I\!I}_*(D,M_D);  \\
 HH^{I\!I\;*}(D,M_D)\rarrow HH^{I\!I\;*}(B,M)
 \quad&\text{and}\quad 
 HH^{I\!I\;*}(D,M_D)\rarrow HH^{I\!I\;*}(C,M_C). \notag \\
\intertext{In particular, we obtain natural isomorphisms}
\label{hoch-B-C-isomorphisms}
 HH^{I\!I}_*(B)\simeq HH^{I\!I}_*(C)
 \quad&\text{and}\quad
 HH^{I\!I\;*}(B)\simeq HH^{I\!I\;*}(C).
\end{align}
 This is a generalization of~\cite[Theorem~2.14]{Seg}.

 When the ring $k$ has finite weak homological dimension,
any $\Gamma$\+graded complex of flat $k$\+modules is flat.
 So if the $\Gamma$\+graded $k$\+modules of morphisms in the category
$B^\#$, and hence also in the category $C^\#$, are flat, then
the complexes of morphisms in the DG\+category $C$ are h\+flat.
 Thus, both the Hochschild (co)homology of the first and the second
kind are defined for the DG\+category $C$, and therefore
the natural maps between the Hochschild (co)homology of the first
and second kind of the DG\+category $C$ with coefficients
in any DG\+module over $C\ot_k C^\op$ are defined.

\Section{Derived Categories of the Second Kind}

 In this section we interpret, under certain homological dimension
assumptions, the $\Ext$ and $\Tor$ of the second kind over
a CDG\+category in terms of the derived categories of the second
kind of CDG\+modules over it.
 This allows to obtain sufficient conditions for an isomorphism
of the Hochschild (co)homology of the first and second kind
for a DG\+category, and in particular, for the DG\+category $C$ of
CDG\+modules over a CDG\+category $B$, projective and finitely
generated as $\Gamma$\+graded $B^\#$\+modules.

\subsection{Conventional derived category}  \label{derived-first-kind}
 Given a DG\+category $D$, the additive category $H^0(D)$ has
a natural triangulated category structure provided that a zero
object and all shift and cones exist in~$D$.
 In particular, for any small DG\+category $C$ the categories
$H^0(C\modld)$ and $H^0(\modrd C)$ are triangulated.
 These are called the \emph{homotopy categories} of (left and right)
DG\+modules over $C$.

 A (left or right) DG\+module $M$ over $C$ is said to be
\emph{acyclic} if the complexes $M(X)$ are acyclic for all
objects $X\in C$, i.~e., $H(M)=0$.
 Acyclic DG\+modules form thick subcategories, closed under both
infinite directs sums and infinite products, in the homotopy
categories of DG\+modules.
 The quotient categories by these thick subcategories are called
the (conventional) derived categories (of the first kind) of
DG\+modules over $C$ and denoted by $\DD(C\modld)$ and
$\DD(\modrd C)$.

 The full subcategory of h\+projective DG\+modules
$H^0(C\modld)_\prj\subset H^0(C\modld)$ is a triangulated subcategory
whose functor to $\DD(C\modld)$ is an equivalence of
categories~\cite{Kel}, and the same applies to the full subcategory
of h\+injective DG\+modules $H^0(C\modld)_\inj\subset H^0(C\modld)$.
 To prove these results, one notices first of all that any projective
object in the exact category $Z^0(C\modld)$ is an h\+projective
DG\+mod\-ule, and similarly for injectives
(see~\ref{ext-tor-first-kind} for the discussion of this exact
category and its projective/injective objects).
 Let $P_\bu$ be a left projective resolution of a DG\+module $M$ in
$Z^0(C\modl)$; then the total DG\+module of $P_\bu$, constructed by
taking infinite direct sums along the diagonals, is an h\+projective
DG\+module quasi-isomorphic to~$M$.
 Similarly, if $J^\bu$ is a right injective resolution of a DG\+module
$M$ in $Z^0(C\modld)$, then the total DG\+module of $J^\bu$,
constructed by taking infinite products along the diagonals, is
an h\+injective DG\+module quasi-isomorphic to~$M$
\cite[Section~1]{Pkoszul}.

{\hbadness=1500
 Furthermore, the full subcategory of h\+flat DG\+modules
$H^0(C\modld)_\fl\subset H^0(C\modld)$ is a triangulated subcategory
whose quotient category by its intersection with thick subcategory of
acyclic DG\+modules is equivalent to $\DD(C\modld)$.
 This follows from the above result for h\+projective DG\+modules
and the fact that any h\+flat DG\+module is h\+projective.
 The same applies to the full subcategory of h\+flat right
DG\+modules $H^0(\modrd C)_\fl\subset H^0(\modrd C)$. \par}

 Let $k$ be a commutative ring and $C$ be a small $k$\+linear
DG\+category.
 Restricting the triangulated functor of two arguments
(see~\eqref{dg-hom}) 
$$
 \Hom^C\:H^0(C\modld)^\op\times H^0(C\modld)\lrarrow
 \DD(k\modl)
$$
to the full subcategory of h\+projective
DG\+modules in the first argument, one obtains a functor that
factors through the derived category in the second argument,
providing the derived functor
$$
 \Ext_C\:\DD(C\modld)^\op\times\DD(C\modld)\lrarrow
 \DD(k\modl).
$$
 Alternatively, restricting the functor $\Hom^C$ to the full
subcategory of h\+injective DG\+modules in the second argument,
one obtains a functor that factors through the derived category
in the first argument, leading to the same derived functor $\Ext_C$.
 The composition of this derived functor with the localization
functor $Z^0(C\modld)\rarrow\DD(C\modld)$ is isomorphic to
the derived functor $\Ext_C$ constructed in~\ref{ext-tor-first-kind}.
 For any left DG\+modules $L$ and $M$ over $C$ there is a natural
isomorphism
$$
 H^*\Ext_C(L,M)\simeq\Hom_{\DD(C\modld)}(L,M[*]).
$$

 Analogously, restricting the triangulated functor of two arguments
(see~\eqref{dg-tensor-product})
$$
 \ot_C\:H^0(\modrd C)\times H^0(C\modld)\lrarrow
 \DD(k\modl)
$$
to the full subcategory of h\+flat
DG\+modules in the first argument one obtains a functor that
factors through the Cartesian product of the derived categories,
providing the derived functor
$$
 \Tor^C\:\DD(\modrd C)\times\DD(C\modld)\lrarrow\DD(k\modl).
$$
 The same derived functor can be obtained by restricting
the functor $\ot_C$ to the full subcategory of h\+flat DG\+modules
in the second argument.
 Up to composing with the localization functors
$Z^0(\modrd C)\rarrow\DD(\modrd C)$ and $Z^0(C\modld)\rarrow
\DD(C\modld)$, this is the same derived functor $\Tor^C$ that
was constructed in~\ref{ext-tor-first-kind}.

\subsection{Derived categories of the second kind}
\label{categories-second-kind}
 Let $B$ be a small CDG\+category.
 As in~\ref{derived-first-kind}, the \emph{homotopy categories of
CDG\+modules} $H^0(B\modlc)$ and $H^0(\modrc B)$ over $B$ are
naturally triangulated.
 Given a short exact sequence $0\rarrow K'\rarrow K\rarrow K''
\rarrow 0$ in the abelian category $Z^0(B\modlc)$, one can consider
it as a finite complex of closed morphisms in the DG\+category
$B\modlc$ and take the corresponding total object in
$B\modlc$ \cite[Section~1.2]{Pkoszul}.

 A left CDG\+module over $B$ is called \emph{absolutely acyclic} if
it belongs to the minimal thick subcategory of $H^0(B\modlc)$
containing the total CDG\+modules of exact triples of CDG\+modules.
 The quotient category of $H^0(B\modlc)$ by the thick subcategory
of absolutely acyclic CDG\+modules is called the \emph{absolute
derived category} of left CDG\+modules over $B$ and denoted by
$\DD^\abs(B\modlc)$ \cite[Section~3.3]{Pkoszul}.

 A left CDG\+module over $B$ is called \emph{coacyclic} if it
belongs to the minimal triangulated subcategory of $H^0(B\modlc)$
containing the total CDG\+modules of exact triples of CDG\+modules
and closed under infinite direct sums.
 The quotient category of $H^0(B\modlc)$ by the thick subcategory
of coacyclic CDG\+modules is called the \emph{coderived category}
of left CDG\+modules over $B$ and denoted by $\DD^\co(B\modlc)$.

 The definition of a \emph{contraacyclic} CDG\+module is dual to
the previous one.
 A left CDG\+module over $B$ is called contraacyclic if it
belongs to the minimal triangulated subcategory of $H^0(B\modlc)$
containing the total CDG\+modules of exact triples of CDG\+modules
and closed under infinite products.
 The quotient category of $H^0(B\modlc)$ by the thick subcategory
of contraacyclic CDG\+modules is called the \emph{contraderived
category} of left CDG\+modules over $B$ and denoted by
$\DD^\ctr(B\modlc)$.

 Coacyclic, contraacyclic, and absolutely acyclic right CDG\+modules
are defined in the analogous way.
 The corresponding exotic derived (quotient) categories are denoted
by $\DD^\co(\modrc B)$, \ $\DD^\ctr(\modrc B)$, and
$\DD^\abs(\modrc B)$.

 We will use the similar notation $\DD^\co(C\modld)$, \
$\DD^\ctr(C\modld)$, etc.,\ in the particular case of the coderived,
contraderived, and absolutely derived categories of DG\+modules
over a small DG\+category $C$.
 Notice that any coacyclic or contraacyclic DG\+module is acyclic.
 The converse is not true~\cite[Examples~3.3]{Pkoszul}.

 Furthermore, given an exact subcategory in the abelian category
of $\Gamma$\+graded $B^\#$\+modules, one can define the class of
absolutely acyclic CDG\+modules with respect to this exact
subcategory (or the DG\+category of CDG\+modules whose underlying
$\Gamma$\+graded modules belong to this exact subcategory).
 For this purpose, one considers exact triples of CDG\+modules whose
underlying $\Gamma$\+graded modules belong to the exact subcategory,
takes their total CDG\+modules, and uses them to generate a thick
subcategory of the homotopy category of all CDG\+modules whose
underlying $\Gamma$\+graded modules belong to the exact subcategory.
 When the exact subcategory is closed under infinite direct sums
(resp.,\ infinite products), the class of coacyclic (resp.,\
contraacyclic) CDG\+modules with respect to this exact subcategory
is defined.
 Taking the quotient category, one obtains the coderived,
contraderived, or absolute derived category of CDG\+modules with
the given restriction on the underlying $\Gamma$\+graded modules.

 We will be particularly interested in the coderived and absolute
derived categories of CDG\+modules over $B$ whose underlying
$\Gamma$\+graded $B^\#$\+modules are flat or have finite
flat dimension (see~\ref{change-grading-group} for the terminology).
 Denote the DG\+categories of right CDG\+modules over $B$ with such
restrictions on the underlying $\Gamma$\+graded modules by
$\modrcfl B$ and $\modrcffd B$, and their absolute derived categories
by $\DD^\abs(\modrcfl B)$ and $\DD^\abs(\modrcffd B)$.
 The coderived category of $\modrcfl B$, defined as explained above,
is denoted by $\DD^\co(\modrcfl B)$.

 The definition of the coderived category $\DD^\co(\modrcffd B)$
requires a little more care because the class of modules of finite
flat dimension is not closed under infinite direct sums; only
the classes of modules of flat dimension not exceeding
a fixed number~$d$ are.
 Let us call a CDG\+module $N$ over $B$ \emph{d\+flat} if its
underlying $\Gamma$\+graded $B^\#$\+module $M^\#$ has flat
dimension not greater than~$d$.
 Define an object $N\in H^0(\modrcffd B)$ to be \emph{coacyclic
with respect to $\modrcffd B$} if there exists an integer $d\ge 0$
such that the CDG\+module $N$ is coacyclic with respect to
the DG\+category of $d$\+flat CDG\+modules over $B$. 
 The coderived category $\DD^\co(\modrcffd B)$ is the quotient
category of the homotopy category $H^0(\modrcffd B)$ by
the thick subcategory of CDG\+modules coacyclic with respect to
$\modrcffd B$.

 Similarly, let $B\modlc_\prj$, \ $B\modlc_\fpd$, $B\modlc_\inj$,
and $B\modlc_\fid$ denote the DG\+cate\-gories of left CDG\+modules
over $B$ whose underlying $\Gamma$\+graded $B^\#$\+modules are
projective, of finite projective dimension, injective, and
of finite injective dimension, respectively.
 The notation for the homotopy categories and exotic derived
categories of these DG\+categories is similar to the above.
 The definition of the coderived category $\DD^\co(B\modlc_\fpd)$
and the contraderived category $\DD^\ctr(B\modlc_\fid)$ involves
the same subtle point as discussed above.
 It is dealt with in the same way, i.~e., the class of CDG\+modules
coacyclic with respect to $B\modlc_\fpd$ or contraacyclic with
respect to $B\modlc_\fid$ is defined as the union of the classes
of CDG\+modules coacyclic or contraacyclic with respect to
the category of modules of the projective or injective dimension
bounded by a fixed integer.

\begin{thm}
\textup{(a)}
 The functors\/ $\DD^\co(\modrcfl B)\rarrow\DD^\co(\modrcffd B)$ 
and\/ $\DD^\abs(\modrcfl B)\allowbreak\rarrow\DD^\abs(\modrcffd B)$
induced by the embedding\/ $\modrcfl B\rarrow\modrcffd B$ are
equivalences of triangulated categories. \par
\textup{(b)}
 The functors $H^0(B\modlc_\prj)\rarrow\DD^\abs(B\modlc_\fpd)
\rarrow\DD^\co(B\modlc_\fpd)$, the first of which is induced by
the embedding\/ $B\modlc_\prj\rarrow B\modlc_\fpd$ and the second
is the localization functor, are equivalences of triangulated
categories. \par
\textup{(c)}
 The functors $H^0(B\modlc_\inj)\rarrow\DD^\abs(B\modlc_\fid)
\rarrow\DD^\ctr(B\modlc_\fid)$, the first of which is induced by
the embedding\/ $B\modlc_\inj\rarrow B\modlc_\fid$ and the second
is the localization functor, are equivalences of triangulated
categories.
\end{thm}

\begin{proof}
 The first equivalence in part~(b) is easy to prove.
 By~\cite[Theorem~3.5(b)]{Pkoszul}, CDG\+modules that are
projective as $\Gamma$\+graded modules are semiorthogonal to any
contraacyclic CDG\+modules in $H^0(B\modlc)$.
 The construction of~\cite[proof of Theorem~3.6]{Pkoszul}
shows that any object of $H^0(B\modlc_\fpd)$ is a cone
of a morphism from a CDG\+module that is absolutely acyclic
with respect to $B\modlc_\fpd$ to an object of $H^0(B\modlc_\prj)$.
 It follows that the functor $H^0(B\modlc_\prj)\rarrow
\DD^\abs(B\modlc_\fpd)$ is an equivalence of triangulated categories.
 Moreover, any object of $H^0(B\modlc_\fpd)$ that is contraacyclic
with respect to $B\modlc$ is absolutely acyclic with respect to
$B\modlc_\fpd$.

 To prove the second equivalence in part~(b), it suffices to show
that any object of $H^0(B\modlc_\prj)$ that is coacyclic with
respect to $B\modlc_\fpd$ is coacyclic with respect to
$B\modlc_\prj$ (as any object of the latter kind is clearly
contractible).
 The proof of this is analogous to the proof of part~(a) below.
 It follows that any CDG\+module coacyclic with respect to
$B\modlc_\fpd$ is absolutely acyclic with respect to
$B\modlc_\fpd$.
 The proof of part~(c) is analogous to the proof of part~(b) 
up to the duality.

 To prove part~(a), notice that the same construction
from~\cite[proof of Theorem~3.6]{Pkoszul} allows to present
any object of $H^0(\modrcffd B)$ as a cone of a morphism from
a CDG\+module that is absolutely acyclic with respect to
$\modrcffd B$ to an object of $H^0(\modrcfl B)$.
 By~\cite[Lemma~1.6]{Pkoszul}, it remains to show that any object
of $H^0(\modrcfl B)$ that is coacyclic (absolutely acyclic) with
respect to $\modrcffd B$ is coacyclic (absolutely acyclic) with
respect to $\modrcfl B$.

 We follow the idea of the proof of~\cite[Theorem~7.2.2]{Psemi}.
 Given an integer $d\ge0$, let us call a $d$\+flat right CDG\+module
$N$ over $B$ \ \emph{$d$\+coacyclic} if it is coacyclic with respect
to the exact category of $d$\+flat CDG\+modules over $B$.
 We will show that for any $d$\+coacyclic CDG\+module $N$ there
exists an $(d-1)$-coacyclic CDG\+module $L$ together with
a surjective closed morphism of CDG\+modules $L\rarrow N$ whose
kernel $K$ is also $(d-1)$-coacyclic.
 It will follow that any $(d-1)$-flat $d$\+coacyclic CDG\+module
$N$ is $(d-1)$-coacyclic, since the total CDG\+module of
the exact triple $K\rarrow L\rarrow M$ is $(d-1)$\+coacyclic, as
is the cone of the morphism $K\rarrow L$.
 By induction we will conclude that any $0$\+flat $d$\+coacyclic
CDG\+module is $0$\+coacyclic.
 The argument for absolutely acyclic CDG\+modules will be similar.

 To prove that a $d$\+coacyclic CDG\+module can be presented as
a quotient of a $(d-1)$-coacyclic CDG\+module by a $(d-1)$-coacyclic
CDG\+submodule, we will first construct such a presentation for
totalizations of exact triples of $d$\+flat CDG\+modules, and
then check that the class of $d$\+flat CDG\+modules presentable in
this form is stable under taking cones and homotopy equivalences.

\begin{lemA}
 Let $N$ be the total CDG\+module of an exact triple of
$d$\+flat CDG\+modules $N'\rarrow N''\rarrow N'''$.
 Then there exists a surjective closed morphism onto $N$ from
a $0$\+coacyclic CDG\+module $P$ with a $(d-1)$\+coacyclic
kernel~$K$.
\end{lemA}

\begin{proof}
 Choose projective objects $P'$ and $P'''$ in the abelian
category of CDG\+modules $Z^0(\modrc B)$
(see~\ref{second-kind-general}) such that there are surjective
morphisms $P'\rarrow N'$ and $P'''\rarrow N'''$.
 Then there exists a surjective morphism from the exact triple
CDG\+modules $P'\rarrow P''= P'\oplus P'''\rarrow P'''$ onto
the exact triple $N'\rarrow N''\rarrow N'''$.
 Let $K'\rarrow K''\rarrow K'''$ be the kernel of this morphism
of exact triples; then the CDG\+modules $P^{(i)}$ are $0$\+flat,
while the CDG\+modules $K^{(i)}$ are $(d-1)$-flat.
 Therefore, the total CDG\+module $P$ of the exact triple $P'\rarrow
P''\rarrow P'''$ is $0$\+coacyclic (in fact, $0$\+flat and
contractible), while the total CDG\+module $K$ of the exact triple
$K'\rarrow K''\rarrow K'''$ is $(d-1)$-coacyclic.
\end{proof}

\begin{lemB}
\textup{(a)}
 Let $K'\rarrow L'\rarrow N'$ and $K''\rarrow L''\rarrow N''$
be exact triples of CDG\+modules such that the CDG\+modules
$K'$, $L'$, $K''$, $L''$ are $(d-1)$-coacyclic, and let 
$N'\rarrow N''$ be a closed morphism of CDG\+modules.
 Then there exists an exact triple of CDG\+modules $K\rarrow L
\rarrow N$ with $N=\cone(N'\to N'')$ and\/ $(d-1)$-coacylic
CDG\+modules $K$ and~$L$. \par
\textup{(b)}
 In the situation of~\textup{(a)}, assume that the morphism
$N'\rarrow N''$ is injective with a $d$\+flat cokernel $N_0$.
 Then there exists an exact triple of CDG\+modules $K_0\rarrow
L_0\rarrow N_0$ with $(d-1)$-coacyclic CDG\+modules $K_0$
and $L_0$.
\end{lemB}

\begin{proof}
 Denote by $L'''$ the CDG\+module $L'\oplus L''$; then there is
the embedding of a direct summand $L'\rarrow L'''$ and
the surjective closed morphism of CDG\+modules $L'''\rarrow N''$
whose components are the composition $L'\rarrow N'\rarrow N''$
and the surjective morphism $L''\rarrow N''$.
 These two morphisms form a commutative square with the morphisms
$L'\rarrow N'$ and $N'\rarrow N''$.
 The kernel $K'''$ of the morphism $L'''\rarrow N''$ is the middle
term of an exact triple of CDG\+modules $K''\rarrow K'''\rarrow L'$.
 Since the CDG\+modules $K''$ and $L'$ are $(d-1)$-coacyclic,
the CDG\+module $K'''$ is $(d-1)$-coacyclic, too.
 Set $L=\cone(L'\to L''')$ and $K=\cone(K'\to K''')$.
 
 To prove part~(b), notice that the above morphisms of CDG\+modules
$L'\rarrow L'''$ and $K'\rarrow K'''$ are injective; denote their
cokernels by $L_0$ and~$K_0$.
 Then the CDG\+module $L_0\simeq L''$ is $(d-1)$-coacyclic.
 In the assumptions of part~(b), the CDG\+module $K_0$ is the kernel of
the surjective morphism $L_0\rarrow N_0$, so it is $(d-1)$-flat.
 Hence it follows from the exact triple $K'\rarrow K'''\rarrow K_0$
that $K_0$ is $(d-1)$-coacyclic.
\end{proof}

\begin{lemC}
 For any contractible $d$\+flat CDG\+module $N$ there exists
an exact triple $K\rarrow P\rarrow N$ with with a $0$\+flat
contractible CDG\+module $P$ and a $(d-1)$-flat
contractible CDG\+module~$K$.
\end{lemC}

\begin{proof}
 It is easy to see using the explicit description of projective
objects in $Z^0(\modrc B)$ given in~\ref{second-kind-general} that
any projective CDG\+module is contractible.
 Let $p\:P\rarrow N$ be a surjectuve morphism onto $N$ from
a projective CDG\+module~$P$.
 Let $t\:N\rarrow N$ be a contracting homotopy for $N$ and
$\theta\:P\rarrow P$ be a contracting homotopy for~$P$.
 Then $p\theta-tp\:P\rarrow N$ is a closed morphism of CDG\+modules
of degree~$-1$.
 Since $P$ is projective and $p$ is surjective, there exists
a closed morphism $b\:P\rarrow P$ of degree~$-1$ such that
$p\theta-tp=pb$.
 Hence $\theta-b$ is another contracting homotopy for $P$ making
a commutative square with the contracting homotopy~$t$ and
the morphism~$p$.
 It follows that the restriction of $\theta-b$ on the kernel $K$
of the morphism~$p$ is a contracting homotopy for the CDG\+module~$K$.
\end{proof}

\begin{lemD}
 Let $N\rarrow N'$ be a homotopy equivalence of $d$\+flat
CDG\+modules, and suppose that there is an exact triple of
CDG\+modules $K'\rarrow L'\rarrow N'$ with $(d-1)$-coacyclic
CDG\+modules $K'$ and~$L'$.
 Then there exists an exact triple of CDG\+modules $K\rarrow L
\rarrow N$ with $(d-1)$-coacyclic CDG\+modules $K$ and~$L$.
\end{lemD}

\begin{proof}
 The cone of the morphism $N\rarrow N'$, being a contractible
$d$\+flat CDG\+module, is the cokernel of an injective morphism
of $(d-1)$-coacyclic CDG\+modules by Lemma~C\hbox{}.
 By Lemma~B(a), the cocone $N''$ of the morphism $N'\rarrow
\cone(N\to N')$ can be also presented in such form.
 The CDG\+module $N''$ is isomorphic to the direct sum of
the CDG\+module $N$ and the cocone $N'''$ of the identity endomorphism
of the CDG\+module~$N'$.
 The CDG\+module $N'''$ can be also presented in the desired form.
 Hence, by Lemma~B(b), so can the cokernel $N$ of the injective
morphism $N'''\rarrow N''$.
\end{proof}

 It is clear that the property of a CDG\+module to be presentable
as the quotient of a $(d-1)$-coacyclic CDG\+module by
a $(d-1)$-coacyclic CDG\+submodule is stable under infinite
direct sums.
 The assertion that all $d$\+coacyclic CDG\+modules can be presented
in such form now follows from Lemmas~A, B(a), and~D.
\end{proof}

 In particular, it follows from part~(b) of Theorem that there is
a natural fully faithful functor $\DD^\co(B\modlc_\fpd)\rarrow
\DD^\ctr(B\modlc)$.
 Indeed, the functor $H^0(B\modlc_\prj)\allowbreak\simeq
\DD^\abs(B\modlc_\fpd)\rarrow\DD^\ctr(B\modlc)$ is fully faithful
by~\cite[Theorem~3.5(b) and Lemma~1.3]{Pkoszul}.
 Similarly, the functor $H^0(B\modlc_\inj)\simeq\DD^\abs(B\modlc_\fid)
\rarrow\DD^\co(B\modlc)$ is fully faithful
by~\cite[Theorem~3.5(a)]{Pkoszul}, so there is
a natural fully faithful functor $\DD^\ctr(B\modlc_\fid)\rarrow
\DD^\co(B\modlc)$.

\subsection{Derived functors of the second kind}
\label{functors-second-kind}
 Let $B$ be a small $k$\+linear CDG\+cate\-gory and $L$ be a left
CDG\+module over $B$ such that the $\Gamma$\+graded left
$B^\#$\+module $L^\#$ has finite projective dimension.
 Then the CDG\+module $L$ admits a finite left resolution $P_\bu$ in
the abelian category $Z^0(B\modlc)$ such that the $\Gamma$\+graded
$B^\#$\+modules $P_i^\#$ are projective.
 This resolution can be used to compute the functor
$\Ext_B^{I\!I}(L,{-})$ as defined in~\ref{second-kind-general}.

 On the other hand, let $P$ denote the total CDG\+module of
the finite complex of CDG\+modules $P_\bu$.
 Then for any left CDG\+module $M$ over $B$ the $k$\+module of
morphisms from $L$ into $M$ in $\DD^\abs(B\modlc)$
or $\DD^\ctr(B\modlc)$ is isomorphic to the $k$\+module of
morphisms from $P$ into $M$ in the homotopy category $H^0(B\modlc)$
\cite[Theorem~3.5(b) and Lemma~1.3]{Pkoszul}.
 Thus,
$$
 H^*\Ext_B^{I\!I}(L,M)\.\simeq\.\Hom_{\DD^\abs(B\modlc)}(L,M[*])
 \.\simeq\.\Hom_{\DD^\ctr(B\modlc)}(L,M[*]).
$$

 Similar isomorphisms
$$
 H^*\Ext_B^{I\!I}(L,M)\.\simeq\.\Hom_{\DD^\abs(B\modlc)}(L,M[*])
 \.\simeq\.\Hom_{\DD^\co(B\modlc)}(L,M[*])
$$
hold if one assumes, instead of the condition on $L^\#$, that
the $\Gamma$\+graded $B^\#$\+module $M^\#$ has finite
injective dimension.

 One can lift these comparison results from the level of cohomology
modules to the level of the derived category $\DD(k\modl)$ in
the following way.
 Consider the functor (see~\eqref{cdg-hom})
$$
 \Hom^B\:H^0(B\modlc)^\op\times H^0(B\modlc)\lrarrow\DD(k\modl)
$$
and restrict it to the subcategory $H^0(B\modlc_\prj)^\op$ in
the first argument.
 This restriction factors through the contraderived category
$\DD^\ctr(B\modlc)$ in the second argument.
 Taking into account Theorem~\ref{categories-second-kind}(b), we
obtain a right derived functor
\begin{equation} \label{ext-second-kind-proj}
 \DD^\co(B\modlc_\fpd)^\op\times\DD^\ctr(B\modlc)\lrarrow
 \DD(k\modl).
\end{equation}
 The composition of this derived functor with the localization
functors $Z^0(B\modlc_\fpd)\allowbreak\rarrow\DD^\co(B\modlc_\fpd)$
and $Z^0(B\modlc)\rarrow\DD^\ctr(B\modlc)$ agrees with the derived
functor $\Ext_B^{I\!I}$ where the former is defined.

 In the same way one can use Theorem~\ref{categories-second-kind}(c)
to construct a right derived functor
\begin{equation} \label{ext-second-kind-inj}
 \DD^\co(B\modlc)^\op\times\DD^\ctr(B\modlc_\fid)\lrarrow\DD(k\modl),
\end{equation}
which agrees with the functor $\Ext_B^{I\!I}$ where the former is
defined, up to composing with the localization functors.

 Analogously, consider the functor (see~\eqref{cdg-tensor-product})
$$
 \ot_B\:H^0(\modrc B)\times H^0(B\modlc)\lrarrow\DD(k\modl)
$$
and restrict it to the subcategory $H^0(\modrcfl B)$ in the first
argument.
 This restriction factors through the Cartesian product
$\DD^\co(\modrcfl B)\times\DD^\co(B\modlc)$.
 Indeed, the tensor product of a CDG\+module that is flat as
a $\Gamma$\+graded module with a coacyclic CDG\+module is clearly
acyclic, as is the tensor product of a CDG\+module coacyclic
with respect to $\modrcfl B$ with any CDG\+module over~$B$.
 Taking into account Theorem~\ref{categories-second-kind}(a), we
obtain a left derived functor
\begin{equation} \label{tor-second-kind}
 \DD^\co(\modrcffd B)\times\DD^\co(B\modlc)\lrarrow\DD(k\modl).
\end{equation}
 Up to composing with the localization functors $Z^0(\modrcffd B)
\rarrow\DD^\co(\modrcffd B)$ and $Z^0(B\modlc)\rarrow
\DD^\co(B\modlc)$, this derived functor agrees with the derived
functor $\Tor^{B,I\!I}$ where the former is defined.
 To see this, it suffices, as above, to choose for a CDG\+module
$N\in \modrcffd B$ a \emph{finite} left resolution $Q_\bu$ in
the abelian category $Z^0(\modrc B)$ such that the $\Gamma$\+graded
$B^\#$\+modules $Q_i^\#$ are flat.

\begin{rem}
 The functor $\Tor^{B,I\!I}$ factors through the Cartesian product of
the absolute derived categories, defining a triangulated
functor of two arguments
$$
 \DD^\abs(\modrc B)\times\DD^\abs(B\modlc)\lrarrow\DD(k\modl)
$$
\cite[Section~3.12]{Pkoszul}.
 This functor agrees with the functor~\eqref{tor-second-kind} in
the sense that the composition of the former with the functor
$\DD^\abs(\modrcffd B)\rarrow\DD^\abs(\modrc B)$ in the first
argument is isomorphic to the composition of the latter with
the functors $\DD^\abs(\modrcffd B)\rarrow\DD^\co(\modrcffd B)$
and $\DD^\abs(B\modlc)\rarrow\DD^\co(B\modlc)$ in the first and
second arguments, respectively.
 Analogously, the functor $\Ext_B^{I\!I}$ descends to a triangulated
functor of two arguments
$$
 \DD^\abs(B\modlc)^\op\times\DD^\abs(B\modlc)\lrarrow\DD(k\modl),
$$
which agrees with the functors \eqref{ext-second-kind-proj}
and~\eqref{ext-second-kind-inj} in the similar sense.
\end{rem}

\subsection{Comparison of the two theories}  \label{comparison-subsect}
 Let $C$ be a small $k$\+linear DG\+category.
 Recall (see~\ref{derived-first-kind}) the notation $H^0(C\modld)_\prj$
for the homotopy category of h\+projective left DG\+modules over $C$.
 As in~\ref{categories-second-kind}, let $C\modld_\prj$ and
$H^0(C\modld_\prj)$ denote the DG\+category of left DG\+modules
over $C$ whose underlying $\Gamma$\+graded $C^\#$\+modules are
projective, and its homotopy category.
 Finally, denote by $H^0(C\modld_\prj)_\prj$ the full triangulated
subcategory in $H^0(C\modld_\prj)$ formed by the h\+projective
left DG\+modules over $C$ whose underlying $\Gamma$\+graded
$C^\#$\+modules are projective.
 The functors
$$
 H^0(C\modld_\prj)_\prj\lrarrow H^0(C\modld)_\prj\lrarrow\DD(C\modld)
$$
are equivalences of triangulated categories.
 Moreover, for any left DG\+module $L$ over $C$ there exists
a DG\+module $P\in H^0(C\modld_\prj)_\prj$ together with
a quasi-isomorphism $P\rarrow L$ of DG\+modules over~$C$
(see \cite{Kel} or~\cite[Section~1]{Pkoszul}).
 
 The equivalence of categories $H^0(C\modld_\prj)_\prj\rarrow
\DD(C\modld)$ factors as the following composition
$$
 H^0(C\modld_\prj)_\prj\lrarrow H^0(C\modld_\prj)\lrarrow
 \DD^\co(C\modld_\fpd)\lrarrow\DD(C\modld),
$$
where the middle arrow is also an equivalence of categories
(by Theorem~\ref{categories-second-kind}(b)).
 Besides, there is the localization functor $\DD^\ctr(C\modld)
\rarrow\DD(C\modld)$.
 This allows to construct a natural morphism
\begin{equation} \label{ext-derived-first-second-proj}
 \Ext_C^{I\!I}(L,M)\lrarrow\Ext_C(L,M)
\end{equation}
in $\DD(k\modl)$ for any objects $L\in\DD^\co(C\modld_\fpd)$ and
$M\in\DD^\ctr(C\modld)$.

 Specifically, for a given DG\+module $L$ choose a DG\+module
$F\in H^0(C\modld_\prj)$ and a closed morphism $F\rarrow L$ with
a cone coacyclic with respect to $C\modld_\fpd$.
 Next, for the DG\+module $F$ choose a DG\+module
$P\in H^0(C\modld_\prj)_\prj$ together with a quasi-isomorphism
$P\rarrow F$.
 Then the complex $\Hom^C(F,M)$ represents the object
$\Ext_C^{I\!I}(L,M)$, the complex $\Hom^C(P,M)$ represents
the object $\Ext_C(L,M)$, and the morphism $\Hom^C(F,M)\rarrow
\Hom^C(P,M)$, induced by the morphism $P\rarrow F$, represents
the desired morphism~\eqref{ext-derived-first-second-proj}.
 This morphism does not depend on the choices of the objects
$F$ and~$P$.

 To see that the comparison
morphism~\eqref{ext-derived-first-second-proj} coincides with
the morphism~\eqref{ext-first-second} constructed
in~\ref{second-kind-general}, choose a projective resolution
$P_\bu$ of the object $L$ in the exact category $Z^0(C\modld)$.
 Then both the resolution $P_\bu$ and its finite canonical truncation
$\tau_{\ge -d}P_\bu$ for $d$~large enough are resolutions of $L$
that can be used to compute $\Ext_C^{I\!I}(L,M)$ by the procedure
of~\ref{second-kind-general}, while the whole resolution $P_\bu$
can also be used to compute $\Ext_C(L,M)$ by the procedure
of~\ref{ext-tor-first-kind}.
 Set $F$ to be the total DG\+module of the finite complex of
DG\+modules $\tau_{\ge -d}P_\bu$ and $P$ the total DG\+module of
the complex of DG\+modules $P_\bu$, constructed by taking
infinite direct sums along the diagonals.
 Then the morphism of complexes $\Hom^C(F,M)\rarrow\Hom^C(P,M)$
represents both the morphisms \eqref{ext-first-second}
and~\eqref{ext-derived-first-second-proj} in $\DD(k\modl)$.

 Analogously, denote by $H^0(C\modld_\inj)_\inj$ the full
triangulated subcategory in $H^0(C\modld_\inj)$ formed by
the h\+injective left DG\+modules over $C$ whose underlying
$\Gamma$\+graded $C^\#$\+modules are injective.
 Here, as above, the notation $H^0(C\modld)_\inj$ for the category
of h\+injective DG\+modules comes from~\ref{derived-first-kind},
while the notation $C\modld_\inj$ and $H^0(C\modld_\inj)$ for
the categories of DG\+modules whose underlying $\Gamma$\+graded
modules are injective is similar to that
in~\ref{categories-second-kind}.
 The functors
$$
 H^0(C\modld_\inj)_\inj\lrarrow H^0(C\modld)_\inj\lrarrow
 \DD(C\modld)
$$
are equivalences of triangulated categories; moreover,
for any left DG\+module $M$ over $C$ there exists a DG\+module
$J\in H^0(C\modld_\inj)_\inj$ together with a quasi-isomorphism
$M\rarrow J$ of CDG\+modules over~$C$.
 
 The equivalence of categories $H^0(C\modld_\inj)_\inj\rarrow
\DD(C\modld)$ factors as the following composition
$$
 H^0(C\modld_\inj)_\inj\lrarrow H^0(C\modld_\inj)\lrarrow
 \DD^\ctr(C\modld_\fid)\lrarrow\DD(C\modld),
$$
where the middle arrow is also an equivalence of categories
(by Theorem~\ref{categories-second-kind}(c)).
 Besides, there is the localization functor $\DD^\co(C\modld)\rarrow
\DD(C\modld)$.
 This allows to construct a natural morphism
\begin{equation} \label{ext-derived-first-second-inj}
 \Ext_C^{I\!I}(L,M)\lrarrow\Ext_C(L,M)
\end{equation}
in $\DD(k\modl)$ for any objects $L\in\DD^\co(C\modld)$ and
$M\in\DD^\ctr(C\modld_\fid)$.

 Specifically, for a given DG\+module $M$ choose a DG\+module
$I\in H^0(C\modld_\inj)$ and a closed morphism $M\rarrow I$ with
a cone contraacyclic with respect to $C\modld_\fid$.
 Next, for the DG\+module $I$ choose a DG\+module $J\in
H^0(C\modld_\inj)_\inj$ together with a quasi-isomorphism
$I\rarrow J$.
 Then the complex $\Hom^C(L,I)$ represents the object
$\Ext_C^{I\!I}(L,M)$, the complex $\Hom_C(L,J)$ represents
the object $\Ext_C(L,M)$, and the morphism $\Hom_C(L,I)\rarrow
\Hom_C(L,J)$ represents the desired
morphism~\eqref{ext-derived-first-second-inj}.
 This comparison morphism agrees with the comparison
morphism~\eqref{ext-first-second} from~\ref{second-kind-general}
where the former is defined.

 Finally, denote by $H^0(\modrdfl C)_\fl$ the full triangulated
subcategory in $H^0(\modrdfl C)$ formed by h\+flat right
DG\+modules over $C$ whose underlying $\Gamma$\+graded
$C^\#$\+mod\-ules are flat.
 As above, $H^0(\modrd C)_\fl$ is the homotopy category
of h\+flat right DG\+mod\-ules over $C$, while $\modrdfl C$ and
$H^0(\modrdfl C)$ denote the DG\+category of right DG\+modules
whose underlying $\Gamma$\+graded $C^\#$\+modules are flat,
and its homotopy category.

 The functors between the quotient categories of
$H^0(\modrdfl C)_\fl$ and $H^0(\modrd C)_\fl$ by their intersections
with the thick subcategory of acyclic DG\+mod\-ules and the derived
category $\DD(\modrd C)$ are equivalences of triangulated categories.
 Moreover, for any right DG\+module $N$ over $C$ there exists
a DG\+module $Q\in H^0(\modrdfl C)_\fl$ together with
a quasi-isomorphism of DG\+modules $Q\rarrow N$
\cite[Section~1.6]{Pkoszul}.

 The localization functor $H^0(\modrdfl C)_\fl\rarrow\DD(\modrd C)$
factors into the composition
$$
 H^0(\modrdfl C)_\fl\lrarrow H^0(\modrdfl C)\lrarrow
 \DD^\co(\modrdffd C)\lrarrow\DD(\modrd C).
$$
(the middle arrow being described by
Theorem~\ref{categories-second-kind}(a)).
 There is also the localization functor $\DD^\co(C\modld)\rarrow
\DD(C\modld)$.
 This allows to construct a natural morphism
\begin{equation} \label{tor-derived-first-second}
 \Tor^C(N,M)\lrarrow\Tor^{C,I\!I}(N,M)
\end{equation}
in $\DD(k\modl)$ for any objects $N\in\DD^\co(\modrdffd C)$ and
$M\in\DD^\co(C\modld)$ in the same way as above.

 Specifically, for a given DG\+module $N$ choose a DG\+module
$F\in H^0(\modrdfl C)$ and a closed morphism $F\rarrow N$ with
a cone coacyclic with respect to $\modrdffd C$.
 Next, for the DG\+module $F$ choose a DG\+module $Q\in
H^0(\modrdfl C)_\fl$ together with a quasi-isomorphism $Q\rarrow F$.
 Then the complex $F\ot_C M$ represents the object
$\Tor^{C,I\!I}(N,M)$, the complex $Q\ot_C M$ represents the object
$\Tor^C(N,M)$, and the morphism $Q\ot_C M\rarrow F\ot_C M$
represents the desired morphism~\eqref{tor-derived-first-second}.
 This comparison morphism agrees with
the morphism~\eqref{tor-first-second} from~\ref{second-kind-general}
where the former is defined.

\begin{prop}
\textup{(a)}
 The natural morphism $\Tor^C(N,M)\rarrow\Tor^{C,I\!I}(N,M)$
is an isomorphism whenever the $\Gamma$\+graded $C^\#$\+module
$N^\#$ has finite flat dimension and there exists a closed morphism
$Q\rarrow N$ into $N$ from a DG\+module $Q\in H^0(\modrdfl C)_\fl$
with a cone that is coacyclic with respect to $\modrdffd C$. \par
\textup{(b)}
 The natural morphism $\Ext_C^{I\!I}(L,M)\rarrow\Ext_C(L,M)$ is
an isomorphism whenever the $\Gamma$\+graded $C^\#$\+module $L^\#$
has finite projective dimension and the object
$L\in H^0(C\modld_\fpd)$ belongs to the triangulated subcategory
generated by $H^0(C\modld_\prj)_\prj$ and the subcategory of
objects coacyclic with respect to $C\modld_\fpd$. \par
 Equivalently, the latter conclusion holds whenever
the $\Gamma$\+graded $C^\#$\+module $L^\#$ has finite projective
dimension and the object $L\in\DD^\ctr(C\modld)$ belongs to
the image of the functor $H^0(C\modld)_\prj\rarrow\DD^\ctr(C\modld)$.
\par\textup{(c)}
 The natural morphism $\Ext_C^{I\!I}(L,M)\rarrow\Ext_C(L,M)$ is also
an isomorphism if the $\Gamma$\+graded $C^\#$\+module $M^\#$ has
finite injective dimension and the object $M\in H^0(C\modld_\fid)$
belongs to the triangulated subcategory generated by
$H^0(C\modld_\inj)_\inj$ and the subcategory of objects
contraacyclic with respect to $C\modld_\fid$. \par
 Equivalently, the latter conclusion holds whenever
the $\Gamma$\+graded $C^\#$\+module $M^\#$ has finite injective
dimension and the object $M\in\DD^\co(C\modld)$ belongs to
the image of the functor $H^0(C\modld)_\inj\rarrow\DD^\co(C\modld)$.
\end{prop}

 Notice that the equivalence of categories $H^0(C\modld)_\prj\simeq
\DD(C\modld)$ identifies the functor $H^0(C\modld)_\prj\rarrow
\DD^\ctr(C\modld)$ with the functor left adjoint to the localization
functor $\DD^\ctr(C\modld)\rarrow\DD(C\modld)$.
 Analogously, the equivalence of categories $H^0(C\modld)_\inj\simeq
\DD(C\modld)$ identifies the functor $H^0(C\modld)_\inj\rarrow
\DD^\co(C\modld)$ with the functor right adjoint to the localization
functor $\DD^\co(C\modld)\rarrow\DD(C\modld)$.

 Before we prove the proposition, let us introduce some more
notation.
 The triangulated category of (C)DG\+modules coacyclic (resp.,\
contraacyclic) with respect to a given DG\+category of (C)DG\+modules
$D$ will be denoted by $\Ac^\co(D)$ (resp.,\ $\Ac^\ctr(D)$).
 So $\Ac^\co(D)$ and $\Ac^\ctr(D)$ are triangulated subcategories
of $H^0(D)$.
 Similarly, $\Ac^\abs(D)$ denotes the triangulated subcategory of
absolutely acyclic (C)DG\+modules.
 Finally, given a DG\+category $C$, we denote by $\Ac(C\modld)$ and
$\Ac(\modrd C)$ the full subcategories of acyclic DG\+modules in
the homotopy categories $H^0(C\modld)$ and $H^0(\modrd C)$.

\begin{proof}
 Part~(a) follows immediately from the above construction of
the morphism~\eqref{tor-derived-first-second}.
 To prove the first assertion of part~(b), notice that any morphism
from an object of $H^0(C\modld_\prj)_\prj$ to an object of
$\Ac^\co(C\modld_\fpd)$ vanishes in $H^0(C\modld)$.
 In fact, any morphism from an object of $H^0(C\modld_\prj)$ to
an object of $\Ac^\co(C\modld_\fpd)$ vanishes, and any morphism from
an object of $H^0(C\modld)_\prj$ to an object of $\Ac(C\modld)$
vanishes in the homotopy category.
 By the standard properties of semiorthogonal decompositions
(see, e.~g., \cite[Lemma~1.3]{Pkoszul}), it follows that any object
$L$ in the triangulated subcategory generated by
$H^0(C\modld_\prj)_\prj$ and $\Ac^\co(C\modld_\fpd)$ in
$H^0(C\modld_\fpd)$ admits a closed morphism $P\rarrow L$ from
an object $P\in H^0(C\modld_\prj)_\prj$ with a cone in
$\Ac^\co(C\modld_\fpd)$.

 To prove the equivalence of the two conditions in part~(b), notice
that, by the same semiorthogonality lemma, a DG\+module
$L\in C\modld_\fpd$ belongs to the triangulated subcategory
generated by $H^0(C\modld_\prj)_\prj$ and $\Ac^\co(C\modld_\fpd)$
in $H^0(C\modld_\fpd)$ if and only if, as an object of
$\DD^\co(C\modld_\fpd)$, it belongs to the image of
$H^0(C\modld_\prj)_\prj$ in $\DD^\co(C\modld_\fpd)$.
 Then use the concluding remarks in~\ref{categories-second-kind}.
 Part~(c) is similar to part~(b) up to the duality.
\end{proof}

 In particular, if the left homological dimension of
the $\Gamma$\+graded category $C^\#$ is finite
(see~\ref{change-grading-group} for the terminology), then the classes
of coacyclic, contraacyclic, and absolutely acyclic left DG\+modules
over $C$ coincide~\cite[Theorem~3.6(a)]{Pkoszul}.
 In this case, for any left DG\+modules $L$ and $M$ over $C$
the morphism of $\Gamma$\+graded $k$\+modules $H^*\Ext_C^{I\!I}(L,M)
\rarrow H^*\Ext_C(L,M)$ is naturally identified with the morphism
$$
 \Hom_{\DD^\abs(C\modld)}(L,M)\lrarrow\Hom_{\DD(C\modld)}(L,M).
$$
 So if the class of absolutely acyclic left DG\+modules also coincides
with the class of acyclic DG\+modules, then the natural morphism
$\Ext_C^{I\!I}(L,M)\rarrow\Ext_C(L,M)$ is an isomorphism
for any DG\+modules $L$ and $M$ over~$C$.

 Analogously, if the weak homological dimension of
the $\Gamma$\+graded category $C^\#$ is finite and the category
$H^0(\modrdfl C)$ coincides with its full subcategory
$H^0(\modrdfl C)_\fl$, then the natural morphism
$\Tor^C(N,M)\rarrow\Tor^{C,I\!I}(N,M)$ is an isomorphism
for any DG\+modules $N$ and $M$ over~$C$.
 This follows from part~(a) of Proposition.

\subsection{Comparison for DG\+category of CDG\+modules}
\label{comparison-dg-of-cdg}
 Let $B$ be a small $k$\+linear CDG\+category and
$C=\modrcfp B$ the DG\+category of right CDG\+modules over $B$,
projective and finitely generated as $\Gamma$\+graded $B^\#$\+modules.
 The results below also apply, mutatis mutandis, to finitely generated
free modules in place of finitely generated projective ones.

 The DG\+categories of (left or right) CDG\+modules over $B$ and
DG\+modules over $C$ are naturally equivalent; let $M_C$ denote
the DG\+module over $C$ corresponding to a CDG\+module $M$ over~$B$
(see \ref{pseudo-equi-subsect} and~\ref{dg-of-cdg-subsect}).
 Denote by $B\modlc_\fp$ the DG\+category of left CDG\+modules over
$B$, finitely generated and projective as $\Gamma$\+graded
$B^\#$\+modules.
 Let $k\spcheck$ be an injective cogenerator of the abelian category
of $k$\+modules; for example, one can take $k\spcheck=k$ when
$k$~is a field, or $k\spcheck=\Hom_\Z(k,\Q/\Z)$ for any ring~$k$.

 Recall that an object $X$ of a triangulated category $T$ with
infinite direct sums is called \emph{compact} if the functor
$\Hom_T(X,{-})$ preserves infinite direct sums.
 A set of compact objects $S\subset T$ generates $T$ as
a triangulated category with infinite direct sums if and only if
the vanishing of all morphisms $X\rarrow Y[*]$ in $T$ for all
$X\in S$ implies vanishing of an object $Y\in T$
\cite[Theorem~2.1(2)]{Neem}.

\begin{thmA}
\textup{(a)} If the $\Gamma$\+graded category $B^\#$ has finite
weak homological dimension and the image of the functor
$H^0(\modrcfp B)\rarrow \DD^\co(\modrc B)$ generates
$\DD^\co(\modrc B)$ as a triangulated category with infinite direct
sums, then for any right DG\+module $N_C$ and left DG\+module $M_C$
over $C$ the natural morphism $\Tor^C(N_C,M_C)\rarrow
\Tor^{C,I\!I}(N_C,M_C)$ is an isomorphism. \par
 The same conclusion holds if the $\Gamma$\+graded category $B^\#$
has finite weak homological dimension and all objects of
$H^0(\modrcfl B)$ can be obtained from objects of $H^0(\modrcfp B)$
using the operations of shift, cone, filtered inductive limit, and
passage to a homotopy equivalent CDG\+module over~$B$. \par
\textup{(b)} If the $\Gamma$\+graded category $B^\#$ has finite
left homological dimension and the image of the functor
$H^0(B\modlc_\fp)\rarrow\DD^\co(B\modlc)$ generates
$\DD^\co(B\modlc)$ as a triangulated category with infinite direct
sums, then for any left DG\+modules $L_C$ and $M_C$ over $C$
the natural morphism $\Ext_C^{I\!I}(L_C,M_C)\rarrow
\Ext_C(L_C,M_C)$ is an isomorphism.
\end{thmA}

\begin{proof}
 The proof is based on the results of~\ref{comparison-subsect}.
 The DG\+category $B\modlc_\fp$ is equivalent to the DG\+category
$C^\op$; the equivalence assigns to a right CDG\+module $F$ the left
CDG\+module $G=\Hom^{B^\op}(F,B)$ and to a left CDG\+module $B$
the right CDG\+module $F=\Hom^B(G,B)$ over~$B$.
 Given a left CDG\+module $M$ over $B$, the corresponding left
DG\+module $M_C$ over $C$ assigns to a CDG\+module $F\in\modrcfp B$
the complex of $k$\+modules $F\ot_B M\simeq\Hom^B(G,M)$.
 Given a right CDG\+module $N$ over $B$, the corresponding right
DG\+module $N_C$ over $C$ assigns to a CDG\+module $F\in\modrcfp B$
the complex of $k$\+modules $\Hom^{B^\op}(F,N)\simeq N\ot_B G$.

 The categories of (left or right) $\Gamma$\+graded modules over
the $\Gamma$\+graded categories $B^\#$ and $C^\#$ are also equivalent.
 $\Gamma$\+graded modules corresponding to each other under these
equivalences have equal flat, projective, and injective dimensions.
 So the (weak, left, or right) homological dimensions of
the $\Gamma$\+graded categories $B^\#$ and $C^\#$ are equal.
 The equivalence between the DG\+categories of (left or right)
CDG\+modules over $B$ and DG\+modules over $C$ preserves the classes
of coacyclic, acyclic, and absolutely acyclic (C)DG\+modules.
 Given a left CDG\+module $M$ over $B$, the DG\+module $M_C$ is
acyclic if and only if the complex $F\ot_BM\simeq\Hom^B(G,M)$ is
acyclic for any CDG\+modules $F\in\modrcfp B$ and $G\in B\modlc_\fp$
(related to each other as above); similarly for a right CDG\+module
$N$ over~$B$.

 For any small CDG\+category $B$ the functor $H^0(B\modlc_\fp)
\rarrow\DD^\co(B\modlc)$ is fully faithful, and the objects in its
image are compact in the coderived category~\cite{KLN}.
 Thus, the classes of acyclic and coacyclic left DG\+modules over $C$
coincide if and only if $\DD^\co(B\modlc)$ is generated by
$H^0(B\modlc_\fp)$ as a triangulated category with infinite
direct sums.
 Now parts (a) and~(b) follow from
Proposition~\ref{comparison-subsect}(a-b); see also the concluding
remarks in~\ref{comparison-subsect}.
 For details related to the proof of the second assertion of
part~(a), see the last paragraph of the proof of Theorem~B below. 
\end{proof}

 The next, more technical result is a generalization of Theorem~A
to the case of $\Gamma$\+graded categories $B^\#$ of infinite
homological dimension.

 Let us denote by $\langle T_i\rangle_\oplus\subset T$ (resp.,\
$\langle T_i\rangle_\sqcap\subset T$) the minimal triangulated
subcategory of a triangulated category $T$ containing subcategories
$T_i$ and closed under infinite direct sums (resp.,\ infinite products).
 Given a class of CDG\+modules $E\subset Z^0(\modrc B)$, we denote by
$\langle E\rangle_\cup\subset Z^0(\modrc B)$ the full subcategory of
all CDG\+modules that can be obtained from the objects of $E$ using
the operations of shift, cone, filtered inductive limit, and passage
to a homotopy equivalent CDG\+module.

\begin{thmB}
\textup{(a)} If for a right CDG\+module $N\in H^0(\modrcffd B)$
there exist a right CDG\+module
$$
 Q\in Z^0(\modrcfl B)\cap\langle Z^0(\modrcfp B)\rangle_\cup,
$$
and a closed morphism $Q\rarrow N$ with a cone in
$\Ac^\co(\modrcffd B)$, then for any left DG\+module $M_C$ over $C$
the natural morphism $\Tor^C(N_C,M_C)\rarrow\Tor^{C,I\!I}(N_C,M_C)$
is an isomorphism. \par
\textup{(b)} If a left CDG\+module $L\in H^0(B\modlc_\fpd)$ belongs to
$$
\langle\. H^0(B\modlc_\fp)\;\Ac^\co(B\modlc_\fpd)\.\rangle_\oplus
\.\subset\. H^0(B\modlc_\fpd),
$$
then for any left DG\+module $M_C$ over $C$ the natural morphism
$\Ext_C^{I\!I}(L_C,M_C)\rarrow \Ext_C(L_C,M_C)$ is an isomorphism. \par
 Equivalently, the same conclusion hold if an object $L\in
\DD^\co(B\modlc_\fpd)$ belongs to the minimal triangulated
subcategory of\/ $\DD^\co(B\modlc_\fpd)$ containing the image of
$H^0(B\modlc_\fp)$ and closed under infinite direct sums. \par
\textup{(c)} If a left CDG\+module $M\in H^0(B\modlc_\fid)$ belongs to
$$
 \langle\.\{\Hom_k(F,k\spcheck)\},\Ac^\ctr(B\modlc_\fid)\.\rangle_
 \sqcap\.\subset\. H^0(B\modlc_\fid),
 \quad F\in H^0(\modrcfp B),
$$
then for any left DG\+module $L_C$ over $C$ the natural morphism
$\Ext_C^{I\!I}(L_C,M_C)\rarrow\Ext_C(L_C,M_C)$ is an isomorphism. \par
 Equivalently, the same conclusion holds if if an object
$M\in\DD^\ctr(B\modlc_\fid)$ belongs to the minimal triangulated
subcategory of\/ $\DD^\ctr(B\modlc)_\fid$ containing the objects
$\Hom_k(F,k\spcheck)$, where $F\in H^0(\modrcfp B)$, and closed
under infinite products.
\end{thmB}

\begin{proof}
 The parts (a-c) follow from the corresponding parts of
Proposition~\ref{comparison-subsect}.

 Indeed, a DG\+module over any small DG\+category $C$ is
h\+projec\-tive if and only if it belongs to the minimal triangulated
subcategory of $H^0(C\modld)$ containing the representable DG\+modules
and closed under infinite direct sums~\cite{Kel,Pkoszul}.
 Representable left DG\+modules over $C$ correspond to the objects of
$B\modlc_\fp$ under the equivalence between the DG\+categories
$C\modld$ and $B\modlc$ (see the proof of Theorem~A).
 It follows that the DG\+module $L_C\in H^0(C\modld_\fpd)$ belongs
to the triangulated subcategory generated by $H^0(C\modld_\fpd)_\fpd$
and the objects coacyclic with respect to $C\modld_\fpd$ if and only if
a CDG\+module $L$ over $B$ belongs to the minimal triangulated
subcategory of $H^0(B\modlc_\fpd)$ containing $H^0(B\modlc_\fp)$
and all objects coacyclic with respect to $B\modlc_\fpd$ and
closed under infinite direct sums.

 Similarly, a left DG\+module over a $k$\+linear DG\+category $C$ is
h\+injective if and only if it belongs to the minimal triangulated
subcategory of $H^0(C\modld)$ containing the DG\+modules
$\Hom_k(R_X,k\spcheck)$, where $R_X$ are the representable right
DG\+modules over~$C$.
 Representable right DG\+modules over $C$ correspond to the objects
of $\modrcfp B$ under the equivalence between the DG\+categories
$\modrd C$ nad $\modrc B$.
 So the DG\+module $M_C\in H^0(C\modld_\fid)$ belongs to
the subcategory generated by $H^0(C\modld_\fid)_\fid$ and
the objects contraacyclic with respect to $C\modld_\fid$ if and
only if a CDG\+module $M$ over $B$ belongs to the minimal
triangulated subcategory of $H^0(B\modlc_\fid)$ containing 
$\Ac^\ctr(B\modlc_\fid)$ and all CDG\+modules $\Hom_k(F,k\spcheck)$
for $F\in H^0(\modrcfp B)$, and closed inder infinite products.

 Finally, a right DG\+module over a DG\+category $C$ is h\+flat
whenever it can be obtained from the representable right DG\+modules
using the operations of shift, cone, filtered inductive limit, and
passage to a homotopy equivalent DG\+module (we do not know whether
the converse is true).
 Indeed, the class of h\+flat DG\+modules is closed under shifts,
cones, filtered inductive limits, and homotopy equivalences, since
these operations commute with the tensor product of DG\+modules
over $C$ and preserve acyclicity of complexes of $k$\+modules.
 Thus, if a right CDG\+module $Q$ over $B$ can be obtained from
objects of $\modrcfp B$ using the operations of shift, cone, filtered
inductive limit, and passage to a homotopy equivalent CDG\+module,
then the corresponding DG\+module $Q_C$ over $C$ is h\+flat.

 The equivalence of the two conditions both in~(b) and in~(c)
follows from the same semiorthogonality arguments as in
the proof of Proposition~\ref{comparison-subsect}.
\end{proof}

 Now assume that the commutative ring $k$ has finite weak homological
dimension and the $\Gamma$\+graded $k$\+modules $B^\#(X,Y)$ are
flat for all objects $X$, $Y\in B$.
 Recall that the DG\+categories of left and right CDG\+modules
over $B\ot_k B^\op$ are naturally isomorphic.
 To any left CDG\+module $L$ and right CDG\+module $N$ over $B$
one can assign the (left) CDG\+module $L\ot_k N$ over
the CDG\+category $B\ot_k B^\op$.

\begin{thmC}
\textup{(a)} If the $\Gamma$\+graded category $B^\#\ot_k B^\#{}^\op$
has finite weak homological dimension and the image of the functor
of tensor product
\begin{equation} \label{bi-cdg-tensor-product}
 \ot_k\:H^0(B\modlc_\fp)\times H^0(\modrcfp B)\lrarrow\DD^\co
 (B\ot_k B^\op\modlc)
\end{equation}
generates\/ $\DD^\co(B\ot_k B^\op\modlc)$ as a triangulated category
with infinite direct sums, then the natural map $HH_*(C,M_C)\rarrow
HH_*^{I\!I}(C,M_C)$ is an isomorphism for any DG\+module $M_C$
over the DG\+category $C\ot_k C^\op$. \par
 The same conclusion holds if the $\Gamma$\+graded category
$B^\#\ot_k B^\#{}^\op$ has finite weak homological dimension
and all objects of $H^0(B\ot_k B^\op\modlc_\fl)$ can be obtained
from objects in the image of~\textup{\eqref{bi-cdg-tensor-product}}
using the operations of shift, cone, filtered inductive limit, and
passage to a homotopy equivalent CDG\+module over $B\ot_k B^\op$. \par
\textup{(b)} If the $\Gamma$\+graded category $B^\#\ot_k B^\#{}^\op$
has finite left homological dimension and the image of
the functor~\textup{\eqref{bi-cdg-tensor-product}} generates\/
$\DD^\co(B\ot_k B^\op\modlc)$ as a triangulated category with
infinite direct sums, then the natural map
$HH^{I\!I\;*}(C,M_C)\rarrow HH^*(C,M_C)$ is an isomorphism for any
DG\+module $M_C$ over the DG\+category $C\ot_k C^\op$.
\end{thmC}

\begin{proof}
 This is a particular case of the next Theorem~D.
\end{proof}

\begin{thmD}
\textup{(a)} Suppose that the $\Gamma$\+graded
$B^\#\ot_k B^\#{}^\op$\+module $B^\#$ has finite flat dimension and
there exists a CDG\+module
$$
 Q\in Z^0(\modrcfl B\ot_k B^\op)\cap \langle\{F\ot_k G\}\rangle_\cup,
 \quad F\in H^0(\modrcfp B), \ G\in H^0(B\modlc_\fp),
$$
and a closed morphism $Q\rarrow B$ of CDG\+modules over $B\ot_k B^\op$
with a cone in $\Ac^\co(\modrcffd B\ot_k B^\op)$.
 Then the natural map $HH_*(C,M_C)\rarrow HH_*^{I\!I}(C,M_C)$ is
an isomorphism for any DG\+module $M_C$ over $C\ot_k C^\op$. \par
\textup{(b)} Suppose that the $\Gamma$\+graded
$B^\#\ot_k B^\#{}^\op$\+module $B^\#$ has finite projective dimension
and the CDG\+module $B$ over $B\ot_k B^\op$ belongs to
$$
 \langle\.\{G\ot_k F\},\Ac^\co(B\ot_k B^\op\modlc_\fpd)\.
 \rangle_\oplus,
 \quad F\in H^0(\modrcfp B), \ G\in H^0(B\modlc_\fp).
$$
 Then the natural map $HH^{I\!I\;*}(C,M_C)\rarrow HH^*(C,M_C)$ is
an isomorphism for any DG\+module $M_C$ over $C\ot_k C^\op$. \par
 Equivalently, the same conclusion holds if the the $\Gamma$\+graded
$B^\#\ot_k B^\#{}^\op$\+module $B^\#$ has finite projective dimension
and the object $B\in\DD^\co(B\ot_k B^\op\modlc_\fpd)$ belongs to
the minimal triangulated subcategory of\/
$\DD^\co(B\ot_k B^\op\modlc_\fpd)$, containing the CDG\+modules
$G\ot_kF$, where $F\in H^0(\modrcfp B)$ and\/ $G\in H^0(B\modlc_\fp)$,
and closed under infinite direct sums.
\end{thmD}

\begin{proof}
 It suffices to notice that CDG\+modules $G\ot F$ over $B\ot_k B^\op$
correspond precisely to representable DG\+modules over $C\ot_k C^\op$
under the equivalence of DG\+categories
$B\ot_k B^\op\modlc\simeq C\ot_k C^\op\modld$.
 The rest of the argument is similar to the proof of Theorem~B
and based on Proposition~\ref{comparison-subsect}(a-b).
\end{proof}

\subsection{Derived tensor product functor}
\label{derived-tensor-product-subsect}
 The following discussion is relevant in connection with the role
that the external tensor products of CDG\+modules play in the above
Theorems~\ref{comparison-dg-of-cdg}.C--D.

 Let $k$ be a commutative ring of finite weak homological dimension,
and let $B'$ and $B''$ be $k$\+linear CDG\+categories such that
the $\Gamma$\+graded $k$\+modules of morphisms in the categories
$B'{}^\#$ and $B''{}^\#$ are flat.
 Consider the functor of tensor product
$$
 \ot_k\:H^0(B'\modlc)\times H^0(B''\modlc)\lrarrow
 H^0(B'\ot_k B''\modlc).
$$
 We would like to construct its left derived functor
$$
 \ot_k^\L\:\DD^\co(B'\modlc)\times\DD^\co(B''\modlc)\lrarrow
 \DD^\co(B'\ot_k B''\modlc).
$$

 Denote by $B'\modlc_\kfl$ the DG\+category of left CDG\+modules
$M'$ over $B'$ for which all the $\Gamma$\+graded $k$\+modules
$M'{}^\#(X)$ are flat, and similarly for CDG\+modules over~$B''$.
 Notice that the natural functor from the quotient category of
$H^0(B'\modlc_\kfl)$ by its intersection with $\Ac^\co(B'\modlc)$
to the coderived category $\DD^\co(B'\modlc)$ is an equivalence
of triangulated categories.
 Indeed, the construction of~\cite[proof of Theorem~3.6]{Pkoszul}
shows that for any left CDG\+module $M'$ over $B'$ there exists
a closed morphism $F'\rarrow M'$, where $F'\in H^0(B'\modlc_\kfl)$,
with a coacyclic cone.
 So it remains to use~\cite[Lemma~3.6]{Pkoszul}.

 Restrict the above functor $\ot_k$ to the subcategory
$H^0(B'\modlc_\kfl)$ in the first argument.
 Clearly, this restriction factors through the coderived category
$\DD^\co(B''\modlc)$ in the second argument.
 Let us show that it also factors through the coderived category
$\DD^\co(B'\modlc)$ in the first argument
(cf.~\cite[Lemma~2.7]{Psemi}).
 Indeed, let $M'$ be an object of $H^0(B'\modlc_\kfl)\cap
\Ac^\co(B'\modlc)$ and $M''$ be a left CDG\+module over~$B''$.
 Choose a CDG\+module $F''\in H^0(B''\modlc_\kfl)$ such that
there is a closed morphism $F''\rarrow M''$ with a coacyclic cone.
 Then the CDG\+module $M'\ot_k F''$ is coacyclic, since $M'$ is
coacyclic and $F''$ is $k$\+flat; at the same time, the cone
of the morphism $M'\ot_k F''\rarrow M'\ot_k M''$ is coacyclic,
the cone of the morphism $F''\rarrow M''$ is coacyclic and
$M'$ is $k$\+flat.
 Thus, the CDG\+module $M'\ot_k M''$ is also coacyclic.

 We have constructed the desired derived functor $\ot_k^L$.
 Clearly, the same derived functor can be obtained by
restricting the functor $\ot_k$ to the subcategory
$H^0(B''\modlc_\kfl)$ in the second argument.

 Analogously, one can construct a derived functor
$$
 \ot_k^\L\:\DD^\co(B'\modlc_\fpd)\times\DD^\co(B''\modlc_\fpd)
 \lrarrow\DD^\co(B'\ot_k B''\modlc_\fpd),
$$
or the similar functor with modules of finite projective
dimension replaced by those of finite flat dimension.
 All one has to do is to restrict the functor $\ot_k$ to
the homotopy category of CDG\+modules whose underlying
$\Gamma$\+graded modules satisfy both conditions of $k$\+flat\-ness
and finiteness of the projective dimension over $B'$ or $B''$.

 In these situations one does not even need the condition that
the weak homological dimension of~$k$ is finite.
 However, one has to use the fact that the tensor product over
$k$ preserves finitness of projective/flat dimensions, provided
that at least one of the $\Gamma$\+graded modules being
multiplied is $k$\+flat.

\Section{Examples}

 The purpose of this section is mainly to illustrate the results
of Section~3.
 Examples of DG\+categories $C$ for which the two kinds of Hochschild
(co)homology are known to coincide are exhibited
in~\ref{zero-differentials}\+-\ref{dga-koszul}.
 Examples of CDG\+algebras $B$ such that the two kinds of Hochschild
(co)homology can be shown to coincide for the DG\+category of
CDG\+modules $C=\modrcfp B$ are considered
in~\ref{cdg-koszul}\+-\ref{matrix-factorizations}.
 Counterexamples are discussed in~\ref{counterexample}
and~\ref{direct-sum}.
 Hochschild (co)homology of matrix factorizations are considered
in~\ref{matrix-factorizations}\+-\ref{direct-sum}.

\subsection{DG\+category with zero differentials}
\label{zero-differentials}
 Let $C$ be a small $k$\+linear DG\+category such that
the differentials in the complexes $C(X,Y)$ vanish for all
objects $X$, $Y\in C$.

\begin{prop}
\textup{(a)}
 If $N$ is a right DG\+module over $C$ such that
the differentials in the complexes $N(X)$ vanish for all objects
$X\in C$ and the $\Gamma$\+graded $C^\#$\+module $N^\#$ has
finite flat dimension, then the natural morphism
$\Tor^C(N,M)\rarrow\Tor^{C,I\!I}(N,M)$ is an isomorphism for
any left DG\+module $M$ over~$C$. \par
\textup{(b)}
 If $L$ be a left DG\+module over $C$ such that the differentials
in the complexes $L(X)$ vanish for all objects $X\in C$ and
the $\Gamma$\+graded $C^\#$\+module $L^\#$ has finite projective
dimension, then the natural morphism $\Ext_C^{I\!I}(L,M)
\rarrow\Ext_C(L,M)$ is an isomorphism for any left DG\+module $M$
over~$C$. \par
\textup{(c)}
 If $M$ is a left DG\+module over $C$ such that
the differentials in the complexes $M(X)$ vanish for all objects
$X\in C$ and the $\Gamma$\+graded $C^\#$\+module $M^\#$ has
finite injective dimension, then the natural morphism
$\Ext_C^{I\!I}(L,M)\rarrow\Ext_C(L,M)$ is an isomorphism for
any left DG\+module $L$ over~$C$.
\end{prop}

\begin{proof}
 To prove part~(a), notice that a finite flat left resolution
$P_\bu$ of the $\Gamma$\+graded $C^\#$\+module $N^\#$, with every
term of it endowed with a zero differential, can be used to compute
both kinds of derived functor $\Tor$ that we are interested in.
 The proofs of parts (b) and (c) are similar.
\end{proof}

\begin{corA}
\textup{(a)} If the $\Gamma$\+graded category $C^\#$ has finite
weak homological dimension, then the natural morphism
$\Tor^C(N,M)\rarrow\Tor^{C,I\!I}(N,M)$ is an isomorphism for any
DG\+modules $N$ and~$M$. \par
\textup{(b)} If the $\Gamma$\+graded category $C^\#$ has finite
left homological dimension, then the natural morphism
$\Ext_C^{I\!I}(L,M)\rarrow\Ext_C(L,M)$ is an isomorphism for any
left DG\+modules $L$ and $M$ over $C$.
\end{corA}

\begin{proof}
 Any DG\+module over a DG\+category with vanishing differentials
is an extension of two DG\+modules with vanishing differentials.
 Indeed, the kernel and image of the differential~$d$ on such
a DG\+module is a DG\+submodule.
 So it remains to use the fact that both kinds of functors $\Ext$
and $\Tor$ assign distinguished triangles to short exact sequences
of DG\+modules in any argument, together with the preceding
proposition.
 Part (b) also follows from the fact that the classes of acyclic
and absolutely acyclic left DG\+modules over $C$ coincide in
its assumptions; see~\cite{KLN}.
\end{proof}

\begin{corB}
 Let $C$ be a DG\+category such that the complexes
$C(X,Y)$ are complexes of flat $k$\+modules with zero differentials
for all objects $X$, $Y\in C$. \par
\textup{(a)} If the $\Gamma$\+graded $C^\#\ot_k C^\#{}^\op$\+module
$C^\#$ has finite flat dimension, then the natural morphism of
Hochschild homology $HH_*(C,M)\rarrow HH_*^{I\!I}(C,M)$
is an isomorphism for any DG\+module $M$ over $C\ot_k C^\op$. \par
\textup{(b)} If the $\Gamma$\+graded $C^\#\ot_k C^\#{}^\op$\+module
$C^\#$ has finite projective dimension, then the natural morphism of
Hochschild cohomology $HH^{I\!I\;*}(C,M)\rarrow HH^*(C,M)$
is an isomorphism for any DG\+module $M$ over $C\ot_k C^\op$.
\end{corB}

\begin{proof}
 This follows directly from Proposition.
\end{proof}

\subsection{Nonpositive DG\+category} \label{nonpositive-subsect}
 Assume that our grading group $\Gamma$ is isomorphic to $\Z$ and
the isomorphism identifies $\boldsymbol{1}$ with~$1$
(see~\ref{grading-group}).

 Let $C$ be a small $k$\+linear DG\+category.
 Assume that the complexes of $k$\+modules $C(X,Y)$ are concentrated
in nonpositive degrees for all objects $X$, $Y\in C$.
 Let us call a (left or right) DG\+module $M$ over $C$ \emph{bounded
above} if all the complexes of $k$\+modules $M(X)$ are bounded
above uniformly, i.~e., there exists an integer~$n$ such that
the complexes $M(X)$ are concentrated in the degree~$\le n$ for
all~$X$.
 DG\+modules \emph{bounded below} are defined in the similar way.

\begin{propA}
\textup{(a)}
 If a right DG\+module $N$ and a left DG\+module $M$
over $C$ are bounded above, then the natural morphism $\Tor^C(N,M)
\rarrow\Tor^{C,I\!I}(N,M)$ is an isomorphism. \par
\textup{(b)}
 If a left DG\+module $L$ over $C$ is bounded above and
a left DG\+module $M$ is bounded below, then the natural morphism
$\Ext_C^{I\!I}(L,M)\rarrow\Ext_C(L,M)$ is an isomorphism.
\end{propA}

\begin{proof}
 The proof is based on the construction of the natural morphisms
(\ref{tor-first-second}\+-\ref{ext-first-second})
from~\ref{second-kind-general}.
 To prove part~(a), notice that there exists a left projective
resolution $Q_\bu$ of the DG\+module $N$ in the exact category
$Z^0(\modrd C)$ consisting of DG\+modules bounded above with
the same constant~$n$ as the DG\+module $N$, and then there is
no difference between the two kinds of totalizations of
the bicomplex $Q_\bu\ot_C M$.
 The proof of part~(b) is similar.
\end{proof}

\begin{propB}
\textup{(a)}
 If a right DG\+module $N$ over $C$ is bounded above
and the graded $C^\#$\+module $N^\#$ has finite flat dimension,
then the natural morphism $\Tor^C(N,M)\rarrow\Tor^{C,I\!I}(N,M)$
is an isomorphism for any left DG\+module $M$ over~$C$. \par
\textup{(b)}
 If a left DG\+module $L$ over $C$ is bounded above
and the graded $C^\#$\+module $L^\#$ has finite projective dimension,
then the natural morphism $\Ext_C^{I\!I}(L,M)\rarrow\Ext_C(L,M)$
is an isomorphism for any left DG\+module $M$ over~$C$. \par
\textup{(c)}
 If a left DG\+module $M$ over $C$ is bounded below
and the graded $C^\#$\+module $M^\#$ has finite injective dimension,
then the natural morphism $\Ext_C^{I\!I}(L,M)\rarrow\Ext_C(L,M)$
is an isomorphism for any left DG\+module $L$ over~$C$.
\end{propB}

\begin{proof}
 Parts~(a-c) follow from the corresponding parts of
Proposition~\ref{comparison-subsect}.

 To prove part~(b), let us choose a finite left resolution $P_\bu$
of the DG\+module $L$ in the abelian category $Z^0(C\modld)$
such that the DG\+modules $P_i$ are bounded above and their
underlying graded $C^\#$\+modules are projective.
 Then the total DG\+module $P$ of $P_\bu$ maps into $L$ with
a cone absolutely acyclic with respect to $C\modld_\fpd$, so
it suffices to show that $P$ is h\+projective.
 Indeed, any left DG\+module $P$ over $C$ that is bounded above and
projective as a graded $C^\#$\+module is h\+projective.
 To prove the latter assertion, one can construct by induction in~$n$
an increasing filtration of $P$ by DG\+submodules such that
the associated quotient DG\+modules are direct summands of direct
sums of representable DG\+modules shifted by the degree determined
by the number of the filtration component.

 The proof of part~(c) is similar up to duality, and to prove
part~(a) one has to show that a right DG\+module $Q$ over $C$ that is
bounded above and flat as a graded $C^\#$\+module is h\+flat.
 This can be done, e.~g., by using (the graded version of)
the Govorov--Lazard flat module theorem to construct a filtration
similar to the one in the projective case, except that the associated
quotient DG\+modules are filtered inductive limits of direct sums of
(appropriately shifted) representable DG\+modules.
\end{proof}

 Now assume that the complexes $C(X,Y)$ are complexes of flat
$k$\+modules concentrated in nonpositive cohomological degrees.

\begin{cor}
\textup{(a)}
 For any DG\+module $M$ over $C\ot_k C^\op$ bounded above,
the natural morphism $HH_*(C,M)\rarrow HH_*^{I\!I}(C,M)$ is
an isomorphism.
 If the graded $C^\#\ot_k C^\#{}^\op$\+module $C^\#$ has finite
flat dimension, then the latter morphism is an isomorphism for any
DG\+module~$M$. \par
\textup{(b)}
 For any DG\+module $M$ over $C\ot_k C^\op$ bounded below,
the natural morphism $HH^{I\!I\;*}(C,M)\rarrow HH^*(C,M)$ is
an isomorphism.
 If the graded $C^\#\ot_k C^\#{}^\op$\+module $C^\#$ has finite
projective dimension, then the latter morphism is an isomorphism
for any DG\+module~$M$.
\end{cor}

\begin{proof}
 Apply Propositions A and~B(a-b) to the DG\+category $C\ot_k C^\op$.
\end{proof}

 So the map $HH_*(C)\rarrow HH_*^{I\!I}(C)$ is an isomorphism under
our assumptions on~$C$.
 The map $HH^{I\!I\;*}(C)\rarrow HH^*(C)$ is an isomorphism provided
that either the DG\+module $C$ over $C\ot_k C^\op$ is bounded
below~\cite[Proposition~3.15]{CT}, or the graded
$C^\#\ot_k C^\#{}^\op$\+module $C^\#$ has finite projective
dimension.

\subsection{Strongly positive DG\+category} 
 As in~\ref{nonpositive-subsect}, we assume that the grading group
$\Gamma$ is isomorphic to $\Z$ and the isomorphism identifies
$\boldsymbol{1}$ with~$1$.

 Let $k$ be a field and $C$ be a $k$\+linear DG\+category such that
the complexes of $k$\+vector spaces $C(X,Y)$ are concentrated in
nonnegative degrees for all objects $X$, $Y\in C$, the component
$C^1(X,Y)$ vanishes for all $X$ and $Y$, the component $C^0(X,Y)$
vanishes for all nonisomorphic $X$ and $Y$, and the $k$\+algebra
$C^0(X,X)$ is semisimple for all~$X$.
 Here a noncommutative ring is called (classically) semisimple
if the abelian category of (left or right) modules over it is
semisimple.
 We keep the terminology from~\ref{nonpositive-subsect} related to
bounded DG\+modules.

\begin{propA}
\textup{(a)}
 If a right DG\+module $N$ and a left DG\+module $M$ over $C$ are
bounded below, then the natural morphism
$\Tor^C(N,M)\rarrow\Tor^{C,I\!I}(N,M)$ is an isomorphism. \par
\textup{(b)}
 If a left DG\+module $L$ over $C$ is bounded below and
a left DG\+module $M$ is bounded above, then the natural morphism
$\Ext_C^{I\!I}(L,M)\rarrow\Ext_C(L,M)$ is an isomorphism.
\end{propA}

\begin{proof}
 The proof uses the construction of the natural morphisms
(\ref{tor-first-second}--\ref{ext-first-second})
from~\ref{second-kind-general}.
 To prove part~(a), one can compute both kinds of $\Tor$
in question using the reduced bar-resolution of the DG\+module
$N$ over $C$ relative to $C^0$, i.~e.,
$$
 \dsb\lrarrow N\ot_{C^0}C/C^0\ot_{C^0}C/C^0\ot_{C^0}C
 \lrarrow N\ot_{C^0}C/C^0\ot_{C^0}C\lrarrow N\ot_{C^0}C.
$$
 Here $C^0$ is considered as a DG\+category with complexes of morphisms
concentrated in degree~$0$ and endowed with zero differentials,
$C/C_0$ is a DG\+module over $C^0\ot_k C^0{}^\op$, and $C$ is
a DG\+module over $C^0\ot_k C^\op$.
 The semisimplicity condition on $C^0$ guarantees projectivity of
right DG\+modules of the form $R\ot_{C^0}C$ as objects of the exact
category $Z^0(\modrd C)$ for all right DG\+modules $R$ over~$C^0$.
 Due to the positivity/boundedness conditions on $C$, $N$, and $M$,
there is no difference between the two kinds of totalizations of
the resulting bar-bicomplex.
 The proof of part~(b) is similar.
\end{proof}

\begin{propB}
\textup{(a)}
 If a right DG\+module $N$ over $C$ is bounded below and the graded
$C^\#$\+module $N^\#$ has finite flat dimension, then the natural
morphism $\Tor^C(N,M)\rarrow\Tor^{C,I\!I}(N,M)$ is an isomorphism
for any left DG\+module $M$ over~$C$. \par
\textup{(b)}
 If a left DG\+module $L$ over $C$ is bounded below and the graded
$C^\#$\+module $L^\#$ has finite projective dimension, then
the natural morphism $\Ext_C^{I\!I}(L,M)\rarrow\Ext_C(L,M)$ is
an isomorphism for any left DG\+module $M$ over~$C$. \par
\textup{(c)}
 If a left DG\+module $M$ over $C$ is bounded above and the graded
$C^\#$\+module $M^\#$ has finite injective dimension, then
the morphism $\Ext_C^{I\!I}(L,M)\rarrow\Ext_C(L,M)$ is
an isomorphism for any left DG\+module~$L$.
\end{propB}

\begin{proof}
 The proof is similar to that of
Proposition~\ref{nonpositive-subsect}.B\hbox{}.
 E.~g., in part~(b) the key is to show that any DG\+module over $C$
that is bounded below and projective as a graded $C^\#$\+module is
h\+projective.
 One constructs an increasing filtration similar to that
in~\ref{nonpositive-subsect} with the only difference that
the associated quotient DG\+modules are projective objects
of the exact category $Z^0(C\modld)$.
\end{proof}

\begin{cor}
\textup{(a)}
 For any DG\+module $M$ over $C\ot_k C^\op$ bounded below the natural
morphism $HH_*(C,M)\rarrow HH_*^{I\!I}(C,M)$ is an isomorphism.
 If the graded $C^\#\ot_k C^\#{}^\op$\+module $C^\#$ has finite
flat dimension, then the latter morphism is an isomorphism for
and DG\+module~$M$. \par
\textup{(b)}
 For any DG\+module $M$ over $C\ot_k C^\op$ bounded above,
the natural morphism $HH^{I\!I\;*}(C,M)\rarrow HH^*(C,M)$ is
an isomorphism.
 If the graded $C^\#\ot_k C^\#{}^\op$\+module $C^\#$ has finite
projective dimension, then the latter morphism is an isomorphism
for any DG\+module~$M$. \qed
\end{cor}

 So the map $HH_*(C)\rarrow HH_*^{I\!I}(C)$ is an isomorphism under
our assumptions on~$C$.
 The map $HH^{I\!I\;*}(C)\rarrow HH^*(C)$ is an isomorphism provided
that either the DG\+module $C$ over $C\ot_k C^\op$ is bounded
above, or the graded $C^\#\ot_k C^\#{}^\op$\+module $C^\#$ has
finite projective dimension.

\subsection{Cofibrant DG\+category}  \label{cofibrant-subsect}
 A small $k$\+linear DG\+category is called \emph{cofibrant} if
it is a retract (in the category of DG\+categories and functors
between them) of a DG\+category $k\langle x_{n,\alpha}\rangle$ of
the following form.
 As a $\Gamma$\+graded category, $k\langle x_{n,\alpha}\rangle$ is
freely generated by a set of homogeneous morphisms $x_{n,\alpha}$,
where $n$ runs over nonnegative integers and $\alpha$~belongs
to some set of indices.
 This means that the morphisms in $k\langle x_{n,\alpha}\rangle$
are the formal $k$\+linear combinations of formal compositions of
the morphisms~$x_{n,\alpha}$.
 It is additionally required that the element $dx_{n,\alpha}$
belongs to the class of morphisms multiplicatively and additively
generated by the morphisms $x_{m,\beta}$ with $m<n$.
 The cofibrant DG\+categories are exactly (up to the zero
object issue) the cofibrant objects in the model category
structure constructed in~\cite{Tab} (see also~\cite{Toen}).

 The following lemmas will be used in conjunction with
the results of~\ref{comparison-subsect} in order to prove
comparison results for the two kinds of $\Ext$, $\Tor$, and
Hochschild (co)homology for cofibrant DG\+categories.

\begin{lemA}
 Let $D$ be a DG\+category of the form
$k\langle x_{n,\alpha}\rangle$ as above. \par
\textup{(a)}
 If a right DG\+module $N$ over $D$ is such that all the complexes
of $k$\+modules $N(X)$ are h\+flat complexes of flat $k$\+modules,
then there exists a closed morphism $Q\rarrow N$, where
$Q\in H^0(\modrdfl D)_\fl$, with a cone absolutely acyclic with
respect to $\modrdffd D$. \par 
\textup{(b)}
 If a left DG\+module $L$ over $D$ is such that all the complexes
of $k$\+modules $L(X)$ are h\+projective complexes of projective
$k$\+modules, then there exists a closed morphism $P\rarrow L$,
where $L\in H^0(D\modld_\prj)_\prj$, with a cone absolutely acyclic
with respect to $D\modld_\fpd$. \par
\textup{(c)}
 If a left DG\+module $M$ over $D$ is such that all the complexes
of $k$\+modules $L(X)$ are h\+injective complexes of injective
$k$\+modules, then there is a closed morphism $M\rarrow J$,
where $J\in H^0(D\modld_\inj)_\inj$, with a cone absolutely acyclic
with respect to $D\modld_\fid$.
\end{lemA}

\begin{lemB}
\textup{(a)} If $C$ is a cofibrant $k$\+linear DG\+category and
the ring $k$ has a finite weak homological dimension, then the weak
homological dimension of the $\Gamma$\+graded category $C^\#$ is
also finite.
 Moreover, the categories $H^0(\modrdfl C)$ and $H^0(\modrdfl C)_\fl$
coincide in this case. \par
\textup{(b)} If $C$ is a cofibrant $k$\+linear DG\+category and
the ring $k$ has finite homological dimension, then the left
homological dimension of the $\Gamma$\+graded category $C^\#$
is finite.
 Moreover, the classes of acyclic and absolutely acyclic left
DG\+modules over $C$ coincide in this case.
\end{lemB}

\begin{proof}[Proof of Lemmas A and~B]
 Let us first prove parts~(b) of both lemmas.
 The following arguments generalize the proof of
\cite[Theorem~9.4]{Pkoszul} to the DG\+category case.
 For any objects $X$, $Y\in D$ denote by $V(X,Y)$ the free
$\Gamma$\+graded $k$\+module spanned by those elements
$x_{n,\alpha}$ that belong to $D(X,Y)$.
 Consider the short exact sequence of $\Gamma$\+graded
$D^\#$\+modules
$$\textstyle
 \bigoplus_{Y,Z\in D} D(X,Y)\ot_k V(Y,Z)\ot_k L(Z)\lrarrow
 \bigoplus_{Y\in D} D(X,Y)\ot_k L(Y)\lrarrow L(X).
$$
 The middle and right term are endowed with DG\+module structures,
so the left term also acquires such a structure.
 There is a natural increasing filtration on the left term induced
by the filtration of $V$ related to the indexes $n$ of
the generators $x_{n,\alpha}$.
 It is a filtration by DG\+submodules and the differentials
on the associated quotient modules are the differentials on
the tensor product induced by the differentials on the factors
$D$ and $L$ (as is the differential on the middle term).

 It follows that whenever all the complexes of $k$\+modules $L(X)$
are coacyclic (absolutely acyclic), both the middle and the left
terms of the exact sequence are coacyclic (absolutely acyclic)
DG\+modules, so $L(X)$ is also a coacyclic (absolutely acyclic)
DG\+module.
 In particular, if the homological dimension of $k$ is finite and
$L$ is acyclic, then it is absolutely acyclic.
 Furthermore, when all the complexes $L(X)$ are h\+projective
complexes of projective $k$\+modules, both the middle and
the left terms belong to $H^0(D\modld_\prj)_\prj$.
 So it suffices to take the cone of the left arrow as
the DG\+module~$P$. 

 It also follows from the same exact sequence considered as
an exact sequence of $\Gamma$\+graded $D^\#$\+modules that
the $\Gamma$\+graded $D^\#$\+module $L^\#$ has the projective
dimension at most~$1$ whenever all $L^\#(X)$ are projective
$\Gamma$\+graded $k$\+modules.
 Since for any projective $\Gamma$\+graded $D^\#$\+module $F^\#$
the $\Gamma$\+graded $k$\+modules $F^\#(X)$ are projective,
the left homological dimension of $D^\#$ can exceed
the homological dimension of~$k$ by at most~$1$.
 
 Since $\Ext$ and $\Tor$ over $\Gamma$\+graded categories are
functorial with respect to $\Gamma$\+graded functors,
the (weak, left, or right) homological dimension of a retract
$C^\#$ of a $\Gamma$\+graded category $D^\#$ does not exceed
that of~$D^\#$.
 To prove the second assertion of Lemma~B(b) for a retract $C$
of a DG\+category $D$ as above, consider DG\+functors
$I\:C\rarrow D$ and $\Pi\:D\rarrow C$ such that $\Pi I=\Id_C$.
 Let $M$ be an acyclic DG\+module over $C$; then the DG\+module
$\Pi^*M$ over $D$ is acyclic, hence absolutely acyclic, and
it follows that $M=I^*\Pi^*M$ is also absolutely acyclic.

 It remains to prove the second assertion of Lemma~B(a).
 If the underlying $\Gamma$\+graded $D^\#$\+module of a right
DG\+module $N$ over $D$ is flat, then the above exact sequence
remains exact after taking the tensor product with $N$ over~$D$.
 Besides, the $\Gamma$\+graded $k$\+modules $N^\#(X)$ are flat,
since the $\Gamma$\+graded $k$\+modules $D^\#(X,Y)$ are.
 If the weak homological dimension of~$k$ is finite, it follows
that the complexes of $k$\+modules $N(X)$ are h\+flat.
 Now if the complexes of $k$\+modules $L(X)$ are acyclic, then
the tensor products of the left and the middle terms with $N$
over $D$ are acyclic, hence the complex $N\ot_D L$ is also acyclic.

 Finally, let us deduce the same assertion for a retract $C$
of the DG\+category~$D$.
 For this purpose, notice that for any DG\+functor $F\:C\rarrow D$
the functor $F^*\:H^0(\modrd D)\rarrow H^0(\modrd C)$ has
a left adjoint functor $F_!$ given by the rule $F_!(N) =
N\ot_CD$.
 In other words, the DG\+module $F_!(N)$ assigns the complex of
$k$\+modules $N\ot_C F^*S_X$ to an object $X\in D$, where
$S_X$ is the left (covariant) representable DG\+module over $D$
corresponding to~$X$.
 The functor $F_!$ transforms objects of $H^0(\modrdfl C)$ to
objects of $H^0(\modrdfl D)$ and h\+flat DG\+modules to
h\+flat DG\+modules, since for any right DG\+module $N$ over $C$
and left DG\+module $M$ over $D$ one has $F_!N\ot_D M\simeq
N\ot_C F^*M$.
 Now if $(I,\Pi)$ is our retraction and $N\in H^0(\modrdfl C)$,
then $I_!N\in H^0(\modrdfl D)$, hence $I_!N$ is h\+flat, and
it follows that $N=\Pi_!I_! N$ is also h\+flat.
\end{proof}

\begin{corC}
 Let $C$ be a cofibrant $k$\+linear DG\+category. \par
\textup{(a)}
 Given a right DG\+module $N$ and a left DG\+module $M$ over $C$,
the natural morphism $\Tor^C(N,M)\rarrow\Tor^{C,I\!I}(N,M)$
is an isomorphism provided that either all the complexes $N(X)$, or
all the complexes $M(X)$ are h\+flat complexes of flat $k$\+modules.
 When the ring $k$ has finite weak homological dimension,
this morphism is an isomorphism for any DG\+modules $N$ and~$M$. \par  
\textup{(b)}
 Given two left DG\+modules $L$ and $M$ over $C$, the natural
morphism $\Ext_C^{I\!I}(L,M)\allowbreak\rarrow\Ext_C(L,M)$ is
an isomorphism provided that either all the complexes $L(X)$ are
h\+projective complexes of projective $k$\+modules, or all
the complexes $M(X)$ are h\+injective complexes of injective
$k$\+modules.
 When the ring $k$ has finite homological dimension,
this morphism is an isomorphism for any DG\+modules $L$ and~$M$.
\end{corC}

\begin{proof}
 Since the morphisms (\ref{tor-first-second}--\ref{ext-first-second})
are functorial with respect to DG\+functors $F\:C\rarrow D$, i.~e.,
make commutative squares with the morphisms
(\ref{tor-first-kind-F-star}--\ref{ext-first-kind-F-star}) and
(\ref{tor-second-kind-F-star}--\ref{ext-second-kind-F-star}),
it suffices to prove the statements of Corollary for a DG\+category
$D=k\langle x_{n,\alpha}\rangle$.
 Now the first assertions in both (a) and~(b) follow from
Lemma~A and Proposition~\ref{comparison-subsect}, while the second
ones follow from Lemma~B and the concluding remarks
in~\ref{comparison-subsect}.
\end{proof}

\begin{lemD}
 Let $D$ be a DG\+category of the form $k\langle x_{n,\alpha}\rangle$.
 Then the $\Gamma$\+graded $D^\#\ot_k D^\#{}^\op$\+module
$D^\#$ has projective dimension at most~$1$.
 There exists an $h$\+projective DG\+module $P$ over $D\ot_k D^\op$
and a closed morphism of DG\+modules $P\rarrow D$ with a cone
absolutely acyclic with respect to $D\ot_k D^\op\modld_\fpd$. 
\end{lemD}

\begin{proof}
 It suffices to consider the short exact sequence
\begin{multline*}
 \textstyle
 \bigoplus_{Y',Y''\in D} D(X,Y')\ot_k V(Y',Y'')\ot_k D(Y'',Z) \\
 \textstyle
 \lrarrow\bigoplus_{Y\in D} D(X,Y)\ot_k D(Y,Z)\lrarrow D(X,Z)
\end{multline*}
and argue as above.
\end{proof}

\begin{corE}
 Let $C$ be a cofibrant $k$\+linear DG\+category.
 Then for any DG\+module $M$ over $C\ot_k C^\op$, the natural
morphisms of Hochschild (co)homology $HH_*(C,M)\rarrow
HH_*^{I\!I}(C,M)$ and $HH^{I\!I\;*}(C,M)\rarrow HH^*(C,M)$
are isomorphisms.
\end{corE}

\begin{proof}
 The assertions for a DG\+category $D=k\langle x_{n,\alpha}\rangle$
follow from Lemma~D and Proposition~\ref{comparison-subsect}(a-b).
 To deduce the same results for a retract $C$ of a DG\+category $D$,
use the fact that the comparison
morphisms~\eqref{hoch-first-second-kind} make commutative squares
with the morphisms \eqref{ho-hoch-second-kind-F-star},
\eqref{coho-hoch-second-kind-F-star} and
\eqref{ho-hoch-first-kind-F-star},
\eqref{coho-hoch-first-kind-F-star}.
\end{proof}

\subsection{DG\+algebra with Koszul filtration}  \label{dga-koszul}
 Let $A$ be a DG\+algebra over a field~$k$ endowed with an increasing
filtration $F_iA$, \ $i\ge0$, such that $F_0A=k$, \ $F_iA\cdot F_jA
\subset F_{i+j}A$, and $dF_iA \subset F_{i+1}A$.
 Assume that the associated graded algebra $\gr_FA$ is Koszul
(in the grading~$i$ induced by the filtration~$F$) and has finite
homological dimension (here we use the Koszulity condition without
the assumption of finite-dimensionality of the components of $\gr_FA$,
see e.~g.~\cite{PVi}).
 Then one can assign to $A$ a coaugmented CDG\+coalgebra $\CC$ endowed
with a finite decreasing filtration~$G$ \cite[Section~6.8]{Pkoszul}
(cf.~\ref{cdg-koszul} below).

\begin{cor}
 Assume that the coaugmented coalgebra $\CC$ is conilpotent
\textup{(}see~\cite{PVi} or~\cite[Section~6.4]{Pkoszul}\textup{)}.
 Then for any right DG\+module $N$ and left DG\+module $M$
over $A$ the natural morphism $\Tor^A(N,M)\rarrow\Tor^{A,I\!I}(N,M)$
is an isomorphism.
 For any left DG\+modules $L$ and $M$ over $A$ the natural morphism
$\Ext_A^{I\!I}(L,M)\allowbreak\rarrow\Ext_A(L,M)$ is an isomorphism.
 For any DG\+module $M$ over $A\ot_k A^\op$, the natural maps
of Hochschild (co)homology $HH_*(A,M)\rarrow HH_*^{I\!I}(A,M)$ and
$HH^{I\!I\;*}(A,M)\rarrow HH^*(A,M)$ are isomorphisms.
\end{cor}

\begin{proof}
 The (left or right) homological dimension of the graded algebra
$A^\#$ is finite, since one can compute the spaces $\Ext$ over it
using the nonhomogeneous Koszul resolution.
 By~\cite[Corollary~6.8.2]{Pkoszul}, the classes of acyclic and
absolutely acyclic DG\+modules over $A$ coincide.
 Hence the first two assertions follow from the concluding remarks
in~\ref{comparison-subsect}.
 To prove the last assertion, notice that the DG\+algebra
$A\ot_k A^\op$ is endowed with the induced filtration having
the same properties as required above of the filtration on~$A$;
the corresponding CDG\+coalgebra is naturally identified
with $\CC\ot_k \CC^\op$.
 Since $\CC$ is conilpotent, so is $\CC\ot_k \CC^\op$.
 Thus, the classes of acyclic and absolutely acyclic DG\+modules
over $A\ot_k A^\op$ coincide, too.
\end{proof}

\subsection{CDG\+algebra with Koszul filtration}
\label{cdg-koszul}
 Let $B=(B,d,h)$ be a CDG\+algebra over a field~$k$ endowed with
an increasing filtration $F_iB$, \ $i\ge0$, such that $F_0B=k$, \
$F_iB\cdot F_jB \subset F_{i+j}B$, \ $dF_iB\subset F_{i+1}B$,
and $h\in F_2B$. 
 Assume that the associated graded algebra $\gr_FB$ is Koszul
and has finite homological dimension.
 Then one can assign to the filtered CDG\+algebra $(B,F)$
a CDG\+coalgebra $\CC$ endowed with a finite decreasing filtration
$G$ \cite[Section~6.8]{Pkoszul}.

 Let $C=\modrcfp B$ be the DG\+category of right CDG\+modules over
$B$, projective and finitely generated as $\Gamma$\+graded
$B^\#$\+modules.
 All the results below will also hold for finitely generated free
modules in place of finitely generated projective ones.

\begin{corA}
 For any right DG\+module $N$ and left DG\+module $M$ over $C$
the natural map $\Tor^C(N,M)\rarrow\Tor^{C,I\!I}(N,M)$ is
an isomorphism.
 For any left DG\+modules $L$ and $M$ over $C$, the natural map
$\Ext_C^{I\!I}(L,M)\rarrow\Ext_C(L,M)$ is an isomorphism.
\end{corA}

\begin{proof}
 The homological dimension of the graded algebra $B^\#$ is finite
(see~\ref{dga-koszul}).
 By \cite[Corollary~6.8.1]{Pkoszul}, the coderived category
$\DD^\co(B\modlc)$ is generated by $H^0(B\modlc_\fp)$ as
a triangulated category with infinite direct sums.
 Thus, the assertions of the corollary follow from
Theorem~\ref{comparison-dg-of-cdg}.A.
\end{proof}

 Let $\CC^\ss$ denote the maximal cosemisimple $\Gamma$\+graded
subcoalgebra of the $\Gamma$\+graded coalgebra $\CC$
\cite[Section~5.5]{Pkoszul}.
 Assume that the differential~$d$ and the curvature linear
function~$h$ on $\CC$ annihilate $C^\ss$, and the tensor product
coalgebra $\CC^\ss\ot_k\CC^{\ss\,\.\op}$ is cosemisimple.
 The latter condition always holds when the field $k$ is perfect
and the grading group $\Gamma$ contains no torsion of the order
equal to the characteristic of~$k$.

\begin{corB}
 Under the above assumptions, the natural maps of Hochschild
(co)homology $HH_*(C,M)\rarrow HH_*^{I\!I}(C,M)$ and
$HH^{I\!I\;*}(C,M)\rarrow HH^*(C,M)$ are isomorphisms for any
DG\+module $M$ over the DG\+category $C\ot_k C^\op$.
\end{corB}

\begin{proof}
 The CDG\+algebra $B\ot_k B^\op$ is endowed with the induced
filtration having the same properties; the corresponding
CDG\+coalgebra is naturally identified with $\CC\ot_k\CC^\op$.
 The coderived category of CDG\+modules $\DD^\co(B\modlc)$
is equivalent to the coderived category of CDG\+comodules
$\DD^\co(\CC\comodlc)$ \cite[Theorem~6.8]{Pkoszul}.
 This equivalence transforms the functor of tensor product
$$
 \ot_k\:\DD^\co(B\modlc)\times\DD^\co(B^\op\modlc)\lrarrow
 \DD^\co(B\ot_k B^\op\modlc)
$$
into the similar functor of tensor product
$$
 \ot_k\:\DD^\co(\CC\comodlc)\times\DD^\co(\CC^\op\comodlc)
 \lrarrow\DD^\co(\CC\ot_k\CC^\op\comodlc).
$$

 When the coalgebra $\CC^\ss\ot_k\CC^{\ss\,\.\op}$ is cosemisimple,
any DG\+comodule over it (considered as a DG\+coalgebra with zero
differential) can be obtained from tensor products of DG\+comodules
over $\CC^\ss$ and $\CC^\ss{}^\op$ using the operations of cone
and passage to a direct summand.
 The coderived category $\DD^\co(\CC\ot_k\CC^\op\comodlc)$ of
CDG\+comodules over $\CC\ot_k\CC^\op$ is generated by DG\+comodules
over $\CC^\ss\ot_k\CC^{\ss\,\.\op}$ as a triangulated category with
infinite direct sums, since the coalgebra without counit
$(\CC\ot_k\CC^\op)/(\CC^\ss\ot_k\CC^{\ss\,\.\op})$ is
conilpotent~\cite[Section~5.5]{Pkoszul}.
 Therefore, the conditions of Theorem~\ref{comparison-dg-of-cdg}.C
are satisfied for the CDG\+algebra~$B$.
\end{proof}

\subsection{Noetherian CDG\+ring} \label{noetherian-cdg-rings}
 Let $B$ be a CDG\+algebra over a commutative ring $k$ and
$C=\modrcfp B$ the DG\+category of right CDG\+modules over $B$,
projective and finitely generated as $\Gamma$\+graded $B^\#$\+modules.

\begin{corA}
 Assume that the $\Gamma$\+graded ring $B^\#$ is graded left
Noetherian and has finite left homological dimension.  Then \par
\textup{(a)} the natural map $\Tor^C(N,M)\rarrow\Tor^{C,I\!I}(N,M)$
is an isomorphism for any right DG\+module $N$ and left DG\+module $M$
over~$C$; \par
\textup{(b)} the natural map $\Ext_C^{I\!I}(L,M)\rarrow\Ext_C(L,M)$
is an isomorphism for any left DG\+modules $L$ and $M$ over~$C$.
\end{corA}

\begin{proof}
 Notice that for a left Noetherian (graded) ring the weak and left
homological dimensions coincide.
 Whenever the graded ring $B^\#$ is left Noetherian, the coderived
category $\DD^\co(B\modlc)$ is compactly generated by CDG\+modules
whose underlying $\Gamma$\+graded modules are finitely generated
(a result of D.~Arinkin, \cite[Theorem~3.11.2]{Pkoszul}).
 Assuming additionally that the left homological dimension of $B^\#$
is finite, it follows easily that $\DD^\co(B\modlc)$ is compactly
generated by $H^0(B\modlc_\fp)$.
 (See the beginning of~\ref{comparison-dg-of-cdg} for a brief
discussion of compact generation.)
 It remains to apply Theorem~\ref{comparison-dg-of-cdg}.A(a-b) to
deduce the assertions of the corollary.
\end{proof}

 Before formulating our next result, let us define yet another
exotic derived category of CDG\+modules.
 Given a small CDG\+category $D$, the \emph{complete derived category}
$\DD^\cmp(D\modlc)$ of left CDG\+modules over $D$ is the quotient
category of the homotopy category $H^0(D\modlc)$ by its minimal
triangulated subcategory, containing $\Ac^\abs(D\modlc)$ and closed
under \emph{both} infinite direct sums and products.
 CDG\+mod\-ules belonging to the latter subcategory are called
\emph{completely acyclic}.

 Now assume that the ring $k$ has finite weak homological dimension
and the $\Gamma$\+graded $k$\+module $B^\#$ is flat.
 Assume further that the $\Gamma$\+graded ring $B^\#$ is both left
and right Noetherian of finite homological dimension,
the $\Gamma$\+graded ring $B^\#\ot_k B^\#{}^\op$ is graded Noetherian
and the $\Gamma$\+graded module $B^\#$ over $B^\#\ot_k B^\#{}^\op$
has finite projective dimension.

\begin{corB}
 Suppose the CDG\+module $B$ over $B\ot_k B^\op$ belongs to
the minimal triangulated subcategory of\/ $\DD^\cmp(B\ot_k B^\op\modlc)$,
closed under infinite direct sums and containing all CDG\+modules of
the form $L\ot_k N$, where $L$ and $N$ are a left and a right
CDG\+module over $B$ and at least one of the $\Gamma$\+graded
$k$\+modules $L^\#$ and $N^\#$ is flat.
 Then for any DG\+module $M$ over $C\ot_k C^\op$ the natural maps
$HH_*(C,M)\rarrow HH_*^{I\!I}(C,M)$ and $HH^{I\!I\;*}(C,M)
\rarrow HH^*(C,M)$ are isomorphisms.
\end{corB}

\begin{proof}
 Let us check the conditions of
Theorem~\ref{comparison-dg-of-cdg}.D\hbox{}.
 In view of \cite[Theorem~3.11.2]{Pkoszul} and the discussion
in~\ref{derived-tensor-product-subsect}, the triangulated subcategory
with infinite direct sums in $\DD^\cmp(B\ot_k B^\op\modlc)$ generated
by the CDG\+modules $L\ot_k N$ with $L$ and $N$ as above coincides
with the triangulated subcategory with infinite direct sums
generated by the CDG\+modules $G\ot_k F$ with $G\in H^0(B\modlc_\fp)$
and $F\in H^0(\modrcfp B)$.
 The construction from~\cite[proof of Theorem~3.6]{Pkoszul} shows
that there exists a closed morphism from a CDG\+module
$P\in H^0(B\ot_k B^\op\modlc_\fp)$ into the CDG\+module $B$ with
the cone absolutely acyclic with respect to $B\ot_k B^\op\modlc_\fpd$.
 The triangulated subcategory with infinite direct sums generated
by $H^0(B\ot_k B^\op\modlc_\fp)$ in $H^0(B\ot_k B^\op)$ is
semiorthogonal to all completely acyclic CDG\+modules, and maps fully
faithfully to $\DD^\co(B\ot_k B^\op\modlc_\fpd)$ and to
$\DD^\cmp(B\ot_k B^\op\modlc)$ \cite{KLN}.
 So the condition that the object $P$ is generated by the objects
$G\ot_k F$ can be equivalently checked in any of these triangulated
categories.
 Notice that since the objects of $H^0(B\ot_k B^\op\modlc_\fp)$ are
compact in these triangulated categories, it does not matter whether
to generate $P$ from $G\ot_K F$ using shift, cones, and infinite
direct sums, or shift, cones, and passages to direct summands only.
\end{proof}

 One can drop the assumption that the $\Gamma$\+graded ring
$B^\#\ot_k B^\#{}^\op$ is graded Noetherian by replacing
the complete derived category $\DD^\cmp(B\ot_k B^\op\modlc)$
with the coderived category $\DD^\co(B\ot_k B^\op\modlc_\fpd)$
in the formulation of Corollary~B\hbox{}.
 Notice also that when the left homological dimension of
$B^\#\ot_k B^\#{}^\op$ is finite, all the exotic derived categories
$\DD^\cmp(B\ot_k B^\op\modlc)$, \ $\DD^\abs(B\ot_k B^\op\modlc)$, \
$\DD^\co(B\ot_k B^\op\modlc_\fpd)$, etc.\
coincide~\cite[Theorem~3.6(a)]{Pkoszul}.

\subsection{Matrix factorizations}  \label{matrix-factorizations}
 Set $\Gamma=\Z/2$.
 Let $R$ be a commutative regular local ring; suppose that $R$ is
also an algebra of essentially finite type over its residue field~$k$.
 Let $w\in R$ be a noninvertible element whose zero locus has
an isolated singularity at the closed point of the spectrum of~$R$.
 Consider the CDG\+algebra $(B,d,h)$ over~$k$, where $B$ is
the algebra $R$ placed in degree~$0$, \ $d=0$, and $h=-w$.
 Let $C=\modrcfp B$ be the corresponding DG\+category of right
CDG\+modules; its objects are conventionally called
the \emph{matrix factorizations} of~$w$.

 The computations in~\cite{Seg} and~\cite{Dyck} show that the two
kinds of Hochschild (co)ho\-mology for the $k$\+linear
DG\+category $C$ are isomorphic.
 The somewhat stronger assertion that the natural maps
$HH_*(C,M)\rarrow HH_*^{I\!I}(C,M)$ and $HH^{I\!I\;*}(C,M)
\rarrow HH^*(C,M)$ are isomorphisms for any DG\+module $M$
over $C\ot_k C^\op$ follows from our
Corollary~\ref{noetherian-cdg-rings}.B\hbox{}.
 Indeed, according to~\cite[Theorem~4.1 and the discussion in
Section~6.1]{Dyck} the assumption of the corollary is satisfied
in this case.

 More generally, let $X$ be a smooth affine variety over
a field~$k$ and $R$ be the $k$\+algebra of regular functions
on~$X$.
 Let $w\in R$ be such a function; consider the CDG\+algebra
$(B,d,h)$ constructed from $R$ and~$w$ as above.
 Let $C=\modrcfp B$ be the DG\+category of right CDG\+modules
over~$B$, projective and finitely generated as $\Gamma$\+graded
$B^\#$\+modules.

\begin{corA}
 Assume that the morphism $w\:X\setminus\{w=0\}\rarrow
\mathbb A^1_k$ from the open complement of the zero locus of~$w$
in $X$ to the affine line is smooth.
 Assume moreover that either
\begin{enumerate}
\renewcommand{\theenumi}{\alph{enumi}}
\item there exists a smooth closed subscheme $Z\subset X$ such that
$w\:X\setminus Z\rarrow\mathbb A^1_k$ is a smooth morphism and
$w|_Z=0$, or
\item the field $k$~is perfect.
\end{enumerate}
Then the natural maps $HH_*(C,M)\rarrow HH_*^{I\!I}(C,M)$ and
$HH^{I\!I\;*}(C,M)\rarrow HH^*(C,M)$ are isomorphisms for any
DG\+module $M$ over $C\ot_k C^\op$.
\end{corA}

\begin{proof}
 The proof is based on Corollary~\ref{noetherian-cdg-rings}.B\hbox{},
Orlov's theorem connecting matrix factorizations with
the triangulated categories of singularities~\cite{Or}, and some
observations from the paper~\cite{LP}.
 We will show that all objects of $H^0(B\ot_k B^\op\modlc_\fp)$ can
be obtained from the objects $G\ot_k F$ with $G\in H^0(B\modlc_\fp)$
and $F\in H^0(\modrcfp B)$ using the operations of cone and passage
to a direct summand.
 By Orlov's theorem, the triangulated categories $H^0(B\modlc_\fp)$
and $H^0(\modrcfp B)$ can be identified with the triangulated
category $\DD_\Sing^\b(X_0)$ of singularities of the zero locus
$X_0\subset X$ of the function~$w$.
 Similarly, the triangulated category $H^0(B\ot_k B^\op\modlc_\fp)$
is identified with the triangulated category $\DD_\Sing^\b(Y_0)$
of singularities of the zero locus $Y_0\subset X\times_k X$ of
the function $w\times 1-1\times w$ on the Cartesian product
$X\times_k X$.

\begin{lemB}
 The equivalences of categories $H^0(B\modlc_\fp)\simeq\DD_\Sing^\b
(X_0)\simeq H^0(\modrcfp\allowbreak B)$ and $H^0(B\ot_k B^\op\modlc_\fp)
\simeq \DD_\Sing^\b(Y_0)$ transform the external tensor product functor
$H^0(B\modlc_\fp)\times H^0(\modrcfp B)\rarrow
H^0(B\ot_k B^\op\modlc_\fp)$ into the functor
$\DD^\b_\Sing(X_0)\times\DD^\b_\Sing(X_0)\rarrow\DD^\b_\Sing(Y_0)$
induced by the composition of the external tensor product of
coherent sheaves on two copies of $X_0$ and the direct image under
the closed embedding $X_0\times_k X_0\rarrow Y_0$.
\end{lemB}

\begin{proof}
 Rather than checking the assertion of the lemma for Orlov's cokernel
functor $\Sigma\:H^0(B\modlc_\fp) \rarrow \DD_\Sing^\b(X_0)$, one can
use the construction of the inverse functor $\Upsilon\:\DD_\Sing^\b
(X_0)\rarrow H^0(B\modlc_\fp)$ given in~\cite{Porl}, for which
the desired compatibility is easy to establish.
 Alternatively, it suffices to use the result
of~\cite[Lemma~2.18]{LP}.

 Let $\DD^\abs(B\modlc_\fg)$ denote the absolute derived category
of left CDG\+modules over $B$ whose underlying $\Gamma$\+graded
$B^\#$\+modules are finitely generated; the notation
$\DD^\abs(\modrcfg B)$ for right CDG\+modules will have
the similar meaning.
 Then the external tensor product of finitely generated CDG\+modules
induces a functor $\DD^\abs(B\modlc_\fg)\times\DD^\abs(\modrcfg B)
\rarrow\DD^\abs(B\ot_k B^\op\modlc_\fg)$; furthermore, the natural
functor $H^0(B\modlc_\fp)\rarrow\DD^\abs(B\modlc_\fg)$ is
an equivalence of categories, since $B^\#$ is Noetherian of
finite homological dimension. {\hfuzz=4.5pt\par}

 Let $M\in H^0(B\modlc_\fp)$; the direct image of the corresponding
coherent sheaf $\Sigma(M)$ on $X_0$ under the closed embedding
$X_0\rarrow X$ can be viewed as an object of $\DD^\abs(B\modlc_\fg)$.
 It is clear from the above-mentioned lemma from~\cite{LP} that
this object is naturally isomorphic to the image of $M$ in
$\DD^\abs(B\modlc_\fg)$.
 Let $N\in H^0(\modrcfp B)$; then the coherent sheaf $\Sigma(N)$
on $X_0$, viewed as an object of $\DD^\abs(\modrcfg B)$, is
isomorphic to $N$.

 Since the external tensor product is well-defined on the absolute
derived categories of finitely generated CDG\+modules, it follows
that the coherent sheaf $\Sigma(M)\bt_k\Sigma(N)$ on $X_0\times_k
X_0$, viewed as an object of $\DD^\abs(B\ot_k B^\op\modlc_\fg)$,
is isomorphic to $M\ot_k N$.
 Applying the same lemma from~\cite{LP} again, we conclude that
$\Sigma(M\ot_kN)\simeq\Sigma(M)\bt_k\Sigma(N)$ in
$\DD^\b_\Sing(Y_0)$.
 The assertion of Lemma~B is proven.
\end{proof}

 Now we can finish the proof of the corollary.
 Recall that in the case~(a) we have a closed subvariety
$Z\subset X_0$; in the case~(b), let $Z=X_0$ (or any closed
subvariety of $X_0$ such that the morphism $w\:X\setminus Z
\rarrow\mathbb A^1_k$ is smooth).
 It suffices to show that the external tensor products of coherent
sheaves on two copies of $Z$, considered as coherent sheaves on
$Y_0$, generate the triangulated category of singularities of~$Y_0$.

 The open complement to $Z\times Z$ in $Y_0$ is a smooth variety.
 Indeed, we have $Y_0=X\times_{\mathbb A^1_k}X$.
 The complement to $Z\times Z$ in $Y_0$ is covered by its open
subschemes $(X\setminus Z)\times_{\mathbb A^1_k}X$ and
$X\times_{\mathbb A^1_k}(X\setminus Z)$, which are both smooth,
since $X$ is smooth over~$k$ and $X\setminus Z$ is smooth
over~$\mathbb A^1_k$.
 By~\cite[Proposition~2.7]{Or2} (see also~\cite[Theorem~3.5]{LP}
and~\cite[Theorem~2.1.5 and/or Lemma~2.6]{Neem}), it follows that
the triangulated category of singularities of $Y_0$ is generated
by coherent sheaves on $Z\times_k Z$.

 It remains to show that the derived category of coherent sheaves
on $Z\times_k Z$ is generated by the external tensor products of
coherent sheaves on the Cartesian factors.
 This is true, at least, (a)~for any smooth affine scheme $Z$ of
finite type over a field~$k$, and (b)~for any affine scheme $Z$ of
finite type over a perfect field~$k$ (the affineness assumption
can be weakened, of course).
 The case~(a) is clear, since $Z\times_k Z$ is (affine and) regular
of finite Krull dimension.
 In the case~(b), any reduced scheme of finite type over~$k$
contains an open dense smooth
subscheme~\cite[Corollaire~(17.15.13)]{EGAIV4}, and~\cite[(proof of)
Theorem~3.7]{LP} applies.

 When $X$ contains connected components on which $w$~is identically
zero, Orlov's theorem is not applicable.
 On such components, one has to consider the $\Z/2$\+graded
derived category of coherent sheaves in place of the triangulated
category of singularities of the zero locus.
 Otherwise, the above argument remains unchanged.
\end{proof}

\subsection{Trivial counterexample}  \label{counterexample}
 The two kinds of Hochschild (co)homology of DG\+algebras cannot
be always isomorphic for very general reasons.

 The Hochschild homology and cohomology of the first kind $HH_*(A)$
and $HH^*(A)$ of a DG\+algebra $A$ are invariant with respect to
quasi-isomorphisms of DG\+algebras (see~\ref{hochschild-subsect}).
 On the other hand, the Hochschild homology and cohomology of
the second kind $HH_*^{I\!I}(A)$ and $HH^{I\!I\;*}(A)$ are invariant
with respect to isomorphisms of DG\+algebras in the category
of CDG\+algebras (since the Hochschild (co)homology of the second
kind are generally functorial with respect to CDG\+functors;
see~\ref{second-kind-general}\+-\ref{hochschild-subsect}).
 These are two incompatible types of invariance properties; indeed,
any two DG\+alge\-bras over a field can be connected by a chain
of tranformations some of which are quasi-isomorphisms and
the other are CDG\+isomorphisms~\cite[Examples~9.4]{Pkoszul}. 

 Here is a rather trivial example of a CDG\+algebra $B$ over
a field~$k$ such that for the corresponding DG\+category
$C=\modrcfp B$ the two kinds of Hochschild (co)homology are
different.
 This example also shows that one cannot drop the conditions on
the differential~$d$ and curvature~$h$ of the CDG\+coalgebra $\CC$
in Corollary~\ref{cdg-koszul}.B, nor can one drop the condition
on the critical values of the potential~$w$ in
Corollary~\ref{matrix-factorizations}.A.

 Set $\Gamma=\Z/2$ and $(B,d,h)=(k,0,1)$.
 So $B$ is the $k$\+algebra $k$ placed in the grading $0\bmod 2$
and endowed with the zero differential and a nonzero curvature
element.
 Then any CDG\+module over $B$ is contractible, any one of the two
components of the differential of a CDG\+module being its
contracting homotopy.
 So the DG\+category $C$ is quasi-equivalent to zero, hence
$HH_*(C)=0=HH^*(C)$.

 On the other hand, the CDG\+algebra $B\ot_k B^\op$ is simply
the $\Z/2$\+graded $k$\+algebra $k$ with the zero differential
and curvature.
 The CDG\+module $B$ over $B\ot_k B^\op$ is the $\Z/2$\+graded
$k$\+module $k$ concentrated in degree $0\bmod 2$.
 So $\Tor^{B\ot_k B^\op}(B,B)\simeq\Tor^{B\ot_k B^\op,I\!I}(B,B)
\simeq k\simeq\Ext_{B\ot_k B^\op}^{I\!I}(B,B)\simeq
\Ext_{B\ot_k B^\op}(B,B)$.
 We conclude that $HH_*^{I\!I}(B)\simeq k\simeq HH^{I\!I\;*}(B)$
and, by the isomorphisms~\eqref{hoch-B-C-isomorphisms}
from~\ref{dg-of-cdg-subsect},
$HH_*^{I\!I}(C)\simeq k\simeq HH^{I\!I\;*}(C)$.

\subsection{Direct sum formula} \label{direct-sum}
 Set $\Gamma=\Z/2$.
 Let $k$ be a commutative ring and $B$ be a small $k$\+linear
CDG\+category such that all $\Gamma$\+graded $k$\+modules
of morphisms in $B^\#$ are $k$\+flat.
 Given a constant $c\in k$, denote by $B_{(c)}$ the $k$\+linear
CDG\+category obtained from the CDG\+category $B$ by adding~$c$
to all the curvature elements in~$B$.
 Then there is a natural (strict) isomorphism $B_{(c)}\ot_k
B_{(c)}^\op\simeq B\ot_k B^\op$, hence a natural isomorphism
between the DG\+categories of CDG\+modules over $B_{(c)}\ot_k
B_{(c)}^\op$ and $B\ot_k B^\op$.
 This isomorphism transforms the diagonal CDG\+bimodule $B_{(c)}$
over $B_{(c)}\ot_k B_{(c)}^\op$ to the diagonal CDG\+bimodule
$B$ over $B\ot_k B^\op$.

 Therefore, we have natural isomorphisms
\begin{equation}
 HH_*^{I\!I}(B_{(c)})\simeq HH_*^{I\!I}(B)
 \quad\textup{and}\quad
 HH^{I\!I\;*}(B_{(c)})\simeq HH^{I\!I\;*}(B),
\end{equation}
and consequently similar isomorphisms for the Hochschild (co)homology
of the second kind of the DG\+categories $C=\modrcfp B$ and
$C_{(c)}=\modrcfp B_{(c)}$.
 Hochschild (co)homology of the first kind of the DG\+categories
$C$ and $C_{(c)}$ are \emph{not} isomorphic in general (and in fact
can be entirely unrelated, as the following example illustrates).

 Let $k$~be an algebraically closed field of characteristic~$0$,
\ $X$ be a smooth affine variety over~$k$, and $w$~be a regular
function on~$X$.
 Let $B$ be the CDG\+algebra associated with $X$ and~$w$ as
in~\ref{matrix-factorizations}.
 The function~$w$ has a finite number of critical values $c_i\in k$.
 When $c$~is not a critical value, the Hochschild (co)homology
of the first kind $HH_*(C_{(c)})$ and $HH^*(C_{(c)})$ vanish,
since the category $H^0(C_{(c)})$ does.
 We have natural maps
$$
 HH_*(C_{(c_i)})\lrarrow HH_*^{I\!I}(C_{(c_i)})\.\simeq\.
 HH_*^{I\!I}(C)
$$
and
$$
 HH^{I\!I\;*}(C)\.\simeq\. HH^{I\!I\;*}(C_{(c_i)})\lrarrow
 HH^*(C_{(c_i)}).
$$

\begin{cor}
 The induced maps
\begin{equation} \label{ho-hoch-direct-sum}
 \textstyle\bigoplus_i HH_*(C_{(c_i})\lrarrow HH_*^{I\!I}(C)
\end{equation}
and
\begin{equation} \label{coho-hoch-direct-sum}
 HH^{I\!I\;*}(C)\lrarrow\textstyle\bigoplus_i HH^*(C_{(c_i)}).
\end{equation}
are isomorphisms.
\end{cor}

\begin{proof}
 It follows from the spectral sequence computation for
the Hochschild (co)ho\-mology of the second kind
in~\cite[proof of Theorem~4.2(b)]{CT} (see also~\cite[Lemma~3.3]{LP}) 
that $HH_*^{I\!I}(C)$ and $HH^{I\!I\;*}(C)$ decompose into
direct sums over the critical values of~$w$.
 More precisely, let $X_{c_i}\subset X$ denote the closed subscheme
defined by the equation $w=c_i$, and let $X_i'\subset X$ denote
the open complement to the union of $X_{c_j}$ over all $j\ne i$.
 Clearly, $X_i'$ is an affine scheme.

 Let $B_i'$ denote the CDG\+algebra associated with $X'_i$ and~$w$
as above, and let $C_i'=\modrcfp B_i'$ be the related DG\+category
of CDG\+modules.
 Then the natural map $HH^{I\!I}_*(C)\rarrow\bigoplus_i
HH^{I\!I}_*(C_i')$ induced by the DG\+functors $C\rarrow C_i'$ is
an isomorphism, since the natural map $HH^{I\!I}_*(B)\rarrow
\bigoplus_i HH^{I\!I}_*(B_i')$ is.

 Furthermore, the diagonal CDG\+bimodule $B_i'$ over $B_i'$ can be
considered as a CDG\+bimodule over $B$ by means of the strict
CDG\+functor $B\rarrow B_i'$, and similarly the diagonal
DG\+bimodule $C'_i$ over $C'_i$ can be considered as
a DG\+bimodule over~$C$.
 The CDG\+bimodule $B_i'$ over $B$ corresponds to the DG\+bimodule
$C_i'$ over $C$ under the equivalence of DG\+categories
$B\ot_k B^\op\modlc\simeq C\ot_k C^\op\modld$
from~\ref{dg-of-cdg-subsect}.
 The natural maps $HH^{I\!I\;*}(C_i')\rarrow HH^{I\!I\;*}(C,C_i')$
are isomorphisms, since the natural maps $HH^{I\!I\;*}(B_i')\rarrow
HH^{I\!I\;*}(B,B_i')$ are; and the map $HH^{I\!I\;*}(C)\rarrow
\bigoplus_i HH^{I\!I\;*}(C,C_i')$ is an isomorphism, for the similar
reason.

 On the other hand, Orlov's theorem~\cite{Or} implies that
the DG\+functor $C_{(c_i)}\rarrow C'_{i\.(c_i)}$, where $C'_{i\.(c_i)}
= \modrcfp B'_{i\.(c_i)}$ and $B'_{i\.(c_i)}$ is the CDG\+algebra
associated with the variety $X'_i$ and the function $w-c_i$,
is a quasi-equivalence.
 Thererefore, the natural maps $HH_*(C_{(c_i)})\rarrow
HH_*(C'_{i\.(c_i)})$ are isomorphisms, as are the natural maps
$HH^*(C_{(c_i)})\rarrow HH^*(C_{(c_i)},C'_{i\.(c_i)})\longleftarrow
HH^*(C'_{i\.(c_i)})$.
 Besides, one has $H^0(C_{i\.(c_j)})=0$, hence
$HH_*((C_{i\.(c_j)}) = 0 = HH^*((C_{i\.(c_j)})$, for all $i\ne j$.

 The
isomorphisms~(\ref{ho-hoch-direct-sum}\+-\ref{coho-hoch-direct-sum})
now follow from Corollary~\ref{matrix-factorizations}.A applied to
the varieties $X_i'$ with the functions $w-c_i$ on them and
commutativity of the natural diagrams.
\end{proof}

\bigskip


\begin{thebibliography}{99}
\medskip

\bibitem{CT}
 A.~Caldararu, J.~Tu.
   Curved A-infinity algebras and Landau-Ginzburg models.
Electronic preprint \texttt{arXiv:1007.2679 [math.KT]}, 42~pp.

\bibitem{Dyck}
 T.~Dyckerhoff.
   Compact generators in categories of matrix factorizations.
Electronic preprint \texttt{arXiv:0904.4713 [math.AG]}, 43~pp.

\bibitem{GJ}
 E.~Getzler, J.~D.~S.~Jones.
   $A_\infty$\+algebras and the cyclic bar complex.
\textit{Illinois Journ.\ of Math.}\ \textbf{34}, \#2, 1990.

\bibitem{EGAIV4}
 A.~Grothendieck, J.~Dieudonn\'e.
   \'El\'ements de g\'eom\'etrie alg\'ebrique IV.
\'Etude locale des sch\'emas et des morphismes de sch\'emas,
Quatri\`eme partie.
\textit{Publ.\ Math.\ IHES} \textbf{32}, p.~5--361, 1967.

\bibitem{HMS}
 D.~Husemoller, J.~C.~Moore, J.~Stasheff.
   Differential homological algebra and homogeneous spaces.
\textit{Journ.\ Pure Appl. Algebra} \textbf{5}, p.~113--185, 1974.

\bibitem{Kel}
 B.~Keller.
   Deriving DG-categories.
\textit{Ann.\ Sci.\ \'Ecole Norm.\ Sup.\ (4)} \textbf{27}, \#1,
p.~63--102, 1994.

\bibitem{KLN}
 B.~Keller, W.~Lowen, P.~Nicol\'as.
   On the (non)vanishing of some ``derived'' categories of curved
dg algebras.
\textit{Journ.\ Pure Appl.\ Algebra} \textbf{214}, \#7, p.~1271--1284,
2010.  \texttt{arXiv:0905.3845 [math.KT]}

\bibitem{LP}
 K.~Lin, D.~Pomerleano.
   Global matrix factorizations.
Electronic preprint \texttt{arXiv:1101.5847 [math.AG]}, 15~pp.

\bibitem{Neem}
 A.~Neeman.
     The Grothendieck duality theorem via Bousfield's techniques and
Brown representability.
\textit{Journ.\ Amer.\ Math.\ Soc.}\ \textbf{9}, p.~205--236, 1996.

\bibitem{Or}
 D.~Orlov.
     Triangulated categories of singularities and D-branes in
Landau--Ginzburg models.
\textit{Proc.\ Steklov Math.\ Inst.}\ \textbf{246}, \#3,
p.~227--248, 2004.  \texttt{arXiv:math.AG/0302304}

\bibitem{Or2}
 D.~Orlov.
   Formal completions and idempotent completions of triangulated
categories of singularities.
\textit{Advances in Math.}\ \textbf{226}, \#1, p.~206--217, 2011.
\texttt{arXiv:0901.1859 [math.AG]}

\bibitem{PV}
 A.~Polishchuk, A.~Vaintrob.
   Chern characters and Hirzebruch-Riemann-Roch formula for matrix
factorizations.
 Electronic preprint \texttt{arXiv:1002.2116 [math.AG]}, 45~pp.

\bibitem{Pcurv}
 L.~Positselski.
   Nonhomogeneous quadratic duality and curvature.
\textit{Funct.\ Anal.\ Appl.}\ \textbf{27}, \#3, p.~197--204, 1993.

\bibitem{Psemi}
 L.~Positselski.
   Homological algebra of semimodules and semicontramodules:
Semi-infinite homological algebra of associative algebraic structures.
Appendix~C in collaboration with D.~Rumynin; Appendix~D in
collaboration with S.~Arkhipov.
Monografie Matematyczne, vol.~70, Birkh\"auser/Springer Basel, 2010,
xxiv+349~pp.  \texttt{arXiv:0708.3398 [math.CT]}

\bibitem{Pkoszul}
 L.~Positselski.
   Two kinds of derived categories, Koszul duality, and
comodule-contramodule correspondence.
\textit{Memoirs Amer.\ Math.\ Soc.} \textbf{212}, \#996,
2011, v+133~pp.  \texttt{arXiv:0905.2621 [math.CT]}

\bibitem{Porl}
 L.~Positselski.
   Coherent analogues of matrix factorizations and relative
singularity categories.
 Electronic preprint \texttt{arXiv:1102.0261 [math.CT]}, 16~pp.

\bibitem{PVi}
 L.~Positselski, A.~Vishik.
     Koszul duality and Galois cohomology.
\textit{Math.\ Research Letters} \textbf{2}, \#6, pp.~771--781,
1995.  \texttt{arXiv:alg-geom/9507010}

\bibitem{Seg}
 E.~Segal.
   The closed state space of affine Landau-Ginzburg B-models.
Electronic preprint \texttt{arXiv:0904.1339 [math.AG]}, 30~pp.

\bibitem{Sch}
 A.~Schwarz.
   Noncommutative supergeometry, duality and deformations.
\textit{Nuclear Physics B} \textbf{650}, p.~475--496, 2003.
\texttt{arXiv:hep-th/0210271}

\bibitem{Tab}
 G.~Tabuada.
   Une structure de cat\'egorie de mod\`eles de Quillen sur
la cat\'egorie des dg-cat\'egories.
   A Quillen model structure on the category of dg categories.
\textit{Comptes Rendus Acad.\ Sci.\ Paris} \textbf{340}, \#1,
p.~15--19, 2005.
\texttt{arXiv:math.KT/0407338}

\bibitem{Toen}
 B.~To\"en.
   The homotopy theory of dg-categories and derived Morita
theory.
\textit{Inventiones Math.} \textbf{167}, \#3, p.~615--667, 2007.
\texttt{arXiv:math.AG/0408337}

\bibitem{Tu}
 J.~Tu.
   Matrix factorizations via Koszul duality.
Electronic preprint \texttt{arXiv:1009.4151 [math.AG]}, 53~pp.
{\hbadness=1500\par}

\end{thebibliography}
\end{document}